\documentclass[12pt]{article}
\usepackage{amssymb}
\usepackage{amsmath}
\usepackage{amsthm}
\begin{document}
\def\Drb{\nobreak\hfill\nobreak\lower 1pt\hbox{$\square $}}
\def\black{\vrule height 6pt width 6pt depth 1pt} 
\newcommand{\Aut}{\mbox{Aut }}
\newcommand{\shift}{\mbox{shift}}
\newcommand{\Ad}{\mbox{Ad }}
\newcommand{\Fun}{\mbox{Fun }}
\newcommand{\ind}{\mbox{ind }}
\newcommand{\diag}{\mbox{diag }}
\newcommand{\Ind}{\mbox{Ind }}
\newcommand{\spann}{\mbox{span }}
\newcommand{\Ob}{\mbox{Ob }}
\newcommand{\Hom}{\mbox{Hom }}
\newcommand{\pfeil}{\longrightarrow}
\newcommand{\Pfeil}{\longmapsto}
\newcommand{\proj}{\mbox{proj }}
\renewcommand{\blacksquare}{\nobreak \hfill \nobreak \black}
\newcommand{\blacksquared}{\nobreak \hfill \nobreak \black}
\newcommand{\tm}{\otimes _M}
\newcommand{\tn}{\otimes _N} 
\newcommand{\tp}{\otimes _P}
\newcommand{\tq}{\otimes _Q}
\newcommand{\ta}{\otimes _A}
\renewcommand{\th}{\otimes _H}
\newcommand{\kreuz}{\rtimes}
\newcommand{\ckreuz}{{\rtimes} _c}
\newcommand{\ockreuz}{{\rtimes} _{\bar{c}}}
\newcommand{\rest}{\! \mid \!}
\newcommand{\C}[1]{{\cal #1}}
\newcommand{\B}[1]{{\bf #1}}
\newcommand{\equi}{\Longleftrightarrow}
\newcommand{\beq}{\begin{equation}}
\newcommand{\eeq}{\end{equation}}
\newcommand{\Kl}[3]{\langle #1 \otimes #2, #3\rangle}
\newcommand{\ska}[2]{\langle #1,\, #2\rangle}
\newcommand{\skab}[2]{\Bigl\langle #1,\, #2\Bigr\rangle}
\newcommand{\skal}[2]{\langle #1,\, #2\rangle _N^l}
\newcommand{\skar}[2]{\langle #1,\, #2\rangle _N^r}
\newcommand{\ov}{\overline}
\newcommand{\newl}{\newline \noindent}
\newcommand{\tr}{\mbox{tr}}
\newcommand{\Tr}{\mbox{Tr}}
\newcommand{\zws}{{\rm II}_1}
\newcommand{\abs}{\vspace{4mm}}
\newcommand{\orho}{\bar{\rho }}
\newcommand{\osigma}{\bar{\sigma }}
\newcommand{\ophi}{\bar{\phi }}
\newcommand{\opsi}{\bar{\psi }}
\newcommand{\otau}{\bar{\tau }}
\newcommand{\ovs}[1]{\ov{\sigma (#1)}}
\newcommand{\oR}{\bar{R}}
\newcommand{\oV}{\bar{V}}
\newcommand{\oW}{\bar{W}}
\newcommand{\oc}{\bar{c}}
\newcommand{\comp}{\mathbb{C}}
\newcommand{\nat}{\mathbb{N}}
\newcommand{\reel}{\mathbb{R}}
\newcommand{\ganz}{\mathbb{Z}}
\newcommand{\rat}{\mathbb{Q}}

 \title{$C^*$-Tensor Categories in the Theory of $\zws $-Subfactors}
\author{Reinhard Schaf\-litzel\\
        Mathematisches Institut\\
        Technische Universit\"at M\"unchen\\
        Arcisstr. 21\\ 
        80290 M\"unchen\\
        GERMANY\\
        E-mail: schafl@mathematik.tu-muenchen.de}
\date {December 2, 1996}

\maketitle
\theoremstyle{plain}
\newtheorem{lem}[subsection]{Lemma}
\newtheorem{cor}[subsection]{Corollary}
\newtheorem{prop}[subsection]{Proposition}
\newtheorem{thm}[subsection]{Theorem}
\theoremstyle{definition}
\newtheorem{defi}[subsection]{Definition}
\newtheorem{quest}[subsection]{Question}
\newtheorem{ass}[subsection]{Assumption}
\setcounter{tocdepth}{1}
\hyphenation{bundles Yamagami}
%
%
%
%
 
\long\def\ig#1{\relax}
\ig{Thanks to Roberto Minio for this def'n.  Compare the def'n of
\comment in AMSTeX.}
 
\setlength{\unitlength}{.01em}%
 
\newcount \coefa
\newcount \coefb
\newcount \coefc
\newdimen\tempdimen
\newdimen\xlen
\newdimen\ylen
 
\newcount\tempcounta
\newcount\tempcountb
\newcount\tempcountc
\newcount\tempcountd
\newcount\tempcounte
\newcount\tempcountf
\newcount\xext
\newcount\yext
\newcount\xoff
\newcount\yoff
\newcount\gap%
\newcount\arrowtypea
\newcount\arrowtypeb
\newcount\arrowtypec
\newcount\arrowtyped
\newcount\arrowtypee
\newcount\height
\newcount\width
\newcount\xpos
\newcount\ypos
\newcount\run
\newcount\rise
\newcount\arrowlength
\newcount\halflength
\newcount\arrowtype
\newsavebox{\tempboxa}%
\newsavebox{\tempboxb}%
\newsavebox{\tempboxc}%

\catcode`@=11 
\def\settoheight#1#2{\setbox\@tempboxa\hbox{#2}#1\ht\@tempboxa\relax}%
\def\settodepth#1#2{\setbox\@tempboxa\hbox{#2}#1\dp\@tempboxa\relax}%
\let\ifnextchar=\@ifnextchar
\catcode`@=12 
 
\def\putbox(#1,#2)#3{\put(#1,#2){\makebox(0,0){#3}}}

\def\setsqparms[#1`#2`#3`#4;#5`#6]{%
\settripairparms[#1`#2`#3`#4`1;#6]%
\width #5
}
 
\def\settriparms[#1`#2`#3;#4]{\settripairparms[#1`#2`#3`1`1;#4]}%

\def\settripairparms[#1`#2`#3`#4`#5;#6]{%
\arrowtypea #1
\arrowtypeb #2
\arrowtypec #3
\arrowtyped #4
\arrowtypee #5
\height #6
\width #6
}
 
\def\resetparms{\settripairparms[1`1`1`1`1;500]\width 500}
 
\def\mvector(#1,#2)#3{
\put(0,0){\vector(#1,#2){#3}}%
\put(0,0){\vector(#1,#2){30}}%
}
\def\evector(#1,#2)#3{{
\arrowlength #3
\put(0,0){\vector(#1,#2){\arrowlength}}%
\advance \arrowlength by-30
\put(0,0){\vector(#1,#2){\arrowlength}}%
}}

\def\horsize#1#2{%
\settowidth{\tempdimen}{$#2$}%
#1=\tempdimen
\divide #1 by\unitlength
}
 
\def\vertsize#1#2{%
\settoheight{\tempdimen}{$#2$}%
#1=\tempdimen
\settodepth{\tempdimen}{$#2$}%
\advance #1 by\tempdimen
\divide #1 by\unitlength
}

\def\vertadjust[#1`#2`#3]{%
\vertsize{\tempcounta}{#1}%
\vertsize{\tempcountb}{#2}%
\ifnum \tempcounta<\tempcountb \tempcounta=\tempcountb \fi
\divide\tempcounta by2
\vertsize{\tempcountb}{#3}%
\ifnum \tempcountb>0 \advance \tempcountb by20 \fi
\ifnum \tempcounta<\tempcountb \tempcounta=\tempcountb \fi
}
 
\def\horadjust[#1`#2`#3]{%
\horsize{\tempcounta}{#1}%
\horsize{\tempcountb}{#2}%
\ifnum \tempcounta<\tempcountb \tempcounta=\tempcountb \fi
\divide\tempcounta by2
\horsize{\tempcountb}{#3}%
\ifnum \tempcountb>0 \advance \tempcountb by60 \fi
\ifnum \tempcounta<\tempcountb \tempcounta=\tempcountb \fi
}
 
\ig{ In this procedure, #1 is the paramater that sticks out all the way,
#2 sticks out the least and #3 is a label sticking out half way.  #4 is
the amount of the offset.}
 
\def\sladjust[#1`#2`#3]#4{%
\tempcountc=#4
\horsize{\tempcounta}{#1}%
\divide \tempcounta by2
\horsize{\tempcountb}{#2}%
\divide \tempcountb by2
\advance \tempcountb by-\tempcountc
\ifnum \tempcounta<\tempcountb \tempcounta=\tempcountb\fi
\divide \tempcountc by2
\horsize{\tempcountb}{#3}%
\advance \tempcountb by-\tempcountc
\ifnum \tempcountb>0 \advance \tempcountb by80\fi
\ifnum \tempcounta<\tempcountb \tempcounta=\tempcountb\fi
\advance\tempcounta by20
}
 
\def\putvector(#1,#2)(#3,#4)#5#6{{%
\xpos=#1
\ypos=#2
\run=#3
\rise=#4
\arrowlength=#5
\arrowtype=#6
\ifnum \arrowtype<0
    \ifnum \run=0
        \advance \ypos by-\arrowlength
    \else
        \tempcounta \arrowlength
        \multiply \tempcounta by\rise
        \divide \tempcounta by\run
        \ifnum\run>0
            \advance \xpos by\arrowlength
            \advance \ypos by\tempcounta
        \else
            \advance \xpos by-\arrowlength
            \advance \ypos by-\tempcounta
        \fi
    \fi
    \multiply \arrowtype by-1
    \multiply \rise by-1
    \multiply \run by-1
\fi
\ifnum \arrowtype=1
    \put(\xpos,\ypos){\vector(\run,\rise){\arrowlength}}%
\else\ifnum \arrowtype=2
    \put(\xpos,\ypos){\mvector(\run,\rise)\arrowlength}%
\else\ifnum\arrowtype=3
    \put(\xpos,\ypos){\evector(\run,\rise){\arrowlength}}%
\fi\fi\fi
}}
 
\def\bfig{\begin{picture}(\xext,\yext)(\xoff,\yoff)}
\def\efig{\end{picture}}

\def\putsplitvector(#1,#2)#3#4{
\xpos #1
\ypos #2
\arrowtype #4
\halflength #3
\arrowlength #3
\gap 140
\advance \halflength by-\gap
\divide \halflength by2
\ifnum \arrowtype=1
    \put(\xpos,\ypos){\line(0,-1){\halflength}}%
    \advance\ypos by-\halflength
    \advance\ypos by-\gap
    \put(\xpos,\ypos){\vector(0,-1){\halflength}}%
\else\ifnum \arrowtype=2
    \put(\xpos,\ypos){\line(0,-1)\halflength}%
    \put(\xpos,\ypos){\vector(0,-1)3}%
    \advance\ypos by-\halflength
    \advance\ypos by-\gap
    \put(\xpos,\ypos){\vector(0,-1){\halflength}}%
\else\ifnum\arrowtype=3
    \put(\xpos,\ypos){\line(0,-1)\halflength}%
    \advance\ypos by-\halflength
    \advance\ypos by-\gap
    \put(\xpos,\ypos){\evector(0,-1){\halflength}}%
\else\ifnum \arrowtype=-1
    \advance \ypos by-\arrowlength
    \put(\xpos,\ypos){\line(0,1){\halflength}}%
    \advance\ypos by\halflength
    \advance\ypos by\gap
    \put(\xpos,\ypos){\vector(0,1){\halflength}}%
\else\ifnum \arrowtype=-2
    \advance \ypos by-\arrowlength
    \put(\xpos,\ypos){\line(0,1)\halflength}%
    \put(\xpos,\ypos){\vector(0,1)3}%
    \advance\ypos by\halflength
    \advance\ypos by\gap
    \put(\xpos,\ypos){\vector(0,1){\halflength}}%
\else\ifnum\arrowtype=-3
    \advance \ypos by-\arrowlength
    \put(\xpos,\ypos){\line(0,1)\halflength}%
    \advance\ypos by\halflength
    \advance\ypos by\gap
    \put(\xpos,\ypos){\evector(0,1){\halflength}}%
\fi\fi\fi\fi\fi\fi
}
 
\def\setpos(#1,#2){\xpos=#1 \ypos#2}
 
\def\putmorphism(#1)(#2,#3)[#4`#5`#6]#7#8#9{{%
\run #2
\rise #3
\ifnum\rise=0
  \puthmorphism(#1)[#4`#5`#6]{#7}{#8}{#9}%
\else\ifnum\run=0
  \putvmorphism(#1)[#4`#5`#6]{#7}{#8}{#9}%
\else
\setpos(#1)%
\arrowlength #7
\arrowtype #8
\ifnum\run=0
\else\ifnum\rise=0
\else
\ifnum\run>0
    \coefa=1
\else
   \coefa=-1
\fi
\ifnum\arrowtype>0
   \coefb=0
   \coefc=-1
\else
   \coefb=\coefa
   \coefc=1
   \arrowtype=-\arrowtype
\fi
\width=2
\multiply \width by\run
\divide \width by\rise
\ifnum \width<0  \width=-\width\fi
\advance\width by60
\if l#9 \width=-\width\fi
\putbox(\xpos,\ypos){$#4$}
{\multiply \coefa by\arrowlength
\advance\xpos by\coefa
\multiply \coefa by\rise
\divide \coefa by\run
\advance \ypos by\coefa
\putbox(\xpos,\ypos){$#5$} }%
{\multiply \coefa by\arrowlength
\divide \coefa by2
\advance \xpos by\coefa
\advance \xpos by\width
\multiply \coefa by\rise
\divide \coefa by\run
\advance \ypos by\coefa
\if l#9%
   \put(\xpos,\ypos){\makebox(0,0)[r]{$#6$}}%
\else\if r#9%
   \put(\xpos,\ypos){\makebox(0,0)[l]{$#6$}}%
\fi\fi }%
{\multiply \rise by-\coefc
\multiply \run by-\coefc
\multiply \coefb by\arrowlength
\advance \xpos by\coefb
\multiply \coefb by\rise
\divide \coefb by\run
\advance \ypos by\coefb
\multiply \coefc by70
\advance \ypos by\coefc
\multiply \coefc by\run
\divide \coefc by\rise
\advance \xpos by\coefc
\multiply \coefa by140
\multiply \coefa by\run
\divide \coefa by\rise
\advance \arrowlength by\coefa
\ifnum \arrowtype=1
   \put(\xpos,\ypos){\vector(\run,\rise){\arrowlength}}%
\else\ifnum\arrowtype=2
   \put(\xpos,\ypos){\mvector(\run,\rise){\arrowlength}}%
\else\ifnum\arrowtype=3
   \put(\xpos,\ypos){\evector(\run,\rise){\arrowlength}}%
\fi\fi\fi}%
\fi\fi
\fi\fi}}

\def\puthmorphism(#1,#2)[#3`#4`#5]#6#7#8{{%
\xpos #1
\ypos #2
\width #6
\arrowlength #6
\putbox(\xpos,\ypos){$#3$\vphantom{$#4$}}%
{\advance \xpos by\arrowlength
\putbox(\xpos,\ypos){\vphantom{$#3$}$#4$}}%
\horsize{\tempcounta}{#3}%
\horsize{\tempcountb}{#4}%
\divide \tempcounta by2
\divide \tempcountb by2
\advance \tempcounta by30
\advance \tempcountb by30
\advance \xpos by\tempcounta
\advance \arrowlength by-\tempcounta
\advance \arrowlength by-\tempcountb
\putvector(\xpos,\ypos)(1,0){\arrowlength}{#7}%
\divide \arrowlength by2
\advance \xpos by\arrowlength
\vertsize{\tempcounta}{#5}%
\divide\tempcounta by2
\advance \tempcounta by20
\if a#8 %
   \advance \ypos by\tempcounta
   \put(\xpos,\ypos){\makebox(0,0){$#5$}}%
\else
   \advance \ypos by-\tempcounta
   \put(\xpos,\ypos){\makebox(0,0){$#5$}}%
\fi
}}
 
\def\putvmorphism(#1,#2)[#3`#4`#5]#6#7#8{{%
\xpos #1
\ypos #2
\arrowlength #6
\arrowtype #7
\settowidth{\xlen}{$#5$}%
\putbox(\xpos,\ypos){$#3$}%
{\advance \ypos by-\arrowlength
\putbox(\xpos,\ypos){$#4$}}%
{\advance\arrowlength by-140
\advance \ypos by -70
\ifdim\xlen>0pt
   \if m#8%
      \putsplitvector(\xpos,\ypos){\arrowlength}{\arrowtype}%
   \else
      \putvector(\xpos,\ypos)(0,-1){\arrowlength}{\arrowtype}%
   \fi
\else
   \putvector(\xpos,\ypos)(0,-1){\arrowlength}{\arrowtype}%
\fi}%
\ifdim\xlen>0pt
   \divide \arrowlength by2
   \advance\ypos by-\arrowlength
   \if l#8%
      \advance \xpos by-40
      \put(\xpos,\ypos){\makebox(0,0)[r]{$#5$}}%
   \else\if r#8%
      \advance \xpos by40
      \put(\xpos,\ypos){\makebox(0,0)[l]{$#5$}}%
   \else
      \putbox(\xpos,\ypos){$#5$}%
   \fi\fi
\fi
}}
 
\def\topadjust[#1`#2`#3]{%
\yoff=10
\vertadjust[#1`#2`{#3}]%
\advance \yext by\tempcounta
\advance \yext by 10
}
 
\def\botadjust[#1`#2`#3]{%
\vertadjust[#1`#2`{#3}]%
\advance \yext by\tempcounta
\advance \yoff by-\tempcounta
}
 
\def\leftadjust[#1`#2`#3]{%
\xoff=0
\horadjust[#1`#2`{#3}]%
\advance \xext by\tempcounta
\advance \xoff by-\tempcounta
}
 
\def\rightadjust[#1`#2`#3]{%
\horadjust[#1`#2`{#3}]%
\advance \xext by\tempcounta
}
 
\def\rightsladjust[#1`#2`#3]{%
\sladjust[#1`#2`{#3}]{\width}%
\advance \xext by\tempcounta
}
 
\def\leftsladjust[#1`#2`#3]{%
\xoff=0
\sladjust[#1`#2`{#3}]{\width}%
\advance \xext by\tempcounta
\advance \xoff by-\tempcounta
}
 
\def\adjust[#1`#2;#3`#4;#5`#6;#7`#8]{%
\topadjust[#1``{#2}]
\leftadjust[#3``{#4}]
\rightadjust[#5``{#6}]
\botadjust[#7``{#8}]}

\def\putsquare(#1)[#2`#3`#4`#5;#6`#7`#8`#9]{%
\setpos(#1)
\puthmorphism(\xpos,\ypos)[#4`#5`{#9}]{\width}{\arrowtyped}b%
\advance\ypos by \height
\puthmorphism(\xpos,\ypos)[#2`#3`{#6}]{\width}{\arrowtypea}a%
\putvmorphism(\xpos,\ypos)[``{#7}]{\height}{\arrowtypeb}l%
\advance\xpos by \width
\putvmorphism(\xpos,\ypos)[``{#8}]{\height}{\arrowtypec}r%
}
 
\def\square[#1`#2`#3`#4;#5`#6`#7`#8]{{
\xext=\width                              
\yext=\height                             
\topadjust[#1`#2`{#5}]
\botadjust[#3`#4`{#8}]
\leftadjust[#1`#3`{#6}]
\rightadjust[#2`#4`{#7}]
\begin{picture}(\xext,\yext)(\xoff,\yoff)
\putsquare(0,0)[#1`#2`#3`#4;#5`#6`#7`{#8}]
\end{picture}
}}

\def\putptriangle(#1,#2)[#3`#4`#5;#6`#7`#8]{%
\xpos=#1 \ypos=#2
\advance\ypos by \height
\puthmorphism(\xpos,\ypos)[#3`#4`{#6}]{\height}{\arrowtypea}a%
\putvmorphism(\xpos,\ypos)[`#5`{#7}]{\height}{\arrowtypeb}l%
\advance\xpos by\height
\putmorphism(\xpos,\ypos)(-1,-1)[``{#8}]{\height}{\arrowtypec}r%
}
 
\def\ptriangle[#1`#2`#3;#4`#5`#6]{{
\width=\height                         
\xext=\width                           
\yext=\width                           
\topadjust[#1`#2`{#4}]
\botadjust[#3``]
\leftadjust[#1`#3`{#5}]
\rightsladjust[#2`#3`{#6}]
\begin{picture}(\xext,\yext)(\xoff,\yoff)
\putptriangle(0,0)[#1`#2`#3;#4`#5`{#6}]%
\end{picture}%
}}

\def\putqtriangle(#1,#2)[#3`#4`#5;#6`#7`#8]{%
\xpos=#1 \ypos=#2
\advance\ypos by\height
\puthmorphism(\xpos,\ypos)[#3`#4`{#6}]{\height}{\arrowtypea}a%
\putmorphism(\xpos,\ypos)(1,-1)[``{#7}]{\height}{\arrowtypeb}l%
\advance\xpos by\height
\putvmorphism(\xpos,\ypos)[`#5`{#8}]{\height}{\arrowtypec}r%
}
 
\def\qtriangle[#1`#2`#3;#4`#5`#6]{{
\width=\height                         
\xext=\width                           
\yext=\height                          
\topadjust[#1`#2`{#4}]
\botadjust[#3``]
\leftsladjust[#1`#3`{#5}]
\rightadjust[#2`#3`{#6}]
\begin{picture}(\xext,\yext)(\xoff,\yoff)
\putqtriangle(0,0)[#1`#2`#3;#4`#5`{#6}]%
\end{picture}%
}}

\def\putdtriangle(#1,#2)[#3`#4`#5;#6`#7`#8]{%
\xpos=#1 \ypos=#2
\puthmorphism(\xpos,\ypos)[#4`#5`{#8}]{\height}{\arrowtypec}b%
\advance\xpos by \height \advance\ypos by\height
\putmorphism(\xpos,\ypos)(-1,-1)[``{#6}]{\height}{\arrowtypea}l%
\putvmorphism(\xpos,\ypos)[#3``{#7}]{\height}{\arrowtypeb}r%
}
 
\def\dtriangle[#1`#2`#3;#4`#5`#6]{{
\width=\height                         
\xext=\width                           
\yext=\height                          
\topadjust[#1``]
\botadjust[#2`#3`{#6}]
\leftsladjust[#2`#1`{#4}]
\rightadjust[#1`#3`{#5}]
\begin{picture}(\xext,\yext)(\xoff,\yoff)
\putdtriangle(0,0)[#1`#2`#3;#4`#5`{#6}]%
\end{picture}%
}}

\def\putbtriangle(#1,#2)[#3`#4`#5;#6`#7`#8]{%
\xpos=#1 \ypos=#2
\puthmorphism(\xpos,\ypos)[#4`#5`{#8}]{\height}{\arrowtypec}b%
\advance\ypos by\height
\putmorphism(\xpos,\ypos)(1,-1)[``{#7}]{\height}{\arrowtypeb}r%
\putvmorphism(\xpos,\ypos)[#3``{#6}]{\height}{\arrowtypea}l%
}
 
\def\btriangle[#1`#2`#3;#4`#5`#6]{{
\width=\height                         
\xext=\width                           
\yext=\height                          
\topadjust[#1``]
\botadjust[#2`#3`{#6}]
\leftadjust[#1`#2`{#4}]
\rightsladjust[#3`#1`{#5}]
\begin{picture}(\xext,\yext)(\xoff,\yoff)
\putbtriangle(0,0)[#1`#2`#3;#4`#5`{#6}]%
\end{picture}%
}}

\def\putAtriangle(#1,#2)[#3`#4`#5;#6`#7`#8]{%
\xpos=#1 \ypos=#2
{\multiply \height by2
\puthmorphism(\xpos,\ypos)[#4`#5`{#8}]{\height}{\arrowtypec}b}%
\advance\xpos by\height \advance\ypos by\height
\putmorphism(\xpos,\ypos)(-1,-1)[#3``{#6}]{\height}{\arrowtypea}l%
\putmorphism(\xpos,\ypos)(1,-1)[``{#7}]{\height}{\arrowtypeb}r%
}
 
\def\Atriangle[#1`#2`#3;#4`#5`#6]{{
\width=\height                         
\xext=\width                           
\yext=\height                          
\topadjust[#1``]
\botadjust[#2`#3`{#6}]
\multiply \xext by2 
\leftsladjust[#2`#1`{#4}]
\rightsladjust[#3`#1`{#5}]
\begin{picture}(\xext,\yext)(\xoff,\yoff)%
\putAtriangle(0,0)[#1`#2`#3;#4`#5`{#6}]%
\end{picture}%
}}

\def\putAtrianglepair(#1,#2)[#3]{\xpos=#1 \ypos=#2%
\putAtrianglepairx[#3]}
\def\putAtrianglepairx[#1`#2`#3`#4;#5`#6`#7`#8`#9]{%
\puthmorphism(\xpos,\ypos)[#2`#3`{#8}]{\height}{\arrowtyped}b%
\advance\xpos by\height
\puthmorphism(\xpos,\ypos)[\phantom{#3}`#4`{#9}]{\height}{\arrowtypee}b%
\advance\ypos by\height
\putmorphism(\xpos,\ypos)(-1,-1)[#1``{#5}]{\height}{\arrowtypea}l%
\putvmorphism(\xpos,\ypos)[``{#6}]{\height}{\arrowtypeb}m%
\putmorphism(\xpos,\ypos)(1,-1)[``{#7}]{\height}{\arrowtypec}r%
}
 
\def\Atrianglepair[#1`#2`#3`#4;#5`#6`#7`#8`#9]{{%
\width=\height
\xext=\width
\yext=\height
\topadjust[#1``]%
\vertadjust[#2`#3`{#8}]
\tempcountd=\tempcounta                       
\vertadjust[#3`#4`{#9}]
\ifnum\tempcounta<\tempcountd                 
\tempcounta=\tempcountd\fi                    
\advance \yext by\tempcounta                  
\advance \yoff by-\tempcounta                 
\multiply \xext by2 
\leftsladjust[#2`#1`{#5}]
\rightsladjust[#4`#1`{#7}]%
\begin{picture}(\xext,\yext)(\xoff,\yoff)%
\putAtrianglepair(0,0)[#1`#2`#3`#4;#5`#6`#7`#8`{#9}]%
\end{picture}%
}}

\def\putVtriangle(#1,#2)[#3`#4`#5;#6`#7`#8]{%
\xpos=#1 \ypos=#2
\advance\ypos by\height
{\multiply\height by2
\puthmorphism(\xpos,\ypos)[#3`#4`{#6}]{\height}{\arrowtypea}a}%
\putmorphism(\xpos,\ypos)(1,-1)[`#5`{#7}]{\height}{\arrowtypeb}l%
\advance\xpos by\height
\advance\xpos by\height
\putmorphism(\xpos,\ypos)(-1,-1)[``{#8}]{\height}{\arrowtypec}r%
}
 
\def\Vtriangle[#1`#2`#3;#4`#5`#6]{{
\width=\height                         
\xext=\width                           
\yext=\height                          
\topadjust[#1`#2`{#4}]
\botadjust[#3``]
\multiply \xext by2 
\leftsladjust[#1`#3`{#5}]
\rightsladjust[#2`#3`{#6}]
\begin{picture}(\xext,\yext)(\xoff,\yoff)%
\putVtriangle(0,0)[#1`#2`#3;#4`#5`{#6}]%
\end{picture}%
}}

\def\putVtrianglepair(#1,#2)[#3]{\xpos=#1 \ypos=#2%
\putVtrianglepairx[#3]}
\def\putVtrianglepairx[#1`#2`#3`#4;#5`#6`#7`#8`#9]{%
\advance\ypos by\height
\putmorphism(\xpos,\ypos)(1,-1)[`#4`{#7}]{\height}{\arrowtypec}l%
\puthmorphism(\xpos,\ypos)[#1`#2`{#5}]{\height}{\arrowtypea}a%
\advance\xpos by\height
\puthmorphism(\xpos,\ypos)[\phantom{#2}`#3`{#6}]{\height}{\arrowtypeb}a%
\putvmorphism(\xpos,\ypos)[``{#8}]{\height}{\arrowtyped}m%
\advance\xpos by\height
\putmorphism(\xpos,\ypos)(-1,-1)[``{#9}]{\height}{\arrowtypee}r%
}
 
\def\Vtrianglepair[#1`#2`#3`#4;#5`#6`#7`#8`#9]{{%
\xoff=0
\yoff=2 
\xext=\height                  
\width=\height                 
\yext=\height                  
\vertadjust[#1`#2`{#5}]
\tempcountd=\tempcounta        
\vertadjust[#2`#3`{#6}]
\ifnum\tempcounta<\tempcountd
\tempcounta=\tempcountd\fi
\advance \yext by\tempcounta
\botadjust[#4``]%
\multiply \xext by2
\leftsladjust[#1`#4`{#7}]%
\rightsladjust[#3`#4`{#9}]%
\begin{picture}(\xext,\yext)(\xoff,\yoff)%
\putVtrianglepair(0,0)[#1`#2`#3`#4;#5`#6`#7`#8`{#9}]%
\end{picture}%
}}

\def\putCtriangle(#1,#2)[#3`#4`#5;#6`#7`#8]{%
\xpos=#1 \ypos=#2
\advance\ypos by\height
\putmorphism(\xpos,\ypos)(1,-1)[``{#8}]{\height}{\arrowtypec}l%
\advance\xpos by\height
\advance\ypos by\height
\putmorphism(\xpos,\ypos)(-1,-1)[#3`#4`{#6}]{\height}{\arrowtypea}l%
{\multiply\height by 2
\putvmorphism(\xpos,\ypos)[`#5`{#7}]{\height}{\arrowtypeb}r}%
}
 
\def\Ctriangle[#1`#2`#3;#4`#5`#6]{{
\width=\height                          
\xext=\width                            
\yext=\height                           
\multiply \yext by2 
\topadjust[#1``]
\botadjust[#3``]
\sladjust[#2`#1`{#4}]{\width}
\tempcountd=\tempcounta                 
\sladjust[#2`#3`{#6}]{\width}
\ifnum \tempcounta<\tempcountd          
\tempcounta=\tempcountd\fi              
\advance \xext by\tempcounta            
\advance \xoff by-\tempcounta           
\rightadjust[#1`#3`{#5}]
\begin{picture}(\xext,\yext)(\xoff,\yoff)%
\putCtriangle(0,0)[#1`#2`#3;#4`#5`{#6}]%
\end{picture}%
}}

\def\putDtriangle(#1,#2)[#3`#4`#5;#6`#7`#8]{%
\xpos=#1 \ypos=#2
\advance\xpos by\height \advance\ypos by\height
\putmorphism(\xpos,\ypos)(-1,-1)[``{#8}]{\height}{\arrowtypec}r%
\advance\xpos by-\height \advance\ypos by\height
\putmorphism(\xpos,\ypos)(1,-1)[`#4`{#7}]{\height}{\arrowtypeb}r%
{\multiply\height by 2
\putvmorphism(\xpos,\ypos)[#3`#5`{#6}]{\height}{\arrowtypea}l}%
}
 
\def\Dtriangle[#1`#2`#3;#4`#5`#6]{{
\width=\height                         
\xext=\height                          
\yext=\height                          
\multiply \yext by2 
\topadjust[#1``]
\botadjust[#3``]
\leftadjust[#1`#3`{#4}]
\sladjust[#2`#1`{#4}]{\height}
\tempcountd=\tempcountd                
\sladjust[#2`#3`{#6}]{\height}
\ifnum \tempcounta<\tempcountd         
\tempcounta=\tempcountd\fi             
\advance \xext by\tempcounta           
\begin{picture}(\xext,\yext)(\xoff,\yoff)
\putDtriangle(0,0)[#1`#2`#3;#4`#5`{#6}]%
\end{picture}%
}}

\def\setrecparms[#1`#2]{\width=#1 \height=#2}%
\def\recurse[#1`#2`#3`#4;#5`#6`#7`#8`#9]{{%
\settowidth{\tempdimen}{#1}
\ifdim\tempdimen=0pt
  \savebox{\tempboxa}{\hbox{#2}}%
  \savebox{\tempboxb}{\hbox{#4}}%
  \savebox{\tempboxc}{\hbox{#7}}%
\else
  \savebox{\tempboxa}{\hbox{$\hbox{#1}\times\hbox{#2}$}}%
  \savebox{\tempboxb}{\hbox{$\hbox{#1}\times\hbox{#4}$}}%
  \savebox{\tempboxc}{\hbox{$\hbox{#1}\times\hbox{#7}$}}%
\fi
\tempcounte=\height
\divide\tempcounte by 2
\tempcountf=\tempcounte
\advance\tempcountf by \width
\xext=\tempcountf \yext=\height
\topadjust[#2`\usebox{\tempboxa}`{#5}]%
\botadjust[#4`\usebox{\tempboxb}`{#9}]%
\sladjust[#3`#2`{#6}]{\tempcounte}%
\tempcountd=\tempcounta
\sladjust[#3`#4`{#8}]{\tempcounte}%
\ifnum \tempcounta<\tempcountd
\tempcounta=\tempcountd\fi
\advance \xext by\tempcounta
\advance \xoff by-\tempcounta
\rightadjust[\usebox{\tempboxa}`\usebox{\tempboxb}`\usebox{\tempboxc}]%
\bfig
{\settriparms[-1`1`1;\tempcounte]%
\putCtriangle(0,0)[`#3`;#6`#7`{#8}]}%
\arrowtypea=-1 \arrowtypeb=0 \arrowtypec=1 \arrowtyped=-1
\putsquare(\tempcounte,0)[#2`\usebox{\tempboxa}`#4`\usebox{\tempboxb};%
#5``\usebox{\tempboxc}`#9]%
\efig
}}

\begin{abstract} 
The article contains a detailed description of the connection between
finite depth inclusions of $\zws $-subfactors and 
finite $C^*$-tensor categories
(i.e. $C^*$-tensor categories
with dimension function for which the number of 
equivalence classes of irreducible objects is finite). The $(N,N)$-bimodules
belonging to a $\zws $-subfactor $N\subset M$ with finite Jones index form
a $C^*$-tensor category with dimension function. Conversely, taking an object
of a finite $C^*$-tensor category $\C{C}$ we construct a subfactor
$A\subset R$ of the hyperfinite $\zws$-factor $R$ with
finite index and finite depth. For this subfactor
we compute the standard invariant and show that the $C^*$-tensor category
of the corresponding $(A,A)$-bimodules is equivalent to a subcategory of
$\C{C}$. We illustrate the results for the $C^*$-tensor category of
the unitary finite dimensional corepresentations of a finite dimensional
Hopf-$*$-algebra.

\end{abstract}

\section*{Introduction} The theory of subfactors was established
by V.F.R. Jones in his famous paper \cite{Jo}. The goal of this paper is to
sketch the connection between $C^*$-tensor categories and subfactors of
type $\zws $-factors. We will use J.E. Roberts' theory of dimension for
$C^*$-tensor categories (see \cite{Ro0}, \cite{Ro}, and \cite{LoRo}), but
our approach to subfactors starting from $C^*$-tensor categories is different
from that in \cite{LoRo}, in contrast to \cite{LoRo} we restrict ourselves 
to factors of type $\zws $. Some parts of this article seem to be known to 
some experts (A. Wassermann for example), but the author could only find the
expositions of J.E. Roberts and R. Longo. This article is the abridged version
of a part of a paper (see \cite{Sch2}) accepted as a Habilitationsschrift
by the Technische Universit\"at M\"unchen.

$C^*$-tensor categories with dimension function (and some other properties)
are called compact $C^*$-tensor categories in this paper, one may regard them
as a concept to deal with more general symmetries than those described by
(compact) groups. Hopf algebras, in particular quantum groups, may also be 
regarded as concepts describing generalized symmetries. As one expects, 
there is a close connection between the concepts of $C^*$-tensor categories
and Hopf algebras.

Compact $C^*$-tensor categories are based on the following idea:
Consider the unitary
finite dimensional representations of a compact group $G$. For these 
representations    
one has intertwining operators. One may introduce subrepresentations,
finite direct sums and tensor products of representations. There is
a conjugate representation for every representation. One can
produce a category with an additional structure from these ingredients.
The objects are  representations, the morphisms are  intertwining
operators. The space of morphisms is endowed 
with an additional structure, which
is related to $C^*$-algebras. Moreover, the tensor product gives a product 
operation on the objects.
A further property of this category is the existence of conjugates. 
Compact $C^*$-tensor categories behave in a similar way as these categories.
While in the $C^*$-tensor category of a compact group the tensor product
of representations is commutative, 
the product need not be commutative for general compact $C^*$-tensor
categories. J.E. Roberts introduced a dimension for the objects of
compact $C^*$-tensor categories which in contrast to the group case need
not be a natural number in general.

S. Doplicher and J.E. Roberts showed the following (see \cite{DR}):
Every compact 
$C^*$-tensor category, for which the tensor product operation is commutative
in a certain strict way, is equivalent to the $C^*$-tensor category
of a suitable compact group.

In the theory of subfactors compact $C^*$-tensor categories appear in a
natural way.  A. Ocneanu had the idea to develop the theory of subfactors
by considering bimodules.
For a pair $N\subset M$ of $\zws $-factors with finite
index, he considered the $(N,N)$-, $(N,M)$-, $(M,N)$- and $(M,M)$-bimodules
contained in $L^2(M_k)$ for some $k\in \nat \cup \{0\}$; here $M_k$
denotes the $\zws $-factor after a $k+1$-fold application of the basic
construction. The occuring $(N,N)$-bimodules form a compact $C^*$-tensor
category, where the product operation for the bimodules is the
$N$-tensor product $\tn $.
\abs

An important subject of this paper is a method for the construction of 
subfactors
of the hyperfinite $\zws $-factor from a given $C^*$-tensor category.
Compact $C^*$-tensor categories, for which the number of the equivalence
classes of the irreducible objects is finite, are called finite
$C^*$-tensor categories. We construct a finite depth subfactor
of the hyperfinite $\zws $-factor for a given object of a finite
$C^*$-tensor category. The method is very general. Important special cases
of this construction are certain subfactors associated with finite 
groups as well as the subfactors due to H. Wenzl. In particular, an
explicit construction of the $C^*$-tensor categories associated
with H. Wenzl's subfactors is given.

Let us now give a more detailed outline of this paper:
\newl
In Section 1 $C^*$-tensor categories and especially compact $C^*$-tensor
categories are introduced and basic facts about these structures
are presented. 

Section 2 contains the proof that the $(N,N)$-bimodules 
associated with a subfactor $N\subset M$ of a $\zws $-factor $M$
with finite index form a compact $C^*$-tensor category.
Moreover
we study the subfactors $N\subset N\kreuz H$, where $N$ is a $\zws $-factor,
$H$ is a finite dimensional Hopf-$*$-algebra acting on $N$ by an outer
action and $N\kreuz H$ is the corresponding crossed product. We show that the
$C^*$-tensor category of the $(N,N)$-bimodules is equivalent to the 
$C^*$-tensor category of the unitary
corepresentations of $H^{cop}$, where $H^{cop}$ is the Hopf-$*$-algebra  
emerging from $H$ by reversing the comultiplication.

In Section 3  we deal with the construction of finite depth subfactors of the
hyperfinite $\zws $-factor from a given object of a finite
$C^*$-tensor category $\C{C}$. For that purpose we use H. Wenzl's 
technics which were 
developed in Chapter 1 of \cite{We} to investigate subfactors generated by a 
ladder of commuting squares. The standard invariant for these subfactors
is computed. This Section also contains a variant of this construction, which 
imitates Wenzl's construction of subfactors in \cite{We} and \cite{We2}.
Furthermore we
illustrate our construction in case $\C{C}$ is the $C^*$-tensor
category of the unitary representations of a finite group or more generally
the category of the unitary corepresentations of a finite dimensional
Hopf-$*$-algebra.

In Section 4 we show that the $C^*$-tensor category of the 
$(A,A)$-bimodules associated with the subfactor $A\subset B$ constructed in
Section 3 is equivalent to a subcategory of the $C^*$-tensor category
$\C{C}$, from which we started to construct the subfactor. The proof
involves a lot of computations in the category $\C{C}$ and in the
$\zws $-factors of the Jones tower for $A\subset B$.

%

In this paper we do not discuss concrete examples of compact 
$C^*$-tensor categories with objects having a non-integer dimension. We intend
to deal with $C^*$-tensor categories associated with H. Wenzl's subfactors
( \cite{We} and \cite{We2}) in a forthcoming paper. Another future project is
to generalize the constructions of this paper by using
$2$-$C^*$-tensor categories instead of $C^*$-tensor categories and including 
the $(N,M)$-, $(M,N)$- and $(M,M)$-bimodules.
\pagebreak
\section {$C^*$-tensor categories}
Mainly following J.E. Roberts, we introduce $C^*$-tensor categories;
in contrast to him we don't restrict ourselves to strict $C^*$-tensor
categories.
As in \cite{McL} we assume that the
objects of a category form a set. We fix a universe and call a set small
if it is an element of this universe (see \cite{McL}, Section I.6).
We suppose that all sets appearing in this paper are small. Henceforth
we will omit all details in this context.

Let $\C{C}$ be a category.
$\Ob \C{C}$ denotes the set of objects of $\C{C}$. The set of morphisms
( = arrows) with source $\rho \in \Ob \C{C}$ and target $\sigma \in \Ob \C{C}$
is denoted by $(\rho ,\sigma )$ and the identity morphism of $\rho $ by 
$\bf{1}_{\rho }$.
 
\begin{defi} \begin{enumerate}
\item[(i)] A category $\C{C}$ is called a $C^*$-category, if
the following hold:
\begin{enumerate}
\item[(a)] The space $(\rho ,\sigma )$ is a complex Banach space
for all objects $\rho,\, \sigma \in \Ob \C{C}$.
\item[(b)] The composition of morphisms gives  a bilinear map
$(S,R) \Pfeil S\circ R$ such that $\|S\circ R\| \leq \|S\| \cdot \|R\|$.
\item[(c)] There is an antilinear
involutive contravariant functor $*:\C{C} \pfeil \C{C}$ 
such that $\sigma =\sigma ^*$ for $\sigma \in \Ob \C{C}$
and  $\|R^*\circ R \|  =\|R\|^2$ for 
$R\in (\rho ,\sigma)$.
\item[(d)] $\dim \, (\rho ,\rho )\geq 1$ for every object
$\rho \in \Ob \C{C}$. 
\end{enumerate}
In particular,
$(\rho ,\rho )$ is a $C^*$-algebra with unit $\B{1}_{\rho }$ for every
object $\rho $ of a $C^*$-category.

\item[(ii)] Let $\C{C}$ be a $C^*$-category and 
$\rho ,\, \sigma \in \Ob \C{C}$.
A morphism $U\in (\rho ,\sigma )$ is called a linear isometry, if 
$U^*\circ U =\B{1}_{\rho }$. A linear isometry satisfying $U\circ U^*=
\B{1}_{\sigma }$ is called unitary.
\item[(iii)] Let $\C{C}$ be a $C^*$-category endowed with the following
product structure:
\begin{enumerate}
\item[(a)] For $\rho ,\, \sigma \in \Ob \C{C}$ there is a product object
$\rho \sigma $, and for morphisms \newl
$T\in (\rho ,\sigma )$ and 
$T'\in (\rho ',\sigma ')$ there is a morphism $T \times T'\in (\rho \rho ',
\sigma \sigma ')$.  The mapping $(T, T') \Pfeil T\times T'$ is  bilinear,
and we have
$(T\times T')^*=T^*\times (T')^*$, $\B{1}_{\rho }\times 
\B{1}_{\sigma }=\B{1}_{\rho \sigma}$ for $\rho ,\sigma \in \Ob \C{C}$ as 
well as the interchange law 
$$S\times S'\circ T\times T'=(S\circ T)\times (S'
\circ T')$$
whenever the left side is defined. (To save brackets we evaluate 
$\times $ before $\circ $.)
\item[(b)] For all objects $\rho ,\, \sigma $ and 
$\tau $ of $\C{C}$ there is a 
unitary operator $a(\rho ,\sigma ,\tau )\in 
(\rho (\sigma \tau ),\,(\rho \sigma )\tau )$
such that the pentagonal diagram
\begin{center}
\xext=1000 \yext=1000
\adjust[\phi (\rho (\sigma \tau ))`;`a(\phi ,\rho \sigma ,\tau );
`a(\phi \rho ,\sigma ,\tau );`a(\phi ,\rho ,\sigma )\times \B{1}_{\tau
}]
\begin{picture}(\xext,\yext)(\xoff,\yoff)
\settriparms[1`1`0;500]
\putAtriangle(0,500)[\phi (\rho (\sigma \tau ))`\phi ((\rho \sigma )\tau )`
(\phi \rho )(\sigma \tau );\B{1}_{\phi }\times a(\rho ,\sigma ,\tau )`
a(\phi ,\rho ,\sigma \tau )`]
\setsqparms[0`1`1`1;1000`500]
\putsquare(0,0)[\phantom{\phi ((\rho \sigma )\tau )}`
\phantom{(\phi \rho )(\sigma \tau )}`(\phi (\rho \sigma ))\tau `
((\phi \rho )\sigma )\tau ;
`a(\phi ,\rho \sigma ,\tau )`a(\phi \rho ,\sigma , \tau )`
a(\phi ,\rho ,\sigma )\times \B{1}_{\tau }]
\end{picture} \end{center}
commutes for objects $\phi ,\,\rho ,\,\sigma $ and $\tau $. Furthermore
$a(\rho , \sigma ,\tau)$ is natural in $\rho ,\sigma , \tau $. 
\item[(c)] There is a distinguished object $\iota $ in
$\C{C}$ (the unit object). For each
object $\rho $ in $\C{C}$ there are unitary operators $l_{\rho }\in
(\iota \rho, \rho )$ and $r_{\rho }\in (\rho \iota ,\rho )$ satisfying
the following properties:
\newl  They are natural 
in $\rho $, it holds $l _{\iota } =r _{\iota }$, and the diagram
\begin{center}
\settriparms[1`1`1;700]
\qtriangle[\rho (\iota \sigma )`(\rho \iota )\sigma `\rho \sigma ;a(\rho ,\iota
,\sigma )`\B{1}_{\rho } \times l_{\sigma }`r_{\rho }\times \B{1}_{\sigma }]
\end{center}
commutes for all $\rho ,\, \sigma \in \Ob \C{C}$.
\end{enumerate}
Then $(\C{C}, \times ,a ,l ,r)$ (or simply $\C{C}$) is said to be a
$C^*$-tensor category.
\item[(iv)] The $C^*$-tensor category $\C{C}$ is called strict if 
$\rho (\sigma \tau )=
(\rho \sigma )\tau $, $a(\rho ,\sigma ,\tau )=\B{1}_{\rho (\sigma \tau )}$,
$\rho \iota =\iota \rho =\rho $ and $l_{\rho }, \,r_{\rho }=\B{1}_{\rho }$
for all objects $\rho ,\,\sigma ,\,\tau $ of $\C{C}$. 
We also write $\rho \sigma \tau $ instead of $\rho (\sigma \tau )$ for objects
$\rho,\, \sigma ,\,\tau $ of a strict
$C^*$-tensor category. 
\item[(v)] A $C^*$-tensor category $\C{D}$ is called a (full) $C^*$-tensor
subcategory of the $C^*$-tensor category $\C{C}$ if $\Ob \C{D} \subset
\Ob \C{C}$, if the morphism spaces $(\rho ,\sigma )$ of $\C{D}$ coincide with
those of $\C{C}$ for $\rho ,\, \sigma \in \Ob \C{D}$, and if one obtains
the remaining structure of $\C{D}$ by restricting the $C^*$-tensor structure
of $\C{C}$ onto $\Ob \C{D}$. 
\end{enumerate} \label{DC}   \end{defi}    

\begin{defi} \begin{enumerate} \item[(i)]
Two objects $\rho $ and $\sigma $ of a $C^*$-category $\C{C}$ are called
equivalent if and only if there is a unitary operator
$U\in (\rho ,\sigma )$. An object $\rho $ is
called irreducible if and only if $(\rho ,\rho )=\comp \B{1}_{\rho }$.
Clearly 'equivalent' defines an equivalence relation on $\Ob \C{C}$. $[\rho ]$
denotes the equivalence class of the object $\rho $,  
$[\C{C}]$ the set of equivalence classes of $\Ob \C{C}$ and 
$[[\C{C}]]$ the set of equivalence classes of irreducible objects of
$\C{C}$. 
\item[(ii)] Let $\rho $ be an object of $\C{C}$ and $E$ a projection in 
$(\rho ,\rho )$. (In this paper projections are assumed to be orthogonal.)
An object $\sigma $ of $\C{C}$ is called a subobject of
$\rho $ corresponding to $E$ if there is a linear isometry 
$V\in (\sigma , \rho )$ such that $V\circ V^*=E$. The notation 
$\sigma \leq \rho $ means
that $\sigma $ is a subobject of $\rho $.
\item[(iii)] An object $\tau $ of $\C{C}$
is called the direct sum $\bigoplus _{i=1}^n \rho _i$ of the objects $\rho _i$
if there are linear isometries $V_i \in (\rho _i,\tau )$ for $i=1,\ldots ,n$
such that  $\sum _{i=1}^n V_i\circ V_i^* =\B{1}_{\tau }$.
\item[(iv)] The $C^*$-category $\C{C}$ is said to have subobjects (or to
be closed under subobjects) if
there is a subobject of $\rho $ corresponding to $E$ for every object
$\rho $ of $\C{C}$ and for every projection $E\in (\rho ,\rho )$.
$\C{C}$ is said to have (finite) direct sums if a direct sum $\rho \oplus
\sigma $ exists for all objects $\rho ,\sigma \in \Ob \C{C}$. 
\end{enumerate}
\label{D121} \end{defi}
 
\subsection{Remarks: \label{Resds} } \begin{enumerate}
\item[(1)] A subobject $\sigma $ of an object $\rho \in \Ob \C{C}$ 
corresponding
to a projection $E\in (\rho ,\rho )$ is unique
up to equivalence. If a linear isometry
$V\in (\sigma ,\rho )$ is chosen as in Definition
~\ref{D121} (ii) then
$$\alpha :(\sigma ,\sigma )\pfeil E\, (\rho ,\rho )\,  E,\,T\Pfeil 
V\circ T\circ
V^*, $$
defines a surjective isomorphism from the $C^*$-algebra $(\sigma ,\sigma )$
onto the $C^*$-algebra $E\, (\rho ,\rho )\, E$.
\item[(2)] Two projections $E$ and $F$ in $(\rho ,\rho )$ belong to the same
subobject $\sigma $ if and only if they are equivalent 
in $(\rho ,\rho )$ (i.e. there is 
a partial isometry $U\in (\rho ,\rho )$ such that $U^*\circ U=E$ and $U\circ
U^*=F$.)
\item[(3)] Let us assume that a direct sum $\tau =\oplus _{i=1}^n \rho _i$
of objects $\rho _i \in \Ob \C{C}$ is defined as in Definition
~\ref{D121} (iii). Let $\tau '$ be another direct sum of the objects 
$\rho _i $ ($i\in 1,\ldots ,n$) and let $V_i'\in (\rho _i,\tau ')$
be the corresponding linear isometries. Then there is a unitary operator
$U\in (\tau ,\tau ')$ such that 
$$U\circ V_i\circ V_i^*\circ U^* =V_i'\circ (V_i')^* \qquad 
   \hbox{($i=1,\ldots ,n$),} $$
in particular a direct sum $\bigoplus _{i=1}^n \rho _i$ is unique up to
equivalence.
\item[(4)] Let $\tau $ be a direct sum $\rho _1 \oplus \rho _2$, and let
$E:=V_1\circ V_1^*$ and
$F:=V_2\circ V_2^*$. There is a unique surjective linear isometric map
$\iota =\iota _{\rho_1 \rho _2}$ from
$(\rho _1,\rho _2)$ onto
$$F\, (\rho _1,\rho _1)\, E:=
\{x\in (\tau ,\tau ):F\circ x\circ E \,=x\} \subset (\tau ,\tau )$$
given by $\iota (T)= V_2\circ T\circ V_1^*$ for $T\in (\rho _1,\rho _2 )$.
We obtain $\iota _{\rho_1 \rho _2}(T)^*=\iota _{\rho_2 \rho _1}(T^*)$.
\item[(5)] Let $\C{C}$ be a $C^*$-tensor category. If $\rho _1$ and $\rho _2$
are objects of $\C{C}$ and if $\sigma _i$ is a subobject of
$\rho _i$ for $i=1,2$ then $\rho _1 \rho _2$ is a subobject of 
$\sigma _1 \sigma _2$.
\end{enumerate}

\abs
The proof of the Remarks is just an easy exercise.

\begin{defi} \begin{enumerate}
\item[(i)] Let $\C{C}$ and $\C{D}$ be $C^*$-categories and $F:\C{C}\pfeil
\C{D}$ a covariant functor. $F$ is called a $C^*$-functor if
$$T\in (\rho ,\sigma ) \Pfeil F(T)\in (F(\rho ), F(\sigma ))$$
is linear
for all objects $\rho ,\sigma \in \Ob \C{C}$ and
$F(T^*) =F(T)^*$ holds for $T\in (\rho , \sigma )$.
\item[(ii)]
 Let $\C{C}$ and $\C{D}$ be $C^*$-tensor categories, 
$\iota _{\C{C}}$ and 
 $\iota _{\C{D}}$ their unit objects and
$F:\C{C}\pfeil \C{D}$ a $C^*$-functor. For objects $\rho ,\sigma \in
\Ob \,\C{C}$ let $U_{\rho \sigma }$ be a unitary operator of
$((F(\rho )F(\sigma ), F(\rho \sigma ))$ and $J$ a unitary operator of 
$(F(\iota _{\C{C}}), \iota _{\C{D}})$ such that the following properties are
satisfied: \begin{enumerate}
\item[(a)] $U_{\rho \sigma }$ is natural in $\rho $ and $\sigma $.
\item[(b)] For $\rho ,\sigma ,\tau \in \Ob \C{C}$ the diagram
\begin{center} \hspace*{-2cm}
\xext=2800 \yext=700
\adjust[`U_{\rho \sigma }\times \B{1}_{\tau };
`\B{1}_{\rho }\times U_{\sigma \tau };
`U_{\rho \sigma }\times \B{1}_{\tau };
`U_{\rho ,\sigma \tau }]
\begin{picture}(\xext,\yext)(\xoff,\yoff)
\setsqparms[1`1`0`1;1600`500]
\putsquare(0,0)[F(\rho )(F(\sigma )F(\tau ))`(F(\rho )F(\sigma ))F(\tau )`  
F(\rho )F(\sigma \tau )`F(\rho (\sigma \tau ));
a(F(\rho ), F(\sigma ), F(\tau ))`\B{1}_{F(\rho )}\times U_{\sigma \tau }``
U_{\rho ,\sigma \tau }] 
\setsqparms[1`0`1`1;1200`500]
\putsquare(1600,0)[\phantom{(F(\rho )F(\sigma ))F(\tau )}`
F(\rho \sigma )F(\tau )`
\phantom{F(\rho (\sigma \tau ))}`
F((\rho \sigma )\tau ) ;
U_{\rho \sigma }\times \B{1}_{F(\tau )} ``U_{\rho \sigma ,\tau }`
F(a(\rho ,\sigma ,\tau ))]
\end{picture} \end{center}
commutes. 
\item[(c)] The diagram
\begin{center}
\setsqparms[1`1`1`1;1000`500]
\square[F(\iota _{\C{C}})F(\rho )`F(\iota _{\C{C}}\rho )` 
\iota _{\C{D}}F(\rho )`F(\rho );
U_{\iota _{\C{C}}\rho }`J\times \B{1}_{F(\rho )}`F(l_{\rho })`l_{F(\rho )}]
\end{center}
commutes for every object $\rho $ in $\C{C}$ just as the corresponding diagram 
for $r_{\rho }$ and $r_{F(\rho )}$.
\end{enumerate}
Then $(F, (U_{\rho \sigma })_{\rho ,\sigma }, J)$ is called a 
$C^*$-tensor functor. Usually we just write $F$ for the $C^*$-tensor
functor.
\item[(iii)] A $C^*$-tensor functor $F:\C{C} \pfeil \C{D}$ is called a 
($C^*$-tensor) equivalence of the $C^*$-tensor categories $\C{C}$ and $\C{D}$
if $F$ is full and faithful 
and if for each object $\tau $ of $\C{D}$ there is an object $\rho $
of $\C{C}$ such that $F(\rho )$ is equivalent to $\tau $.
\item[(iv)] Let $\C{C}$ and $\C{D}$ be strict $C^*$-tensor categories.
A covariant functor  \linebreak $F:\C{C} \pfeil \C{D}$ is called a strict
$C^*$-tensor functor if $F(\rho \sigma )=F(\rho )F(\sigma )$ for
$\rho,\, \sigma \in \Ob \C{C}$, $F(\iota _{\C{C}})=\iota _{\C{D}}$ and if
$(F,(\B{1}_{F(\rho \sigma )})_{\rho ,\sigma },\B{1}_{\iota _{\C{D}}})$
is a $C^*$-tensor functor. \end{enumerate}
\label{DCequ} \end{defi}

\subsection{Remarks: } \begin{enumerate}
\item[(1)] Let $F:\C{C}\pfeil \C{D}$ be a $C^*$-functor.
We obtain   $\|F(T)\| \leq \|T\|$ for 
all morphisms $T$ of $\C{C}$. If $F$ is 
faithful, then  $\|F(T)\| =\|T\|$ for all morphisms $T$ of $\C{C}$. 
(Use the corresponding results for homomorphisms of $C^*$-algebras.)
\item[(2)] Let $\C{C}$ be a $C^*$-tensor category. Remark (1) implies
$$\|\B{1}_{\rho }\times T\|\leq \|T\| \qquad \hbox{and} \qquad
 \|T\times \B{1}_{\rho }\| \leq \|T\|$$
for any object $\rho $ of $\C{C}$ and any morphism $T$ of $\C{C}$.
Hence
$\|S\times T\| \leq \|S\|\cdot \|T\|$ for
all morphisms $S$ and $T$ of $\C{C}$. \end{enumerate}

\abs
An equivalence $F:\C{C} \pfeil \C{D}$ induces a bijection $[F]$ from
the set $[\C{C}]$  onto the set $[\C{D}]$.
Equivalences satisfy the expected properties:

\begin{prop} \begin{enumerate}
\item[(i)] Let $\C{C}, \C{D}$ and $\C{E}$ be $C^*$-tensor categories. 
If \newl
$(F,(U_{\rho \sigma })_{\rho ,\sigma },J)$ is an equivalence from
$\C{C}$ to $\C{D}$ and $(G,(V_{\phi \psi})_{\phi ,\psi },K)$ is an equivalence
from $\C{D}$ to $\C{E}$, then
$(G\circ F, (G(U_{\rho \sigma })\circ V_{F(\rho ),F(\sigma )})_{\rho ,\sigma },
K\circ G(J) )$ is an equivalence from $\C{C}$ to $\C{E}$.
\item[(ii)] If $F:\C{C} \pfeil \C{D}$ is an equivalence of the $C^*$-tensor 
categories $\C{C}$ and $\C{D}$ then there is a $C^*$-tensor equivalence 
$G:\C{D} \pfeil \C{C}$ such that
$[G]$ is the inverse of $[F]$. \end{enumerate}
\label{PCTK} \end{prop}

Part (i) is easy to verify. Since the author failed to find a 
suitable reference for
(ii), a proof of (ii) is presented in the Appendix.

\abs
\begin{prop} Every $C^*$-tensor category is equivalent to a strict
$C^*$-tensor category. \label{P121} \end{prop}

First this result was proved in \cite{McL2} for tensor categories without
$C^*$-structure. (A more accessible reference is \cite{Ka}, Theorem XI.5.3.)
The transfer to $C^*$-tensor categories is obvious. 
The Proposition allows us to work only with strict $C^*$-tensor categories.
This will abbreviate many computations, as we don't need the
operators $a(\rho ,\sigma ,\tau )$, $l_{\rho }$ and $r_{\rho }$.

\begin{defi} \begin{enumerate}
\item[(i)]  Let $\rho $ be an object of the $C^*$-tensor 
category $\C{C}$. An object $\bar {\rho}$ is said to be conjugate to
$\rho $ if there are $R\in (\iota ,\ov{\rho}\rho )$ and 
$\ov{R}\in (\iota ,\rho\ov{\rho } )$ such that
\begin{eqnarray} l_{\rho }\circ \ov{R}^*\times \B{1}_{\rho } \circ 
a(\rho ,\orho ,\rho )\circ \B{1}_{\rho }\times R \circ r_{\rho }^*
   &=& \B{1}_{\rho }
                \qquad \mbox{and} \label{Rr1}\\  
l_{\orho }\circ R^*\times \B{1}_{\ov{\rho }} \circ a(\orho ,\rho ,\orho )\circ
\B{1}_{\ov{\rho } }\times \ov{R} \circ r_{\orho }^*&=& \B{1}_{\ov{\rho } }
               \label{Rr2} \end{eqnarray}
hold. $(R, \ov{R})$ is called a pair of conjugation operators for $\rho $ and
$\ov{\rho }$. If the $C^*$-tensor category $\C{C}$ is strict, (~\ref{Rr1})
and (~\ref{Rr2}) mean 
\begin{eqnarray} \ov{R}^*\times \B{1}_{\rho } \circ \B{1}_{\rho }
             \times R &=& \B{1}_{\rho }
                \qquad \mbox{and} \label{RR1}\\  
               R^*\times \B{1}_{\ov{\rho }} \circ \B{1}_{\ov{\rho } }
              \times \ov{R} &=& \B{1}_{\ov{\rho } }.
               \label{RR2} \end{eqnarray}
From now on, if $\rho $ is an object of a $C^*$-tensor category $\C{C}$
$\orho $ denotes an object of $\C{C}$ conjugate to $\rho $.
\item[(ii)] A strict $C^*$-tensor category is called 
a regular $C^*$-tensor category if it has subobjects and finite
direct sums and $(\iota ,\iota )=\comp \B{1}_{\iota }$ is satisfied for the
unit object $\iota $.
\item[(iii)] A regular $C^*$-tensor category $\C{C}$ is called a compact 
$C^*$-tensor category if every object 
$\rho $ of $\C{C}$ 
possesses a conjugate.
\item[(iv)] A compact $C^*$-tensor category $\C{C}$ 
is called a finite $C^*$-tensor
category if $[[\C{C}]]$ is finite.  \end{enumerate} 
\label{Dco}\end{defi} 

\begin{lem} Let $\C{C}$ and $\C{D}$ be $C^*$-tensor categories, let
$(F,(U_{\rho
\sigma })_{\rho ,\sigma },J)$ be a $C^*$-tensor functor from $\C{C}$ to 
$\C{D}$,
and let $\rho $ be an object of $\C{C}$ having a conjugate 
$\orho $. Then $F(\orho )$
is conjugate to $F(\rho )$. If $(R,\oR )$ is a pair of conjugation operators
for $\rho $ and $\orho $ then
$$(S, \ov{S})= (U_{\orho \rho }^*\circ F(R)\circ J^*,\, 
U_{\rho \orho }^*\circ F(\oR )\circ J^*)$$ is a pair of conjugation
operators for $F(\rho )$ and $F(\orho )$. \label{L12Rr}\end{lem}

\B{Proof: } We verify Equation (~\ref{Rr1}) for the pair $(S, \ov{S})$ and omit
the similar proof of (~\ref{Rr2}). (~\ref{Rr1}) for $(R, \oR)$ implies
\begin{equation}  F(l_{\rho })\circ F(\ov{R}^*\times \B{1}_{\rho }) \circ 
F(a(\rho ,\orho ,\rho ))\circ F(\B{1}_{\rho }\times R) \circ F(r_{\rho })^*
   = \B{1}_{F(\rho )}. \label{Rr3} \end{equation}
Using Definition ~\ref{DCequ} (ii), we replace
\begin{tabbing}
$F(a(\rho ,\orho ,\rho ))$ \qquad \= \kill 
$F(l_{\rho })$ \> by \qquad 
$l_{F(\rho )}\circ J\times \B{1}_{F(\rho )}\circ
U_{\iota _{\C{C}} \rho}^*$, \\
$F(\ov{R}^*\times \B{1}_{\rho })$ \> by 
\qquad $U_{\iota _{\C{C}} \rho} \circ  
F(\ov{R}^*)\times \B{1}_{F(\rho )} \circ U_{\rho \orho ,\rho }^*$,\\
$F(a(\rho ,\orho ,\rho ))$ \> by \qquad 
$U_{\rho \orho ,\rho }\circ 
U_{\rho \orho }\times \B{1}_{F(\rho )} \circ a(F(\rho ), F(\orho ), F(\rho ))
\circ $\\
\> \qquad \qquad \qquad \qquad \qquad \qquad
$\circ \B{1}_{F(\rho )}\times U_{\orho \rho }^*\circ 
U_{\rho ,\orho \rho }^*$, \\
and so on.
\end{tabbing}
In this way Equation (~\ref{Rr3}) yields
$$\displaylines{ 
l_{F(\rho )}\circ J\times \B{1}_{F(\rho )}\circ
F(\ov{R}^*)\times \B{1}_{F(\rho )} \circ
U_{\rho \orho }\times \B{1}_{F(\rho )} \circ a(F(\rho ), F(\orho ), F(\rho ))
\circ \hfill \cr
\hfill \circ \B{1}_{F(\rho )}\times U_{\orho \rho }^*\circ \B{1}_{F(\rho )}
\times F(R) \circ \B{1}_{F(\rho )}
\times J^* \circ r_{F(\rho )}^* \, =\, \B{1}_{F(\rho )},  }$$
and Equation (~\ref{Rr1}) has been shown for $S$ and $\ov{S}$. \blacksquare 

\abs
Definition ~\ref{Dco} (ii) is due to J.E. Roberts for strict 
$C^*$-tensor categories (see \cite{Ro}). He introduced a dimension 
for
objects possessing a conjugate in \cite{Ro}. In the following
we will give a brief summary (without proofs) of his theory.

In the following let $\C{C}$ be a regular $C^*$-tensor category and let
$\rho $ be
an object of $\C{C}$.
We start with the Frobenius reciprocity (see
Lemma 2.1 in \cite{LoRo}): 

\begin{lem} Assume that $\rho $ has a conjugate $\orho $, let $(R,\ov{R})$ be 
a pair of conjugation
operators for $\rho $ and $\orho $ and let $\sigma $ and 
$\tau $ be objects of $\C{C}$.
\newl $S\Pfeil  \B{1}_{\orho }\times S\circ R\times \B{1}_{\sigma }$ is a
bijective linear map from $(\rho \sigma ,\tau )$ onto $(\sigma ,\orho \tau )$.
The inverse is given by $S'\Pfeil \oR ^*\times \B{1}_{\tau } \circ 
\B{1}_{\rho }
\times S'$.

Similarly
$T\in (\sigma \rho , \tau ) \Pfeil T\times \B{1}_{\orho }\circ 
\B{1}_{\sigma }\times \oR \in (\sigma , \tau \orho )$ and 
\newl $T'\in (\sigma ,\tau \orho ) \Pfeil \B{1}_{\tau }\times R^* \circ
 T'\times \B{1}_{\rho } \in (\sigma \rho ,\tau )$ are inverse linear
maps.  \label{LFrob}\end{lem}

We note the following results:
\begin{prop} \begin{enumerate}
\item[(i)] Let $\orho $ be conjugate to $\rho $. An object $\tau $ of $\C{C}$
is conjugate to $\rho $ if and only if $\orho $ and $\tau $ are 
equivalent.
\item[(ii)] If $\rho $ is irreducible and $\orho $ conjugate to $\rho $, 
$\ov{\rho }$ is irreducible, too.
\item[(iii)] If $\rho $ has a conjugate the $C^*$-algebra $(\rho ,\rho )$ is 
finite dimensional, in particular
$\rho $ is a finite direct sum of irreducible objects.  
\item[(iv)] Let $\sigma $ and $\rho $ be irreducible objects of $\C{C}$. 
$\rho $ is conjugate to $\sigma $ if and only if there are operators
$R\in (\iota ,\sigma \rho )$ and $\oR \in (\iota ,\rho \sigma )$ such that
\linebreak
$\oR ^*\times \B{1}_{\rho } \circ \B{1}_{\rho }  \times R \not= 0$. 
\end{enumerate} \label{P12x}\end{prop}

R. Longo and J.E. Roberts proved Proposition ~\ref{P12x} in \cite{LoRo}
((i), (ii) after Lemma 2.1, (iii) in Lemma 3.2 and (iv) in Lemma 2.2).

\subsection{Remarks:} \begin{enumerate}
\item[(1)] Each minimal projections of $(\rho ,\rho )$ corresponds
to an
irreducible subobject $\sigma $ of $\rho $. Two minimal projections
$E$ and $F$ of $(\rho ,\rho )$ correspond to the same subobject $\sigma $
of $\rho $, if and only if they are equivalent in $(\rho ,\rho )$. 
If $A$ is a simple direct summand of $(\rho ,\rho )$ and if $\phi $ is 
an irreducible subobject of $\rho $ belonging to a 
minimal projection $E$ of $A$, then we attach the equivalence class
$[\phi ] \in [[\C{C}]]$ to $A$. In this way we get a bijective 
correspondence between the equivalence classes of the irreducible objects of
$\C{C}$ contained in $\rho $ and the simple
direct summands of $(\rho ,\rho )$.

\item[(2)] Let $\rho $ be an object of $\C{C}$ having a conjugate.
$\B{1}_{\rho }$ has a decomposition 
$\B{1}_{\rho }=\sum _{i=1}^n P_i$ into minimal projections of $(\rho ,\rho )$.
For $i=1,\ldots ,n$ let $\rho _i$ denote an irreducible object corresponding
to $P_i$ ($i=1,\ldots ,n$). Then $\rho $ is a direct sum $\bigoplus _{i=1}^n
\rho _i$ of these objects.

If $\B{1}_{\rho }=\sum _{j=1}^m Q_j$ is another decomposition of 
$\B{1}_{\rho }$ into minimal projections, then $m=n$ and there exists a
permutation $\pi $ of $\{1,\ldots ,n\}$ such that $P_i$ is equivalent to
$Q_{\pi (i)}$ for $i=1,\ldots ,n$. Hence the decomposition of $\rho $
into irreducible objects is unique in the following meaning:
If $\rho =\bigoplus _{i=1}^n \rho _i$ and $\rho 
=\bigoplus _{j=1}^m \rho _j'$ are two
decompositions of $\rho $ into irreducible objects, then $n=m$ and there is
a permutation $\pi $ of $\{1,\ldots ,n\}$ such that $\rho _i$ and 
$\rho _{\pi (i)}'$ are equivalent for $i=1,\ldots ,n$.
\end{enumerate}

\subsection{The statistical dimension \label{S12d}}  
We assume that $\rho $ has a conjugate $\orho $. For the moment let us 
suppose that 
$\rho $ is irreducible. Let $(S,\ov{S})$ be a
pair of conjugation operators for $\rho $ and $\orho $.
In general the operators $S^*\circ S$ and $\ov{S}^*\circ \ov{S}$ are not equal.
But there is a complex number $\alpha \in \comp \setminus \{0\}$ such that 
$R=\alpha S$ and $\oR 
= \frac{1}{ \ov{\alpha }}\, \ov{S}$ are a pair of standard conjugation 
operators, i.e. that 
$$R^*\circ R= 
\ov{R}^*\circ \ov{R}$$ holds. 
By identifying $(\iota ,\iota )=\comp \B{1}_{\iota }$ with $\comp $,
we may regard $R^*\circ R$ as a positive number,
which is called the (statistical) dimension $d(\rho )$ of $\rho $. 
The statistical
dimension $d(\rho )$ does not depend on the choice of the conjugate $\orho $ 
and the standard conjugation operators $R$ and $\oR $ (see below).
We point out that the dimension is not a natural number in general.

We are going to define standard conjugation operators and the dimension, if
$\rho $ is not necessarily irreducible. There are irreducible
objects $\sigma _1,\ldots , \sigma _n$ and linear isometries
$W_i\in (\sigma _i,\rho )$, $\ov{W_i}\in (\ov{\sigma _i}, \ov{\rho })$
($i=1,\ldots ,n)$ 
satisfying
$$\B{1}_{\rho } =\sum _{i=1}^n W_i\circ W_i^* \qquad
\hbox{and} \qquad \B{1}_{\ov{\rho }}=
\sum _{i=1}^n \ov{W_i}\circ \ov{W_i}^*.$$ Moreover, we have
standard conjugation operators 
$R_i\in (\iota ,\ov{\sigma }_i\sigma _i)$
and \newl $\ov{R}_i\in
(\iota ,\sigma _i\ov{\sigma }_i)$ for $\sigma _i$ and
$\osigma _i$ ($i=1,\ldots ,n$).
Put $R:=\sum _{i=1}^n \ov{W_i} \times W_i \circ R_i$ and $\ov{R}:=
\sum _{i=1}^n W_i \times \ov{W_i} \circ \oR _i$.
$(R,\, \oR )$ is a pair of conjugation operators for $\rho $ and $\orho $. 
Pairs of 
conjugation operators are called standard, 
when they are obtained in  this way. 
We get $R^*\circ R =
\ov{R}^*\circ \ov{R}$, and this number is called the statistical
dimension $d(\rho )$ of
$\rho $. $d(\rho )$ is independent of the choice of the standard pair
$(R, \,\oR )$ according to
the following Lemma (see Lemma 2 in \cite{Ro}):

\begin{lem} Let $\tilde{\rho }$ be any conjugate of $\rho $ and
let $T\in (\iota ,\tilde{\rho } \rho )$ be given. There is an operator
$\ov{T}\in (\iota, \rho \tilde{\rho } )$ such that $(T,\ov{T})$ is 
a standard pair of 
conjugation operators, if and only if there is a unitary operator
$U\in (\orho ,\tilde{\rho } )$ such that $(U\times \B{1}_{\rho }) \circ R =T$.
\label{Lconj}
\end{lem}

The conjugate object and the statistical
dimension have the following properties:

\begin{thm} Let $\rho $ and $\sigma $ be objects of $\C{C}$
having conjugates $\orho $ and $\osigma $. \begin{enumerate}
\item[(i)] $\rho \oplus \sigma $ is conjugate to $\orho \oplus \osigma $ and
$d(\rho \oplus \sigma )= d(\rho ) +d(\sigma )$.
\item[(ii)] $\rho $ is also conjugate to $\orho $ and
$d(\rho ) = d(\orho )$ holds.
\item[(iii)] Every subobject of $\rho $ has a conjugate.
\item[(iv)]  The product $\bar{\sigma }\bar{\rho }$ is conjugate to
$\rho \sigma $. If $(R_{\rho }, \oR _{\rho })$ is a standard pair of
conjugation operators for $\rho $ and $\orho $ and $(R_{\sigma }, 
\oR _{\sigma })$ for $\sigma $ and $\osigma $, then
\begin{equation} (R_{\rho \sigma }, \oR _{\rho \sigma }) =
(\B{1}_{\osigma }\times R_{\rho } \times \B{1}_{\sigma } \circ R_{\sigma },
\B{1}_{\rho }\times \oR _{\sigma }\times \B{1}_{\orho } \circ \oR _{\rho })
      \label{conpr}              \end{equation}
is a standard pair of conjugation operators for $\rho \sigma $ and
$\osigma \orho $. The equation
$d(\rho  \sigma )= d(\rho ) \cdot d(\sigma )$ is satisfied. 
\end{enumerate}\label{Dim}
\end{thm}
  
Part (i) and (ii) are almost obvious. Part (iii) follows from
Theorem 2.4 and Part (iv) from Lemma 3.6 and Corollary 3.10 (a) in
\cite{LoRo}.

\subsection {The $C^*$-tensor category of the unitary corepresentations 
of a finite dimensional Hopf-$*$-algebra}
The most natural example for a compact $C^*$-tensor category is the category
$\C{U}_G$ of the unitary finite dimensional representations of a 
compact group $G$. In \cite{LoRo} R. Longo and J.E. Roberts deal with the 
more general case of the unitary finite dimensional corepresentations of a 
compact matrix pseudogroup in the sense of S.~Woronowicz (\cite{Wor}). In this
article we only need the special case of a finite dimensional Hopf-$*$-algebra.

We assume
that the reader is familiar with Hopf-$*$-algebras and mention \cite{Sw},
\cite{Ka}, \cite{Koo}, and \cite{Wor} as basic references. 
Let us introduce the finite $C^*$-tensor
category $\C{U}_{H}$ of the unitary corepresentations
of the finite dimensional Hopf-$*$-algebra $H$.
$\Delta :H \pfeil H\otimes H$ denotes the comultiplication of $H$,
$\epsilon :H\pfeil \comp $ the counit of $H$ and $S:H\pfeil H$ the antipode
of $H$.
The objects of $\C{U}_H$ are the unitary 
corepresentations $\sigma :V_{\sigma }\pfeil 
V_{\sigma } \otimes H$ of $H$ on a (small) finite dimensional
Hilbert space $V_{\sigma }$. We recall that a unitary corepresentation 
$\sigma :V_{\sigma }\pfeil 
V_{\sigma } \otimes H$ is a linear map satisfying
$$(\sigma \otimes id_H)\circ \sigma \,= 
  (id_{V_{\sigma }} \otimes \Delta )\circ \sigma \qquad \hbox{and} \qquad
  (id_{V_{\sigma }} \otimes \epsilon ) \circ \sigma \,=id_{V_{\sigma }} $$
as well as
 \begin{equation}\sum _{(v),(w)} \ska{v^{(1)}}{w^{(1)}}\, w^{(2)*}v^{(2)}
 \, = \ska{v}{w}\,  \B{1} \qquad \hbox{for $v,w \in V_{\sigma }$}
   \label{Uni}       \end{equation}
where $\sigma (v)=\sum _{(v)} v^{(1)} \otimes v^{(2)}$ (Sweedler notation
for corepresentations).
If $\rho :V_{\rho }\pfeil V_{\rho }\otimes H$ and
$\sigma :V_{\sigma }\pfeil 
V_{\sigma } \otimes H$
are objects of  $\C{U}_H$
$(\rho ,\sigma )$ is the Banach space of the intertwining operators
$T:V_{\rho }\pfeil V_{\sigma }$ of the corepresentations
$\rho $ and $\sigma $ (i.e. $T$ satisfies 
$(T\otimes id_H)\circ \rho = \sigma \circ T$.)
The composition $\circ $ is the composition of maps and 
$T^*\in (\sigma ,\rho )$ is the adjoint operator of 
$T\in (\rho ,\sigma )$.
The tensor product $\rho \sigma $ 
is the tensor product $\rho \otimes
\sigma $ of the corepresentations $\rho $ and $\sigma $
defined by
$$\rho \otimes \sigma :V_{\rho }\otimes V_{\sigma } \pfeil 
V_{\rho }\otimes V_{\sigma } \otimes H,\,
  (\rho \otimes \sigma )(v\otimes w)\, =\sum _{(v),(w)}
  v^{(1)} \otimes w^{(1)}\otimes v^{(2)}w^{(2)},$$
 and $T\times T'$ is the tensor product
$T\otimes T'$ of linear operators for morphisms $T\in (\rho ,\sigma )$ and
$T'\in (\rho ',\sigma ')$. The unit object is the identity corepresentation
$\iota :\comp \pfeil \comp \otimes H,\, \iota (\lambda )=\lambda \otimes 
 \B{1}$.
The associativity constraint 
$a(\rho ,\sigma ,\tau )$ is the canonical unitary operator
from $V_{\rho }\otimes  (V_{\sigma }\otimes V_{\tau })$
onto $(V_{\rho }\otimes  V_{\sigma })\otimes V_{\tau }$,
and $l_{\rho }$ (resp. $r_{\rho }$) is the canonical unitary
operator from $\comp \otimes V_{\rho }$ (resp. $V_{\rho }\otimes \comp $)
onto $V_{\rho }$. Using
the usual identifications of these Hilbert spaces, we may regard $\C{U}_H$ as 
a strict
$C^*$-tensor category. 

The conjugate object for an object $\sigma \in \Ob \C{U}_H$ is the 
contragredient corepresentation $\osigma $ of $\sigma $. Let 
$\C{B}=(v_1,\ldots ,v_n)$ be an orthonormal base 
of $V_{\sigma }$ and $(a_{ij})_{i,j=1}^n$ be the matrix coefficients of
$\sigma $ with respect to $\C{B}$, i.e.
$\sigma (v_j)=\sum _{i=1}^n v_i \otimes a_{ij}$ for $i=1,\ldots ,n$. 
\vspace*{3pt} Then
$\osigma $ is given on the dual Hilbert space 
$\ov{V_{\sigma }}$
by $\osigma (\ov{v_j})=\sum _{i=1}^n \ov{v_i} \otimes S(a_{ji})$,
where $(\ov{v_1},\ldots ,\ov{v_n})$ is the base dual to $\C{B}$.
$\osigma $ is a unitary corepresentation of $H$ with respect to the inner
product on $V_{\sigma }$. (According to \cite{Koo}, a corepresentation
$\sigma $ fulfils the unitarity property (~\ref{Uni}) if and only if
\[ S(a_{ij}) = a_{ji}^* \qquad \hbox{for $i,j=1,\ldots ,n$.}
\label{Suni} \]
Using this characterization of the unitarity and applying $S^2=id_H$,
one concludes that $\osigma $ is unitary.)

Observe that $\sigma $ is unitary if and only if the matrix coefficients
satisfy 
$$\sum _{j=1}^n a_{ij}a_{kj}^* =\delta _{ik} {\bf
1} \qquad \hbox{for $i,k=1,\ldots ,n$}$$
(see \cite{Koo}). Hence
$$(\sigma \otimes \osigma )\sum _{j=1}^n v_j\otimes \ov{v_j} \,=
\sum _{i,k=1}^n v_i\otimes \ov{v_k} \otimes \sum _{j=1}^n
a_{ij} \cdot a_{kj}^* =
\sum _{i=1}^n v_i\otimes \ov{v_i} \otimes \B{1}.$$ Thus
the restriction of $\sigma \otimes \osigma $ onto 
$\comp \sum _{j=1}^n v_j\otimes \ov{v_j}$ is equivalent to the identity
corepresentation $\iota $ and 
$$\oR :\comp \pfeil V_{\sigma }\otimes \ov{V_{\sigma }},\, 
\lambda \Pfeil \lambda \sum_{j=1}^n
  v_j \otimes \ov{v_j},$$
intertwines the corepresentations $\iota $ and $\sigma \otimes \osigma $.
Similarly,  
$$R :\comp \pfeil \ov{V_{\sigma }}\otimes V_{\sigma },\, \lambda 
\Pfeil \lambda \sum_{j=1}^n
  \ov{v_j} \otimes v_j,$$
is an intertwining operator for the corepresentions
$\iota $ and $\osigma \otimes \sigma $.
Now we get 
\begin{eqnarray*}
 (\oR ^*\otimes \B{1}_{\sigma } \circ \B{1}_{\sigma }\otimes R) \,w
& =& (\oR ^*\otimes \B{1}_{\sigma }) \,w \otimes \sum _{j=1}^n
 \ov{v_j}\otimes v_j\,\, = \\
\sum _{j=1}^n \skab{w\otimes \ov{v_j}}{\, \sum _{i=1}^n v_i\otimes \ov{v_i}}
\,v_j &= &\sum _{j=1}^n \ska{w}{v_j}\, v_j\,\, =\qquad w  \end{eqnarray*}
for $w\in V_{\sigma }$
and conclude $\oR ^*\otimes \B{1}_{\sigma }\circ \B{1}_{\sigma } \otimes R
=\B{1}_{\sigma }$. The same argument yields 
$R^* \otimes \B{1}_{\osigma }\circ \B{1}_{\osigma } \otimes \oR
=\B{1}_{\osigma }$.
As $R^*\circ R =n = \oR ^*\circ \oR$, we have
$d(\sigma )=n=\dim V_{\sigma }$, if $\sigma $ is irreducible. 
By applying Theorem
~\ref{Dim} (i) we obtain $d(\sigma )=\dim V_{\sigma }$ for every object 
$\sigma \in
\Ob \C{U}_H$. We point out that the equation $d(\sigma )=\dim V_{\sigma }$ 
is not valid for general compact quantum groups, as $\osigma $ need not be
unitary with respect to the dual inner product on
$\ov{V_{\sigma }}$. (Observe that $S^2 =id_H$ is not satisfied in general.)

\section{The $C^*$-tensor category of the $(N,N)$-bimodules for
subfactors $N\subset M$}
First let us fix some notations concerning subfactors and the basic
construction. In order to avoid
subtleties we will assume that every von Neumann algebra
appearing in this paper acts on a separable Hilbert space.

\subsection{The basic construction and the Jones' tower}
Let $A\subset B$ be an inclusion of two finite von Neumann algebras $A$ and
$B$ with the same unit, and let $\tr$ be a faithful normal normalized
trace of $B$. $L^2(B)$ is defined as the completion of $B$
with respect to the inner product
$(x,y)\Pfeil \tr (xy^*)$. We denote an element $x$ of $B$
by $\ov{x}$ if $x$ is regarded as an element of $L^2(B)$.
The Hilbert space $L^2(A)$ (defined by the trace $\tr \rest A$)
is a closed subspace of
$L^2(B)$. $B$ acts normally and faithfully on $L^2(B)$ by left multiplication.
So we may consider $B$ as a von Neumann algebra on $L^2(B)$.

Let $e$ denote the orthogonal projection from $L^2(B)$ onto
$L^2(A)$. The von Neumann algebra $B_1$ generated by $B$ and $e$ is called
the basic construction for $A\subset B$ (and for the trace $\tr $). $e$ is 
called the Jones projection. 
 
$e$ maps $B$ onto $A$, the restriction 
$$E=e\rest B:B \pfeil A$$ of $e$ is a normal faithful conditional expectation
from $B$ onto $A$, in particular we have $E(a_1ba_2)=a_1E(b)a_2$ for 
$a_1,\,a_2\in A$ and
 $b\in B$
and $E(b)\geq 0$ for 
$b\geq 0$. For $b\in B$, $E(b)$ is the unique element of $A$ satisfying 
\begin{equation}
\tr (E(b)\, a)=\tr (b\, a) \qquad \text{for every $a\in A$}. \label{cond}
\end{equation}
$E$ is called the conditional expectation from $B$ to $A$ corresponding
to the trace $\tr $.

Let $B^{op}$ be the von Neumann algebra opposite to B. ($B^{op}$ is equal to
$B$ as a complex vector space, the multiplication law is reversed, that
means $b\stackrel{\rm op} \circ c:=c\cdot b$ for $b,c\in B$, and the
involution $*$ is the same as in $B$.) We have a normal representation
$\rho :B^{op} \pfeil \B{L}(L^2(B))$ of $B^{op}$ given by
$\rho (b)\, \ov{x} =\ov{xb}$ for $b,x\in B$ (multiplication from right).   
The basic construction satisfies the relation
\begin{equation} B_1=\rho (A^{op})'. \label{fund} \end{equation}
Let $J_B$ be the involutive antiunitary operator on $L^2(B)$ given by 
$J_B\, \ov{x} =\ov{x^*}$
for $x\in B$. Often we abbreviate $J_B$ to $J$.
We get $\rho (b) =Jb^*J$ for $b\in B$. From equation (~\ref{fund})
we conclude
\begin{equation} B_1 =JA'J. \label{fundJ} \end{equation}

For a $\zws $-factor $L$, let $\tr _L$ denote the unique normalized 
trace of $L$.
Now let $N\subset M$ be an inclusion of $\zws $-factors (with the same unit).
The basic construction $M_1$ (with respect to $\tr _M$)
is a factor of type $\zws $ if and only if the Jones index
$\beta :=[M:N] < \infty $. From now on let us assume
that $\beta $ is finite. 
In this case $M_1$ is the $\comp $-vector space generated by
$m_1em_2,\,\, m_1,\, m_2\in M$ (see \cite{GHJ}, Theorem 3.6.4).
The subfactor $M\subset M_1$ is called the subfactor dual to $N\subset M$.
We obtain  $[M_1:M] =[M:N]
< \infty $, so we are able to repeat the basic construction 
infinitely many times and get the so called Jones
tower
$$N=M_{-1}\subset M=M_0 \subset M_1 \subset M_2 \subset M_3 \subset
 \ldots  ,$$
where $M_{k+1}$ is the basic construction for the subfactor
$M_{k-1} \subset M_k$. 
For $k\in \nat \cup \{0\}$ let
$e_k \in M_{k+1}$ denote the orthogonal
projection from $L^2(M_k)$ onto $L^2(M_{k-1})$ and 
$E_k:M_k \pfeil M_{k-1}$ the corresponding 
conditional expectation. We have the following relations:
\begin{eqnarray} e_k\, e_{k\pm 1}\,e_k \,&=& \beta ^{-1}e_k \qquad \hbox{and}
                     \label{E1}            \\
                  e_k\,e_l\,&=& e_l\,e_k                    
              \qquad \hbox{for $|k-l|\geq 2$.}\label{E2}
\end{eqnarray}
Moreover we have
\begin{eqnarray} e_k\,x\,e_k\,&=&E_k(x)e_k \qquad \hbox{for $x\in M_k$\ \ \
 and} \label{E3} \\ 
                   E_k(e_{k-1})&=& \beta ^{-1}\B{1}.  \label{E4}\end{eqnarray}

The trace $\tr _{M_{k+1}}$ on $M_{k+1}$ satisfies the Markov property
\begin{equation} \beta \, \tr _{M_{k+1}}(xe_k)=\tr _{M_{k+1}}(x)  \qquad 
\hbox{for $x\in M_k$.} \label{markov}  \end{equation}

\subsection{$C^*$-tensor categories and bimodules}
Let $N$ be a $\zws $-factor. An easy access to $(N,N)$-bimodules may be found
in \cite{Sch3} (see also \cite{DH1}, \cite{DH2}, \cite{Y1} and \cite{Y2}). 
We assume the notation from \cite{Sch3} and recall that an $(N,N)$-bimodule
$_N\C{H}_N$ is called regular if the intersection of the left and of the 
right bounded elements of $\C{H}$ is dense in $\C{H}$.
For every $(N,N)$-bimodule $_N{\cal H}_N$ there is the conjugate bimodule
$_N\bar{\cal H}_N$, where
 $\bar{\cal H}$ is equal to $\cal H$ as a real vector space, but the inner 
product and the scalar
multiplication are conjugate. The left action $\bar {\lambda }$ of $N$ and 
the right action $\bar {\rho }$
of $N$ on $\bar{\cal H}$ are given by
$\bar{\lambda }(n)=\rho (n^*)$  and $\bar{\rho }(n)=\lambda (n^*)$
for $n\in N$.

We introduce the $C^*$-tensor category 
$\C{B}_N$ of all
(small) regular $(N,N)$-bimodules. The objects are the (small) regu\-lar 
$(N,N)$-bi\-mo\-dules.
For two objects $\rho ={_N\C{H}_N}$ and $\sigma ={_N\C{K}_N}$
of $\C{B}_N$, the 
$(N,N)$-linear
continuous operators from $\C{H}$ into $\C{K}$ 
(endowed with the operator norm) form the 
complex Banach space
$(\rho ,\sigma )$. Endowed with the usual $*$-operation
$*:(\rho ,\sigma ) \pfeil (\sigma ,\rho )$, $\C{B}_N$ is a 
$C^*$-category with subobjects and
finite direct sums. 

The product structure is given by the tensor product $\tn $. We define
$T\times T':=T\tn T'$ for $T\in (\rho ,\sigma )$ and 
$T'\in (\rho ',\sigma ')$. The unit object $\iota $ is $_NL^2(N)_N$. The
associativity constraint $a(\rho ,\sigma ,\tau )$ as well as
the maps $l_{\rho }$ and $r_{\rho }$ are defined as 
in \cite{Sch3}, Section 1.

Let $N\subset M$ be an inclusion of type $\zws $-factors such that 
$[M:N]< \infty $. Let $\C{B}_{N\subset M}$ denote the full $C^*$-tensor
subcategory of $\C{B}_N$ whose object set consists of all (small) 
$(N,N)$-bimodules which are equivalent to a finite direct sum of
$(N,N)$-bimodules contained in $_NL^2(M_k)_N$ for some $k\in \nat \cup \{0\}$.
Observe $L^2(M_{k-1}) \tn L^2(M_{l-1}) \cong L^2(M_{k+l-1})$ 
(for example see \cite{Sch3}), hence
the $\tn $-tensor product of two objects of $\C{B}_{N\subset M}$
is an object of $\C{B}_{N\subset M}$.

From now on we will treat the $C^*$-tensor category $\C{B}_{N\subset M}$, 
as if it were
strict. (This is allowed by Proposition \ref{P121}.)
The following result is probably known to some experts, but the author
could not find any reference.

\begin{prop} The $C^*$-tensor category $\C{B}_{N\subset M}$ is compact.
It is finite if and only if the subfactor $N\subset M$ has finite depth.

If $\C{H}$ is an $(N,N)$-bimodule of $\C{B}_{N\subset M}$ 
the conjugate bimodule $\bar{\C{H}}$ 
is conjugate to $\C{H}$ in the
sense of Definition ~\ref{Dco}. We have 
\begin{equation} d(\C{H})=\sqrt{c(\lambda (N),\C{H})\cdot
c(\rho (N),\C{H})}\label{dbim} \end{equation}
for any irreducible $(N,N)$-bimodule $_N\C{H}_N \in \Ob \C{B}_{N\subset M}$
(where $c(\lambda (N),\C{H})$ denotes the coupling constant 
of $\lambda (N)$ in \C{H}). 
\label{PBim}
\end{prop}

\subsection{Remarks: } \begin{enumerate}
\item[(1)] For any bimodule $\C{H}\in \C{B}_{N\subset M}$,
$d(\C{H})$ is equal to the minimal dimension of $\C{H}$ which is
the square root of the minimal index
$[\rho _{\C{H}}(N)':\lambda _{\C{H}}(N)]_{min}$
for the subfactor $\lambda _{\C{H}}(N)\subset \rho _{\C{H}}(N)'$. 
(We will not define
the minimal index here, see \cite{Hav}, \cite{Hi} or \cite{Lo1} for
the definition). 
\item[(2)] For the so-called extremal subfactors 
$N\subset M$ (see \cite{P}), we have 
 $$d(\C{H}) =c(\lambda (N),\C{H})= c(\rho (N),\C{H})$$ 
for every $\C{H}\in \Ob \C{B}_{N\subset M}$. 
We point out that
finite depth subfactors
and irreducible subfactors are extremal.
\item[(3)] If $N\neq M$ every object of
$\C{B}_{N\subset M}$ is indeed equivalent to a sub-bimodule of $L^2(M_k)$ 
for a suitable
$k\in \nat $. \end{enumerate}

The Remarks (1) and (2) are shown at the end of the proof of the Proposition, 
where we make
free use of S. Popa's definitions and results in \cite{P}. The easy proof of
(3) is left to the reader.

\abs
\B{Proof of the Proposition:} It suffices to show the assertions
for a subbimodule $\C{H}$ of $_NL^2(M)_N$, 
as it is possible to carry out the same proof for $M_k$ ($k\in \nat $) instead
of $M$.

Let $\C{H}$ be $pL^2(M)$ where $p$ is an orthogonal projection of 
$N'\cap M_1$. 
$M$ possesses a finite base
as a right $N$-module, as it was shown in \cite{PP}, Proposition 1.3.
A version of this result may be found in \cite{GHJ}, Theorem 3.6.4, we will 
modify the proof there in
order to get a right $N$-module base also for $pM$ ($\subset pL^2(M)$):
\newl As in the proof in \cite{GHJ} there are finitely many 
partial isometries $w_1,\ldots ,w_m\in M_1\setminus \{0\}$
such that 
\begin{equation} w_j^*\, w_j \leq e_0    \label{pis}\end{equation}
for $j=1,\ldots ,m$ and $\sum _{j=1}^m w_j\, w_j^* =p$. 
Similarly we find partial 
isometries
$w_{m+1},\ldots ,w_n \in M_1\setminus \{0\}$ such that (~\ref{pis}) 
holds for $j=m+1,\ldots ,n$ and
$\sum _{j=m+1}^n w_j\, w_j^* =\B{1}-p$. 
For $j=1,\ldots ,n$ there is a unique $v_j\in M$ such that
$w_j=v_je_0$. As in \cite{GHJ} we see that $(v_i: \,i=1,2,\ldots ,n)$ fulfils
the properties of a Pimsner-Popa basis:
\begin{itemize}
\item[(a)] $E_0(v_i^*v_j)=0$ for $i\not= j$.
\item[(b)] $f_i:=E_0(v_i^*v_i)$ is a projection satisfying 
$v_if_i=v_i$.
\item[(c)] Every $x\in M$ has a unique expansion 
$x=\sum _{j=1}^n v_jy_j$ with $y_j\in N$. In fact,
we have $y_j=E_0(v_j^*x)$.
\end{itemize}
Furthermore, we obviously get the following relations:
\begin{itemize}
\item[(d)] $\sum _{j=1}^m v_je_0v_j^* =p$ \qquad and
\item[(e)] $p\, . \ov{x} =\sum _{j=1}^m v_je_0v_j^*\, .\ov{x}\,=
\ov{\sum _{j=1}^m v_jE_0(v_j^*x)}$ \qquad for $x=\sum _{j=1}^n v_jy_j\in M$. 
\end{itemize}
From (e) we conclude that
\begin{itemize}
\item[(f)] $pL^2(M)\cap M$ is dense in $pL^2(M)$. 
\end{itemize}
The conjugate bimodule
$\ov{\C{H}}$ is isomorphic to $p^{op}L^2(M)$, where $p^{op}=JpJ$ and $J$ is 
the antiunitary operator $J_M$. We get 
$p^{op} L^2(M) = [\, \ov{x^*}:\, x\in pL^2(M)\cap M]$ 
and identify $\ov{\C{H}}$ and 
$p^{op} L^2(M)$.

Using the canonical isomorphism between $L^2(M_1)$ and
$L^2(M)\tn L^2(M)$, we regard $\C{K}:=pL^2(M)\tn p^{op}L^2(M)$ as an 
$(N,N)$-sub-bimodule of $L^2(M_1)$. $\ov{p}\in L^2(M_1)$ belongs to $\C{K}$,
as $\ov{p}= \beta ^{-1/2} \sum_{j=1}^m \ov{v_j} \tn \ov{v_j^*}$ by Property
(d) and $\ov{v_j}\in p\, L^2(M)$ for $j=1,\ldots ,m$ 
according to  Property (e). 
Since $p$ and $N$ commute,
$$\oR:L^2(N)\pfeil  \C{K}=pL^2(M)\tn p^{op} L^2(M),\,\,\ov{n}\Pfeil 
\beta ^{1/2}\, \ov{np}\in \C{K} \subset L^2(M_1),$$
defines a continuous $(N,N)$-linear operator. For $x\in \C{K}\cap M_1$ we get
\begin{eqnarray*} \ov{(\B{1}-p)x} & =& (\B{1}-p)\, .\ov{x} 
 \\ 
                   &\in &(\B{1}-p)\, . \,(pL^2(M)\tn L^2(M)) \,\,= 
                   \\
                   & = & (\B{1}-p)p L^2(M)\tn L^2(M) =\{0\} \end{eqnarray*}
(by applying Theorem 2.2 and Corollary 2.3 in \cite{Sch3}). 
Hence $x=px$. Moreover we have  
$$\displaylines{ \ska{\oR ^*\, \ov{x}}{\ov{n}} = \ska{\ov{x}}{\oR \,\ov{n}}  =
\beta ^{1/2} \,\tr _{M_1}(x(pn)^*) = \hfill\cr
      \beta ^{1/2} \,\tr _{M_1}(pxn^*) =\beta ^{1/2} \, \tr _{M_1}(xn^*) =
\beta 
^{1/2} \, \tr _{M_1}(E_0(E_1(x))\, n^*)\hfill }$$
for $n\in N$, as $E_0\circ E_1:M_1 \pfeil N$ is the conditional expectation 
onto $N$ associated
with the unique normalized trace $\tr _{M_1}$ of $M_1$.
Hence $\oR ^*\, \ov{x} =\beta ^{1/2}\, \ov{E_0(E_1(x))}$.
Let $R$ be defined as $\oR $ where $p$ is replaced by $p^{op}$. We compute
\begin{eqnarray*}
(R^* \times \B{1}_{\ov{\C{H} }})\circ (\B{1}_{\ov{\C{H} }} 
\times \oR )\,\, \ov{m} &=&
    \beta ^{1/2}(R^* \times \B{1}_{\ov{\C{H} }}) \,\,\ov{m} \tn \ov{p}
      \,\,\,  =\\
    \beta ^{3/2}\, (R^* \times \B{1}_{\ov{\C{H} }}) \,\sum _{j=1}^m
    \ov{m e_0e_1 v_je_0v_j^*} &=&
    \beta ^{1/2} \sum _{j=1}^m R^*\,\, \ov{me_0v_j} \tn \ov{v_j^*} 
  \,\,\, = \\
    \beta \, \sum _{j=1}^m \ov{E_0(E_1(me_0v_j))\, v_j^*} & =& 
 \qquad \qquad \qquad \hbox{(by Equation (~\ref{E4}))} \\
\sum _{j=1}^m \ov{(v_j\,E_0(v_j^*m^*))^* } &=&
    J \,(p.\ov{m^*})\, =\qquad  p^{op}.\ov{m} =\qquad \ov{m} 
            \end{eqnarray*}    
for every 
$m\in p^{op}L^2(M)\cap M$. (In line 2 we used the canonical isomorphisms from
$L^2(M_2)$ onto $L^2(M)\tn L^2(M_1)$ and $L^2(M_1)\tn L^2(M)$, see \cite{Sch3} 
or consider the proof of Theorem 4.1.)
By applying (f) with $p^{op}$ instead of $p$
and reversing the parts of $R$ and
$\oR $, we see that $R$ and $\oR $ are conjugation operators for $\C{H}$.
For $\ov{\B{1}}\in N$ we have
\begin{eqnarray*}
(\oR ^*\circ \oR) \,\ov{\B{1}} \,\,=\,\,\beta \,\ov{E_0(E_1(p))}&=&
\beta \, \tr _N(E_0(E_1(p))) \,\ov{\B{1}}\,\, =\,\,\beta \, \tr _{M_1}(p) 
\,\ov{\B{1}} =\\
c(\lambda (N), L^2(M))\cdot \tr _{M_1}(p) \,\ov{\B{1}} &=&
c(\lambda (N), pL^2(M))\,\ov{\B{1}} \end{eqnarray*}
(observe $\oR ^*\circ \oR\in (\iota ,\iota )=\comp $ for the second equation
and see Proposition 3.2.5 (e) in \cite{GHJ} for the last equation).
Just so we get 
$$ R^*\circ R=c(\lambda (N), p^{op}L^2(M)).$$ 
It is easy to see that $c(\lambda (N), \ov{\C{H}}) =c(\rho (N), \C{H})$.
If $\C{H}$ is irreducible, then
$$\sqrt[4]{ \frac{c(\lambda (N), \C{H})}{c(\rho (N), \C{H})}}
\cdot R \qquad\hbox{and} \qquad  
\sqrt[4]{\frac{c(\rho (N), \C{H})}{c(\lambda (N), \C{H})}}
\cdot \oR $$
are standard conjugation operators and we get formula (~\ref{dbim}).
 
\abs
We prove the Remarks (1) and (2): Let $_N\C{H}_N$ be an irreducible element of 
$\Ob \C{B}_{N\subset M}$.
Since $\lambda (N)' \cap \rho (N)' =\comp \B{1}$, the inclusion 
$\lambda (N)\subset \rho (N)'$ is 
extremal and the minimal index and the Jones' index agree 
(see \cite{P}, 1.2.5). We get
$$d(\C{H}) = \left( \frac{c(\lambda (N),\C{H})}{c(\rho (N)',\C{H})} 
\right) ^{1/2}   =
[\rho (N)':\lambda (N)]^{1/2} =([\rho (N)':\lambda (N)]_{min})^{1/2}. $$
So we obtain Remark (1) for irreducible bimodules.
As the minimal dimension
and $d$ are additive for direct sums (see \cite{Lo2} for the minimal 
dimension), we get Remark (1) for any bimodule of $\Ob \C{B}_{N\subset M}$. 

If $N\subset M$ is extremal $N\subset M_{2k-1}$ is
extremal for every $k\in \nat $ by 1.2.5 (iii) and (iv) in \cite{P}.
This implies $\tr _{M_{2k-1}}(p)=
\tr _{M_{2k-1}}(p^{op})$ for every minimal projection
$p\in N'\cap M_{2k-1}$ by 1.2.5 (i) in \cite{P}. We obtain
\refstepcounter{equation} \label{E12E1}
$$\displaylines{c(\rho (N), pL^2(M_{k-1}) )\, =\tr _{M_{2k-1}}(p)\, \cdot \,
c(\rho (N), L^2(M_{k-1}) )\, =\hfill \cr
\tr _{M_{2k-1}}(p)\, \cdot \, [M:N]^k\, =
\tr _{M_{2k-1}}(p^{op})\cdot [M:N]^k\, =\hfill\cr
c(\rho (N), p^{op}L^2(M_{k-1}) )\,=
c(\lambda (N), pL^2(M_{k-1}) ), \hfill (~\ref{E12E1})     }$$
and Remark (2) follows for every irreducible bimodule of 
$\Ob \C{B}_{N\subset M}$. Now the observation that the dimension $d$ and the
coupling constant are additive for direct sums completes the proof of
Remark (2).
\blacksquare

\abs
We consider the example $N\subset N\kreuz H$, where $H$ is 
a finite dimensional Hopf-$*$-algebra acting on $N$ by an outer action.
Observe that every irreducible subfactor $N\subset M$ with finite index and
depth 2 is isomorphic to a subfactor of that kind.
(This result was announced by A. Ocneanu, for complete proofs see
\cite{Da} or \cite{Sz}.) First we recall some definitions and results
for the reader's convenience (see \cite{Sz2} for example).
 
\subsection{Hopf-$*$-algebras and their actions on $\zws $-factors}
 (1) Let $N$ be a unital $*$-algebra and $H$ be a finite 
dimensional
Hopf-$*$-algebra. A bilinear map 
$$H\times N \pfeil N,\,  (a,x)\Pfeil \alpha (a)x,$$
is called an action of $H$ on $N$, if
\begin{itemize}
\item[(a)] $\alpha $ is a nondegenerate representation of the algebra $H$
on the linear space $N$,
\item[(b)] $\alpha (a)(xy) =\sum _{(a)} \alpha (a^{(1)})x \cdot 
\alpha (a^{(2)})y$ for $a \in H,\, x,y\in N$, where
$\Delta (a) =\sum _{(a)} a^{(1)}\otimes a^{(2)}$ (Sweedler notation), 
\item[(c)] $\alpha (a)\B{1}\,= \epsilon (a) \B{1}$ for $a\in H$ and
\item[(d)] $(\alpha (a)x)^*\, =\alpha (S(a)^*)x^*$ for $a\in H$ and $x\in N$.
\end{itemize}
(We use the convention $\alpha (a)x \cdot y:=( \alpha (a)x)\, y$.)

\abs
\noindent (2) The algebra
$$N^H:=\{x\in N:\, \alpha (a)x =\epsilon (a) x \hbox{\ for\ } a\in H\}$$
is called the fixed point algebra under the action $\alpha $.

\abs
\noindent (3) The crossed product $N\kreuz _{\alpha } H =N\kreuz
 H$ is a unital $*$-algebra that as a linear space coincides with the 
algebraic tensor product $N\otimes H$. The multiplication is given by
$(x\otimes a)\cdot (y\otimes b)= \sum _{(a)} x\cdot
\alpha (a^{(1)})y \otimes
a^{(2)}b$ and the $*$-operation by
$(x\otimes a)^*=\sum _{(a)} \alpha (a^{(1)*})x^* \otimes a^{(2)*}$.

\abs \noindent (4)
There is a unique faithful normalized trace $\tau $ (see \cite{Koo})
of $H$ called the
Haar trace such that 
$$(\tau \otimes id_H)(\Delta (a)) \, =\tau (a)\B{1} \,=
  (id_H \otimes \tau )(\Delta (a))$$
for any $a\in H$. If $p$ is a minimal central projection of $H$ such that
\linebreak $pH\cong \B{L}(\comp ^d)$, then $\tau (f)=d/ \dim H$
for every minimal projection $f$ of $pH$. One easily concludes that the
comultiplication $\Delta :H\pfeil H\otimes H$ is a unitary 
corepresentation with respect to
the inner product $(a,b)\Pfeil \tau (ab^*)$.

\abs
\noindent (5) If $N$ is a factor of type $\zws $ and $\alpha $ an action of 
the finite dimensional Hopf-$*$-algebra $H$ on $N$,
then $N^H$ as well as $N\kreuz H$ are von Neumann algebras (see \cite{Sz2}).
$N$ is embedded into $N\kreuz H$ by $x\Pfeil x\otimes \B{1}$.
$\tr_N \otimes \tau $ is a 
normalized finite trace of the von Neumann algebra
$N\kreuz H$ (see \cite{Sz2}, Proposition 2.7).

\abs
\noindent (6) An action $\alpha $ of $H$ on a 
factor $N$ of type $\zws $ is called 
outer
if $(N^H)' \cap N  =\comp \B{1}$ and if $\alpha $ is a faithful 
representation of $H$. If $\alpha $ is an outer action 
$N\kreuz H$ is a $\zws$-factor.
Furthermore, we have 
\begin{equation} N'\cap (N\kreuz H)=\comp \B{1}   \label{F13a} \end{equation} 
 (see \cite{Sz2}, 
Proposition 3.2).
According to (5) we may identify $L^2(N\kreuz H)$
with $L^2(N)\otimes H$ where $H$ is endowed with the inner product given by
the Haar trace $\tau $.

\abs
\noindent (7)
Let $H$ be a finite dimensional Hopf-$*$-algebra. If we use the reversed
comultiplication 
$$\Delta ^{cop}:H \pfeil H\otimes H,\, a\Pfeil \sum _{(a)}a^{(2)}\otimes
  a^{(1)},$$
instead of $\Delta $ on $H$ and do not change the remaining structure on $H$,
$H$ is a Hopf-$*$-algebra again (see \cite{Ka}, Corollary
III.3.5). We denote this Hopf-$*$-algebra by $H^{cop}$. 
(For general Hopf-algebras the antipode
$S$ has to be replaced by $S^{-1}$,
but $S=S^{-1}$ is satisfied for finite dimensional Hopf-$*$-algebras.) 
  
Similarly, we are able to introduce another Hopf-$*$-algebra structure
on $H$ by 
using the reversed multiplication $a\circ b:=ba$ (and leaving
the remaining structure unchanged). This Hopf-$*$-algebra
is called $H^{op}$. 
 
\abs
The connection between finite dimensional Hopf-$*$-algebras
and subfactors has also been studied by R. Longo (see \cite{Lo3}),
whose approach is based on index theory for infinite factors,
in particular on sector theory.

\begin{prop} Let $H$ be a finite dimensional Hopf-$*$-algebra and let $N$ be
a factor of type $\zws $ endowed with an outer action $\alpha $ of $H$.
\begin{enumerate}
\item[(i)] Let $\sigma :V_{\sigma }\pfeil V_{\sigma } \otimes H^{cop}$  
be a (small)
unitary finite dimensional corepresentation of $H^{cop}$. On the Hilbert
space $\C{H}_{\sigma }=L^2(N)\otimes V_{\sigma }$ there are a left action
$\lambda $ of $N$ determined by $\lambda (m)\, \ov{n} \otimes v \,=
\ov{mn} \otimes v$ and a right action 
$\rho $ of $N$ determined by
$$\rho (m)\, \ov{n} \otimes v\, = \sum _{(v)} \ov{n\cdot \alpha (v^{(2)})m} 
\otimes v^{(1)}$$ 
($m,n\in N,\, v\in V_{\sigma }$). Endowed with these actions
$\C{H}_{\sigma }$ is an $(N,N)$-bimodule belonging to the objects of the
$C^*$-tensor category $\C{B}_{N\subset N\kreuz H}$.
\item[(ii)] The following data determine an equivalence
$(F,(U_{\rho \sigma })_{\rho ,\sigma },J)$
of the $C^*$-tensor categories $\C{U}_{H^{cop}}$ and
$\C{B}_{N\subset N\kreuz H}$: \nopagebreak
\newl $F:\C{U}_{H^{cop}}\pfeil   
\C{B}_{N\subset N\kreuz H}$ is the $C^*$-tensor functor given by 
$F(\sigma ) =\C{H}_{\sigma }$, $F(T)=id_{L^2(N)}\otimes T$ for
$T\in (\rho ,\sigma )$ and $\rho ,\, \sigma \in \Ob \C{U}_{H^{cop}}$,
\newl $U_{\rho \sigma }:\C{H}_{\rho } \tn \C{H}_{\sigma }
\pfeil \C{H}_{\rho \otimes \sigma }$ is defined by 
\begin{equation} U_{\rho \sigma }\,
(\ov{n_1}\otimes v) \tn (\ov{n_2}\otimes w)=\sum _{(v)}
\ov{n_1\cdot \alpha (v^{(2)})n_2} \otimes v^{(1)} \otimes w \label{F13j} 
\end{equation}
for $n_1,\, n_2 \in N,\, v\in V_{\rho }$ and $w\in V_{\sigma }$,\newl 
$J:\C{H}_{\iota }=L^2(N) \otimes \comp \pfeil L^2(N)$ is defined by
$J(\ov{n}\otimes \gamma ) =\gamma\, \, \ov{n}$ for $\gamma \in \comp $ and
$n\in N$. \end{enumerate} \label{PH1} \end{prop} 

\B{Remark: }Let $G$ be a finite group and $H$ the Hopf-$*$-algebra 
$\Fun (G)$ of all complex valued functions on $G$. $N\kreuz H$ coincides with 
the usual crossed product $N\kreuz G$, where $G$ acts on $N$ by an outer
action in the usual meaning. The Hopf-$*$-algebras $\Fun (G)$ and
$\Fun (G)^{cop}$ are isomorphic (where the isomorphism is given by
$f\Pfeil \hat{f}, \hat{f}(g)=f(g^{-1})$). Hence the $C^*$-tensor category
$\C{B}_{N\subset N\kreuz G}$ is equivalent to the $C^*$-tensor category 
$\C{U}_G$ of all unitary finite dimensional representations of $G$.

\B{Proof:\ } (i)
Obviously $\lambda $ defines a left action of $N$.
For $m\in N$, $\rho (m)$ is well defined on the algebraic
tensor product $N\otimes V_{\sigma }$. It is easy to see
that $\rho (m)$ is continuous
on $N\otimes V_{\sigma }$ with respect to the Hilbert space norm.
Hence $\rho (m)$ has a
unique extension to a continuous linear operator
on $L^2(N) \otimes V_{\sigma }$. We easily get $\rho (\B{1}) =\B{1}$.
We  show
\begin{equation}
\rho (m_1m_2)= \rho (m_2)\,
 \rho (m_1) \qquad \hbox{for all $m_1,\, m_2\in N$} \label{Rmult} 
\end{equation}
and
\begin{equation} \rho (m^*)=\rho (m)^* \qquad \hbox{for $m\in N$,} 
\label{Rstern} \end{equation} 
which implies that $\rho $ is a right 
action of $N$ on $\C{H}_{\sigma }$. Obviously
$\rho $ and $\lambda $ commute such that $\C{H}_{\sigma }$ is actually an
$(N,N)$-bimodule. For $\C{H}_{\sigma }\in \Ob \C{B}_{N\subset N\kreuz H}$
see Part (ii) (b) of the proof.

The computation
$$\displaylines{
\rho (m_2)\, \rho (m_1) \, \ov{n}\otimes v \, =
\hfill\cr
\sum _{(v)}\sum _{(v^{(1)})} \ov{n\cdot \alpha (v^{(2)})m_1 \cdot
 \alpha (v^{(1)(2)}) m_2} \otimes v^{(1)(1)}\, =\hfill\cr
\sum _{(v)}\sum _{(v^{(2)})} \ov{n\cdot \alpha (v^{(2)(2)})m_1 \cdot 
\alpha (v^{(2)(1)})
m_2} \otimes v^{(1)}\, = \qquad \hbox{(*)}\hfill\cr
\sum _{(v)} \ov{n\cdot \alpha (v^{(2)})(m_1m_2) } \otimes v^{(1)}\,=\qquad
\qquad \qquad \rho (m_1m_2) \, \ov{n}\otimes v  \hfill      }$$
($n\in N$ and $v\in V_{\sigma }$) shows Equation (~\ref{Rmult}).
(Observe that the Sweedler notation for $H^{cop}$ is used in the
computation, but in (*)
the comultiplication of $H$ is needed.) 

For $n_1,\, n_2
\in N$ and $v,w\in V_{\sigma }$ we get
\begin{equation} \ska{\rho (m^*)\, \ov{n_1}\otimes v}{\ov{n_2}\otimes w} \, =
\sum _{(v)} \tr_N\Bigl(n_1 \cdot\alpha (v^{(2)})m^*\cdot n_2^*\Bigr)\, 
\ska{v^{(1)}}{w}
\label{EH5} \end{equation}
as well as
\begin{eqnarray}
\ska{\ov{n_1}\otimes v}{\rho (m)\, \ov{n_2}\otimes w} & =&
\sum _{(w)} \tr_N \Bigl(n_1\cdot (n_2 \cdot \alpha (w^{(2)})m)^*\Bigr)\,
\ska{v}{w^{(1)}}  
\nonumber \\
&=& \sum _{(w)}\tr _N\Bigl(n_2^*\cdot n_1\cdot 
\alpha (S(w^{(2)*}))m^* \Bigr)\, 
\ska{v}{w^{(1)}}.  \hspace{5mm} \label{EH6}      \end{eqnarray}
Since $\sigma $ is unitary, we know 
$$\sum _{(v)} \ska{v^{(1)}}{w}\, S(v^{(2)})
\, = \sum _{(w)} \ska{v}{w^{(1)}}\, w^{(2)*} \qquad 
 \hbox{for $v,w\in V_{\sigma }$} $$
(see \cite{Koo}, (1.42)). An application of the antipode $S$ yields    
\begin{equation} \sum _{(v)}\ska{v^{(1)}}{w}\, v^{(2)} \, =
 \sum _{(w)} \ska{v}{w^{(1)}}\, S(w^{(2)*}) \label{EH1} \end{equation}
(observe $S^2=id_H$), and the right sides of (~\ref{EH6}) and (~\ref{EH5})
coincide. So Equation (~\ref{Rstern}) has been shown.

\abs
\noindent (ii) (a) Let $\sigma $ and $\pi $ be unitary corepresentations of 
$H^{cop}$. We will prove that a continuous linear operator $S$ from 
$\C{H}_{\sigma }$
into $\C{H}_{\pi }$ is $(N,N)$-linear if and only if there is an
operator $T\in (\sigma ,\pi )$ such that $S=\B{1}_{L^2(N)}\otimes T$.

Obviously, $(\B{1} \otimes T)\, \lambda _{\sigma }(m)\, =      
\lambda _{\pi }(m)\, (\B{1}\otimes T)$ for every $T\in
(\sigma ,\pi )$ and $m\in N$, where $\lambda _{\sigma }$ denotes the left 
action of $N$ on $\C{H}_{\sigma }$. The corresponding equation for the right
action is shown by the following computation:
\begin{eqnarray*}
\rho _{\pi }(m)\, (\B{1}\otimes T)\, \ov{n} \otimes v &=&
\sum _{(Tv)} \ov{n\cdot \alpha ((Tv)^{(2)})m } \otimes (Tv)^{(1)} \,\,\, =\\
\sum _{(v)} \ov{n\cdot \alpha (v^{(2)})m} \otimes Tv^{(1)} &=&
(\B{1}\otimes T) \rho _{\sigma }(m) \, \ov{n}\otimes v
\qquad (n\in N,\, v\in V_{\sigma }).  \end{eqnarray*}

Now let $S$ be any continuous $(N,N)$-linear operator in $\C{H}_{\sigma }$.
Let $\lambda _0$ (resp. $\rho _0$) denote the canonical left action
(resp. right action) of $N$ on $L^2(N)$. The following
computation shows $\rho (m) =\rho _0(m) \otimes \B{1}$ for every
$m\in N^H$:
$$\rho (m)\, \ov{n}\otimes v \,= \sum _{(v)} \ov{n\cdot \epsilon (v^{(2)})m}
\otimes v^{(1)}\, = \ov{nm}\otimes v$$
for every $n\in N$ and $v\in V_{\sigma }$. Hence
$$\displaylines{
S\in \Bigl((\lambda _0(N) \cup \rho _0(N^H)) \otimes \comp \B{1}\Bigr)'\, \,= 
\Bigl(\lambda _0(N)' \cap \rho _0(N^H)'\Bigr)\otimes \B{L}(V_{\sigma })
 \,=\hfill\cr
\qquad \Bigl(\rho _0(N) \cap \rho _0(N^H)'\Bigr)\otimes \B{L}(V_{\sigma }) \,=
\comp \B{1} \otimes \B{L}(V_{\sigma }). \hfill      }$$
There is an operator $T\in \B{L}(V_{\sigma })$ such that $S=\B{1}\otimes T$.
$S\rho (m)\,\, \ov{\B{1}}\otimes v \,= \rho (m)S\, \, \ov{\B{1}}\otimes v $ 
for $m\in N$ and $v\in V_{\sigma }$ implies
$$\sum _{(v)} \ov{\alpha (v^{(2)})m } \otimes Tv^{(1)} \,=
\sum _{(Tv)} \ov{\alpha ((Tv)^{(2)})m } \otimes (Tv)^{(1)}.$$
Since $\alpha $ is a faithful action, we get
$$\sum _{(v)} Tv^{(1)}\otimes v^{(2)}\, =\sum _{(Tv)} (Tv)^{(1)} 
  \otimes (Tv)^{(2)}$$
for every $v\in V_{\sigma }$. It follows $T \in (\sigma ,\sigma )$.

At last we consider a continuous $(N,N)$-linear operator $S$ from 
$\C{H}_{\sigma }$ into 
$\C{H}_{\pi }$. $S$ may be regarded as an $(N,N)$-linear operator in
$\C{H}_{\sigma }\oplus \C{H}_{\pi }= \C{H}_{\sigma \oplus \pi }$, and
we are able to apply the considerations from above 
(with $\sigma \oplus \pi $ instead of $\sigma $) and find that there is
an operator $T\in (\sigma ,\pi )\subset (\sigma \oplus \pi ,\sigma 
\oplus \pi )$
such that $S=\B{1}\otimes T$ (observe Remark 1.3 (4)).
 \abs

\noindent (b) $\Delta ^{cop}$ is a unitary corepresentation of $H^{cop}$
with respect to the Haar trace $\tau $ of $H$. 
Obviously the 
the $(N,N)$-bimodule $\C{H}_{\Delta ^{cop}}$ is equal to the 
$(N,N)$-bimodule $L^2(N\kreuz H)$ (=$L^2(N) \otimes H$
as a Hilbert space).
From \cite{Koo} we conclude that every irreducible corepresentation
$\sigma $ of $H^{cop}$ is contained in $\Delta ^{cop}$ and that
$\C{H}_{\sigma }\in \Ob \C{B}_{N\subset N\kreuz H}$. Hence 
$F$ is a full and faithful $C^*$-functor. 
Let $\C{H}$ be an irreducible object of $\C{B}_{N\subset N\kreuz H}$.
As $N\subset N\kreuz H$ is a subfactor of depth 2, $\C{H}$ is equivalent
to a sub-bimodule of $_NL^2(N\kreuz H)_N \,=\C{H}_{\Delta ^{cop}}$.
Therefore there is an irreducible corepresentation $\sigma $ of 
$H^{cop}$ such that
$\C{H}$ is equivalent to $\C{H}_{\sigma }$.

\abs
\noindent (c) Let $\rho $ and $\sigma $ be unitary corepresentations
of $H^{cop}$. 
It is easy to see that 
$\xi \tn (\ov{n}\otimes w) \Pfeil (\xi \tn \ov{n}) \otimes w$
($\xi \in \C{H}_{\rho },\, n\in N,\, w\in V_{\sigma }$)
defines a unitary operator $W_{\rho \sigma }$ from 
$\C{H}_{\rho }\tn \C{H}_{\sigma }=\C{H}_{\rho } \tn
(L^2(N) \otimes V_{\sigma })$ onto  
$(\C{H}_{\rho }\tn
L^2(N)) \otimes V_{\sigma }$. Now we easily conclude that
$U_{\rho \sigma }:=(r_{\C{H}_{\rho }}\otimes id_{V_{\sigma }})\circ
 W_{\rho \sigma }$ 
is the unique continuous linear operator satisfying (~\ref{F13j}), 
where $r_{\C{H}_{\rho }}$ is the canonical unitary operator from
$\C{H}_{\rho }\tn L^2(N)$ onto $\C{H}_{\rho }$.
Obviously $U_{\rho \sigma }$ is 
unitary and left $N$-linear.

We get
$$\displaylines{
U_{\rho \sigma }\circ \rho (m)\, (\ov{n}\otimes v)\tn (\ov{\B{1}}\otimes w)\, =
U_{\rho \sigma }\sum _{(w)} (\ov{n}\otimes v)\tn 
\Bigl(\ov{\alpha (w^{(2)})m}\otimes w^{(1)}\Bigr)\, =
\hfill\cr
\sum _{(v),(w)} \ov{n \cdot \alpha (v^{(2)}w^{(2)})m } \otimes v^{(1)}
\otimes w^{(1)} \,=\hfill\cr
\rho (m) \, \ov{n} \otimes v \otimes w =\rho (m)\circ  U_{\rho \sigma }\, 
(\ov{n}\otimes v)\tn (\ov{\B{1}}\otimes w) \hfill        }$$
for all $m,n\in N,\, v\in V_{\rho },\, w\in V_{\sigma }$. Since
$$\spann \{(\ov{n} \otimes v)\tn (\ov{\B{1}}\otimes w):\, v,w\in V_{\rho },
\, n\in N\}$$
is dense in $\C{H}_{\rho }\tn \C{H}_{\sigma }$, we 
have verified that $U_{\rho \sigma }$ is right $N$-linear.

Routine arguments show that the relations (a), (b) and (c) of
Definition ~\ref{DCequ} (ii) are satisfied. (The computation for (b) is easier,
if one uses the density argument from above.)
Hence $F$ is a $C^*$-tensor equivalence. 
 \blacksquare 

\section{Subfactors defined by $C^*$-tensor categories}
Let a finite $C^*$-tensor category $\C{C}$ be given. The goal of this
Section is to construct a subfactor belonging to an object $\sigma $ of
$\C{C}$ and to discuss the properties of this subfactor.

For an object $\rho $ of $\C{C}$ $(R_{\rho },\, \ov{R}_{\rho })$ denotes a
standard pair of conjugation operators for $\rho $ and $\ov{\rho }$.

\subsection {Some observations \label{obs}}
Let $\rho ,\, \sigma $ and $\tau $ be objects of $\C{C}$.
\newl (1) There is a distinguished
normalized faithful trace $\tr _{\rho }$ of the finite dimensional
$C^*$-algebra $(\rho ,\rho )$. $\tr _{\rho }$ is uniquely determined
by its values on minimal projections. Let $E$ be a minimal projection
and $\pi $ a subobject of $\rho $ corresponding to $E$. Then
\begin{equation} \tr _{\rho }(E):= \frac{d(\pi )}{d(\rho )}.   \label{E21E1} 
\end{equation}
Conversely this relation determines a faithful trace, as equivalent 
minimal projections have the same positive number.

Obviously the trace $\tr _{\rho }$ is faithful. Theorem ~\ref{Dim} (i) shows
that $\tr _{\rho }$ is normalized. If
$A,B \in (\rho ,\sigma )$ then
\begin{equation} d(\rho ) \, \tr_{\rho }(B^*A) =d(\sigma ) \, 
\tr _{\sigma }(AB^*). 
 \label{zwtr} \end{equation}
(Observe that $d(\rho ) \, \tr _{\rho }$ is the trace introduced
in \cite{LoRo} after Lemma 3.7; thus (~\ref{zwtr}) is a consequence of Lemma  
3.7 in \cite{LoRo}.)

\abs
\noindent (2) The
$*$-algebras $(\rho ,\rho )$ and $(\sigma ,\sigma )$ are embedded into the
$*$-algebra $(\rho \sigma ,\rho \sigma )$ by
$T\Pfeil T\times \B{1}_{\sigma }$ respectively 
$T\Pfeil \B{1}_{\rho } \times T$.  From
\begin{eqnarray*} \B{1}_{\rho} \times \oR ^*_{\sigma } \,\circ \,
  (T\times \B{1}_{\sigma })\times \B{1}_{\osigma } \,\circ \,
  \B{1}_{\rho} \times \oR _{\sigma } \,&=& 
  (\oR _{\sigma }^*\circ \oR _{\sigma }) \cdot T 
\qquad \hbox{and} \\
  R_{\rho }^* \times \B{1}_{\sigma } \, \circ \,
 \B{1}_{\orho } \times (\B{1}_{\rho } \times T) \,\circ \,   
 R_{\rho } \times \B{1}_{\sigma } \,&=&
 (R_{\rho }^* \circ R_{\rho })\cdot T 
\end{eqnarray*}
we conclude that the embeddings are injective.
R. Longo and J.E. Roberts proved (see Corallary 3.10 in \cite{LoRo}) 
\[ d(\rho )\, \tr_{\rho }(S_1T_1^*)\cdot d(\sigma )\, \tr_{\sigma }(S_2T_2^*)
\,= d(\rho \sigma ) \, \tr _{\rho \sigma }((S_1\times S_2)(T_1\times T_2)^*)
\]
for all $S_1,\, T_1\in (\rho ,\rho ),\, S_2,\, T_2\in (\sigma ,\sigma )$.
Using $d(\rho \sigma )=d(\rho )d(\sigma )$ we get 
\begin{equation} \tr _{\rho \sigma }\rest (\rho ,\rho ) = \tr _{\rho } 
\qquad \hbox{and}
\qquad \tr_{\rho \sigma }\rest (\sigma ,\sigma ) = \tr _{\sigma }.
\label{Trrest}\end{equation}

\abs
\noindent (3) The linear maps
$$\displaylines{
 \Phi ^{\rho }_{\sigma }: (\rho \sigma ,\rho \sigma )\pfeil 
   (\sigma ,\sigma ),\,  X\Pfeil \frac{1}{d(\rho )}   
   \, R_{\rho}^*\times \B{1}_{\sigma } \, \circ  \,
   \B{1}_{\rho } \times X \, \circ \, R_{\rho }\times \B{1}_{\sigma }, 
   \hfill \hbox{and} \cr
  \Psi ^{\rho }_{\sigma }: (\sigma \rho ,\sigma \rho )\pfeil 
   (\sigma ,\sigma ),\, X\Pfeil \frac{1}{d(\rho )}   
   \, \B{1}_{\sigma } \times \oR_{\rho}^* \, \circ  \,
   X \times \B{1}_{\rho }\, \circ \, \B{1}_{\sigma } \times \oR_{\rho }, 
    \hfill }$$
are conditional expectations.
(In \cite{LoRo} R. Longo and J.E. Roberts used the more 
general concept
of left and right inverses, which we do not need here.) 
$\Phi ^{\rho }_{\sigma }$ and $\Psi ^{\rho }_{\sigma }$ 
 do not depend on the
choice of $R_{\rho }$ and $\oR _{\rho }$ (see Lemma 3.3 in \cite{LoRo}).
R. Longo and J.E. Roberts showed that
$\Phi _{\iota }^{\rho }:(\rho ,\rho ) \pfeil (\iota ,\iota )=\comp $ and
$\Psi _{\iota }^{\rho }:(\rho ,\rho ) \pfeil (\iota ,\iota )=\comp $ 
are just $\tr _{\rho }$ (see Lemma 3.3 in \cite{LoRo}).
We have   
\begin{eqnarray*} 
\Phi _{\tau }^{\sigma } \circ \Phi _{\sigma \tau }^{\rho } &=& 
\Phi _{\tau }^{\rho \sigma } \qquad   \hbox{as well as} \\
\Psi _{\tau }^{\sigma } \circ \Psi _{\tau \sigma }^{\rho } &=& 
\Psi _{\tau }^{\sigma \rho }. \end{eqnarray*}
Using Equation (~\ref{conpr}) we get 
$$\displaylines{\Phi _{\tau }^{\rho \sigma }(X) =
\frac{1}{d(\rho \sigma )} 
\Bigl(R_{\sigma }^*\circ \B{1}_{\osigma }\times R_{\rho }^*\times
\B{1}_{\sigma }\Bigr) \times \B{1}_{\tau } \,\circ \,
 \B{1}_{\osigma \orho }   \times X \,\circ \,
 \Bigl(\B{1}_{\osigma }\times R_{\rho }\times \B{1}_{\sigma } \circ R_{\sigma }
\Bigr)
 \times \B{1}_{\tau } = \hfill\cr
 \frac{1}{d(\rho )d(\sigma )}\, R_{\sigma }^* \times \B{1}_{\tau } \,\circ \,
   \B{1}_{\osigma }\times \Bigl(R_{\rho }^* \times \B{1}_{\sigma \tau } \circ 
   \B{1}_{\orho }\times X \circ R_{\rho } \times \B{1}_{\sigma \tau }\Bigr) \,
   \circ \, R_{\sigma }\times \B{1}_{\tau} =\hfill \cr
   \hfill =\Phi ^{\sigma }_{\tau }(\Phi ^{\rho }_{\sigma \tau }(X)) }$$
for any $X\in (\rho \sigma \tau ,\rho \sigma \tau )$. The same method also 
shows
the second relation. We note
\begin{equation} \tr _{\sigma } \circ \Phi ^{\rho }_{\sigma } = 
\tr _{\rho \sigma }
 \qquad \hbox{and} \qquad \tr _{\sigma } \circ \Psi ^{\rho }_{\sigma } =
\tr _{\sigma \rho }. \label{trexp} \end{equation}
as special cases. Equation (~\ref{trexp}) 
means that 
$\Phi ^{\rho }_{\sigma }$ (resp. $\Psi ^{\rho }_{\sigma }$) is a conditional
expectation corresponding to the trace $\tr _{\rho \sigma }$ (resp.
$\tr _{\sigma \rho }$).  

\abs     
\noindent (4) \begin{equation}  \begin{array}{ccc}
(\rho \sigma , \rho \sigma ) &\subset &(\rho \sigma \tau , 
  \rho \sigma \tau )\\
 \cup  & &\cup \\
(\sigma , \sigma ) & \subset &(\sigma \tau , \sigma \tau )
\end{array} \label{commC} \end{equation}
is a commuting square of finite dimensional von Neumann algebras
(in the meaning of \cite{GHJ}, Section 4.2)
 with respect to the trace $\tr _{\rho \sigma \tau }$,
as
$\Phi ^{\rho }_{\sigma \tau }(X\times \B{1}_{\tau }) =
\Phi ^{\rho }_{\sigma }(X) \times \B{1}_{\tau }$ holds for
any $X\in (\rho \sigma , \rho \sigma )$.

\abs
\noindent (5) The Bratteli diagram for the inclusion $(\rho , \rho )\subset 
(\rho \sigma , \rho \sigma )$ is described as follows:
\newl The vertices in the lower line (resp. upper line)
are in one to one 
correspondence to the simple direct summands of $(\rho ,\rho )$ 
(resp. $(\rho \sigma ,\rho \sigma )$) and are labelled by 
the equivalence
classes $[\phi ]$ (resp. $[\psi ]$)  of the irreducible subobjects 
$\phi $ (resp. $\psi $) of $\rho $ (resp. 
$\rho \sigma $). 
The number of edges between the vertices belonging
to $[\phi ]$ and to $[\psi ]$ is the dimension of $(\psi ,\phi \sigma )$,
i.e. the number of times that $\psi $ is contained in $\phi \sigma $.
We leave the easy proof to the reader.

One can deal with the inclusion 
$(\rho ,\rho )\subset (\sigma \rho ,\sigma \rho )$ in the same way,
one merely has to replace $\rho \sigma $ by $\sigma \rho $ and
$\phi \sigma $ by $\sigma \phi $.

Let $[[\C{C},\rho ]]:= \{[\phi ]\in [[\C{C}]]:\phi \leq \rho \}$.
Obviously the Bratteli diagrams for $(\rho ,\rho ) \subset (\rho \sigma ,
\rho \sigma )$ and $(\tau ,\tau )\subset (\tau \sigma , \tau \sigma )$ are 
the same if $[[\C{C},\rho ]] =[[\C{C}, \tau ]]$.   

\pagebreak
\begin{lem} Let us consider the inclusion
$$A:=(\rho ,\rho ) \subset B:=(\rho \sigma , \rho \sigma ) \subset
C:=(\rho \sigma \osigma , \rho \sigma \osigma )$$
of finite dimensional $*$-algebras.
Moreover, let $B(A,B)$ be the basic construction for $A\subset B$
and for the trace $\tr _{\rho \sigma }$, let $e$ be the Jones projection
of $B(A,B)$, and let
$$ f:=\frac{1}{d(\sigma )} \B{1}_{\rho } \times \Bigl(\oR _{\sigma }\circ
\oR _{\sigma }^*\Bigr) \in C. $$
$D:= \hbox{\rm span\ } BfB $ is a (two-sided) ideal of $C$, and there is an
isomorphism $\alpha $ from the basic construction $B(A,B)$ 
onto $D$ such that
$\alpha \rest B = id _B$ and $\alpha (e)=f$.
If every irreducible subobject of $\rho \sigma \osigma $  is
a subobject of $\rho $ then $D=C$.

The trace $\tr _{\rho \sigma }$ satisfies the Markov relation
\begin{equation} \tr _{\rho \sigma \osigma }\Bigl((X\times 
\B{1}_{\osigma })\circ f\Bigr) =
\frac{1}{d(\sigma )^2} \,\tr _{\rho \sigma }(X) \label{Mark} \end{equation}
for $X\in (\rho \sigma ,\rho \sigma )$.

The Lemma is also true for $A=(\rho ,\rho )$, $B=(\sigma \rho ,\sigma \rho )$,
$C=(\osigma \sigma \rho ,\osigma \sigma \rho )$, and 
$f=\frac{1}{d(\sigma )} (R_{\sigma }\circ R_{\sigma }^*) \times \B{1}_{\rho }$.
The Markov relation for this case is 
$$\tr _{\osigma \sigma \rho }\Bigr((\B{1}_{\osigma }\times X)\circ f\Bigl)=
 \frac{1}{d(\sigma )^2} \, \tr _{\sigma \rho }(X).$$   \label{Lfund}  \end{lem}

\B{Proof:\ } We will only prove the first case.
$$x\in (\rho ,\rho )\Pfeil xf =x \times \frac{1}
{d(\sigma )}(\oR _{\sigma }\circ \oR _{\sigma }^*) \in 
(\rho \sigma \osigma ,\rho
\sigma \osigma )$$
is an injective homomorphism.
For $b\in B$ we get $fbf =\Psi ^{\sigma }_{\rho }(b) f$, as
the following computation using the interchange law shows:
$$\displaylines{
\Psi ^{\sigma }_{\rho }(b)\, f=\frac{1}
{d(\sigma )^2}\, (\B{1}_{\rho }\times \oR _{\sigma }^* \,\circ b \times     
\B{1}_{\osigma }\circ \B{1}_{\rho }\times \oR _{\sigma }) \times
\B{1}_{\sigma \osigma } \, \circ \, \B{1}_{\rho }    
\times (\oR _{\sigma }\circ \oR _{\sigma }^*) =\hfill\cr
\frac{1}
{d(\sigma )^2} \, \B{1}_{\rho }\times \oR_{\sigma } \circ 
(\B{1}_{\rho }\times \oR _{\sigma }^* \,\circ b \times     
\B{1}_{\osigma }\circ \B{1}_{\rho }\times \oR _{\sigma }) \circ
\B{1}_{\rho }\times \oR _{\sigma }^*= \hfill 
 fbf.    }$$
Now Proposition 2.6.9 in \cite{GHJ} implies that there is an isomorphism
$\alpha $ having the desired properties, and that $D$ is an ideal of the
$*$-algebra $\ska{B}{f}$ generated by $B$ and $f$. We have to show
that $D$ is an ideal of $C$.         

The subobject of $\sigma \osigma $ corresponding to the projection
$f_0:=\frac{1}
{d(\sigma )}\oR _{\sigma }\circ \oR _{\sigma }^*\in $\linebreak
$(\sigma \osigma ,
\sigma \osigma )$ is the unit object $\iota $.
In particular every irreducible
subobject of $\rho \iota =\rho $ is also a subobject of 
$\rho \sigma \osigma $. Let $I$ be the ideal of $C$ consisting of those
simple direct summands of $C$ for which the associated irreducible subobjects 
of
$\rho \sigma \osigma $ are also subobjects of $\rho $. We will
show $I=D$.

Let $\{q_1,q_2,\ldots ,q_n\}$ be a maximal set of pairwise non equivalent 
minimal projections of $A$ and let $\{\pi _1,\ldots ,\pi _n\}$ 
be the associated
irreducible objects of $\C{C}$. The projection $q_jf =q_j\times f_0$ 
(regarded as an element
of $C$) corresponds to the same irreducible object $\pi _j$ of $\C{C}$ as the
projection $q_j$
(regarded as an element of $A$), in particular $q_jf$ is a minimal projection
of $C$.  
Let $I_j$ (resp. $D_j$)
be the simple direct summand of $I$ (resp. $D$) containing
$q_jf$. Every minimal projection $p$ of $D_j$ is equivalent to $q_jf$ in
$C$, too , and belongs to the direct summand $I_j$ of $I$. So $D_j$ is 
contained in
$I_j$ for $j=1,\ldots ,n$. As $D=\oplus _{j=1}^n D_j$, $D$ is a subalgebra
of $I$. On the other hand the inclusion matrices for $B\subset D$ and
$B\subset I$ are the same:

Let $B=\oplus _{k=1}^m B_k$ be the decomposition of $B$ into a direct
sum of minimal simple ideals and let $\tau _k$ be the irreducible object of 
$\C{C}$
belonging to $B_k$. The vertices corresponding to $B_k$ and $I_j$ are 
connected by $\dim \, (\pi _j, \tau _k \osigma )$ edges. Since $D$ is the basic
construction for $A\subset B$, $B_k$ and $D_j$ are connected by
$\dim \, (\tau _k, \pi _j \sigma )$ edges
(see \cite{GHJ}, Proposition 2.4.1 (b)).
The Frobenius reciprocity law
(see Lemma ~\ref{LFrob}) yields
$\dim \, (\tau _k, \pi _j \sigma )= \dim \, (\pi _j, \tau _k \osigma )$. 

So $\dim D =\dim I$, and $D=I$ follows.
It remains to verify the Markov relation:
$$\displaylines{
\tr _{\rho \sigma \osigma }\, \Bigl((X\times \B{1}_{\osigma }) \,\cdot
   \frac{1}{d(\sigma )} (\B{1}_{\rho } \times (\oR _{\sigma }
     \circ \oR _{\sigma }^*))
        \Bigr)=
  \qquad \hbox{(by Equation (~\ref{zwtr}))} \hfill\cr
\frac{1}{d(\sigma )^3}\, \tr _{\rho } (\B{1}_{\rho } \times 
 \oR _{\sigma }^* \circ X\times \B{1}_{\osigma } \circ \B{1}_{\rho } \times 
 \oR _{\sigma }) \,= \hfill \cr
\frac{1}{d(\sigma )^2} \,\tr _{\rho  }(\Phi ^{\sigma }_{\rho }(X)) \,=
\qquad \frac{1}{d(\sigma )^2} \, \tr _{\rho \sigma }(X)
 \qquad \hspace{8mm} \hbox{for $X\in (\rho \sigma ,\rho \sigma )$.}  
\blacksquared }$$
  
\subsection {The subfactors \label{SSf}}
For an object $\sigma $ of $\C{C}$ we introduce the following notation: 
$$\displaylines{ \sigma (0)=\iota ,\,\, \sigma (1)=\sigma ,\,\,
\sigma (2) =\sigma \osigma,\,\, \sigma (3)=\sigma \osigma \sigma ,\,\,
\sigma (4) = \sigma \osigma \sigma \osigma , \ldots \hfill\hbox{and}\cr
\ov{\sigma (0)}=\iota ,\,\, \ov{\sigma (1)}=\osigma ,\,\,
\ov{\sigma (2)} =\sigma \osigma,\,\, \ovs{3}=\osigma \sigma \osigma ,\,\,
\ovs{4} = \sigma \osigma \sigma \osigma ,\, \ldots . \hfill }$$
We will consider the tower
\begin{equation} \hspace*{-8mm} \begin{array}{ccccccccc}
(\sigma ,\sigma ) & \subset &( \osigma \, \sigma   ,\, \osigma\, \sigma )&
\subset
 & (\ovs{2}\, \sigma  ,\, \ovs{2} \, \sigma )   &\subset &
( \ovs{3}\, \sigma , \ovs{3}\, \sigma )
& \subset &\ldots \\
\cup & &\cup & &\cup & &\cup & & \\
(\iota , \iota )& \subset &(\osigma ,\,\osigma )& \subset &
(\ovs{2},\, \ovs{2})  &\subset &  (\ovs{3},\, \ovs{3})
&\subset &\ldots .
\end{array} \label{E21E2} \end{equation}
of finite dimensional $*$-algebras.
In the following we abbreviate \newline
$\Bigl(\ovs{n-1},\, \ovs{n-1}\Bigr)$ to
$A^n$ and $\Bigl(\ovs{n-1}\, \sigma ,
\, \ovs{n-1}\, \sigma \Bigr)$ to
$B^n$. There is a trace $\tr $ of the $*$-algebra $B^{\infty }:=\bigcup _{n=1}^
{\infty } B^n$ 
given by $\tr \rest B^n  =\tr _{
  \ovs{n-1} \, \sigma }$.
Observation ~\ref{obs} (1) and (2) show that $\tr $ is well defined and 
faithful and that $\tr \rest A^n =\tr _{\, \ovs{n-1}}$.

If $d(\sigma )>1$, the subsequent considerations will show that
 the tower (~\ref{E21E2}) (endowed with 
the trace $\tr $) fulfils the periodicity assumptions used in H. Wenzl's
subfactor construction (Theorem 1.5 of \cite{We}). In particular this implies
the following: Let $\pi $ be the GNS representation
of the state $\tr $ of the $*$-algebra $B^{\infty }$ and let
$A=\pi (A^{\infty })''$ (with $A^{\infty }=\bigcup _{n=1}^{\infty }A^n$)
and $B=\pi (B^{\infty })''$. Then $A$ and $B$ are 
$\zws $-factors.

According to Observation ~\ref{obs} (4), the squares appearing in the tower
(~\ref{E21E2}) commute.
Observe that every subobject of $\ov{\sigma (n-1)}$ is also a subobject
of $\ov{\sigma (n+1)}$.
As the $C^*$-tensor category $\C{C}$ is finite, Observation ~\ref{obs} (5)
implies that the Bratteli diagram for $A^{n-1}\subset A^n$ is the same as that
for $A^{n+1}\subset A^{n+2}$ if $n$ is sufficiently large. 
The same holds for $(B^n)_{n\in \nat }$.
Using Observation ~\ref{obs} (5) again, we see
that the inclusion matrices for $A^n \subset B^n$ and $A^{n+2}\subset B^{n+2}$
coincide for a sufficiently large $n$.

By induction, we see that the Bratteli diagram for $A^n \subset A^{n+1}$ is
connected. (That means: For any two vertices of the graph there is a sequence
of edges connecting these vertices.) The case $n=0$ is obvious. The step
$n-1 \rightarrow n$ is a conclusion from the proof of Lemma ~\ref{Lfund}.
The Bratteli diagram for $A^n \subset A^{n+1}$ is the mirror image of the 
Bratteli diagram for $A^{n-1} \subset A^n$, where new vertices for $A^{n+1}$  
are added and connected with some vertices of $A^n$ by new 
edges. Now one easily concludes that the inclusion matrix for
$A^n\subset A^{n+2}$ is primitive for every $n\in \nat $. 

We have to exclude the case $d(\sigma )=1$. If $d(\sigma )>1$ we get 
$B^n\subsetneqq B^{n+2}$ for every $n\in\nat $.

\abs
For $m\in \nat $ we consider the tower
$$\displaylines{ B_m^1 =\Bigl(\sigma (m+1),\, \sigma (m+1)\Bigr) \subset 
\hfill\cr 
\qquad \subset  B_m^2 =\Bigl(\osigma\, \sigma (m+1) ,\,  \osigma \, \sigma (m+1)\Bigr) \subset 
\hfill\cr
\qquad \qquad \subset B_m^3 =
   \Bigl(\ovs{2}\, \sigma (m+1),\, \ovs{2} \, \sigma (m+1)\Bigr) \subset \hfill\cr
\qquad \qquad \qquad \subset 
B_m^4 =\Bigl(\ovs{3}\, \sigma (m+1),\ovs{3}\, \sigma (m+1)\Bigr) 
\subset \ldots .\hfill } $$
Let $B_m$ be the unique $\zws $-factor containing 
$B_m^{\infty } = \bigcup _{n=1}^ {\infty } B_m^n$ as an ultra-weakly dense
subalgebra. $B_{m-1}$ is canonically embedded into $B_m$.    
 
\begin{thm} Let $d(\sigma )>1$. \begin{enumerate}
\item[(i)] $[B:A] =d(\sigma )^2$.
\item[(ii)] There is an isomorphism $j$ from the basic construction
$B(A,B)$ for \linebreak $A\subset B$ onto the $\zws $-factor 
$B_1$ such that $j\rest B =id_B$ and
$j(e_0)=f_0$, where $e_0$ is the Jones projection in $B(A,B)$ and
$f_0= \frac{1}{d(\sigma )} \oR _{\sigma }\circ \oR _{\sigma }^* \in B_1^1 =
( \sigma \osigma , \sigma \osigma ) \subset B_1$.
\item[(iii)] The Jones tower for $A\subset B$ may be identified with
$$B_{-1}=A \subset B_0 =B \subset B_1 \subset B_2 \subset B_3 \subset
 \ldots .$$ The Jones projection $f_m\in B_{m+1}$ for
$B_{m-1}\subset B_m$ is \newl
$f_m= \frac{1}{d(\sigma )} 
 \B{1}_{\sigma (m)}\times \Bigl(\oR _{\sigma }\circ \oR _{\sigma }^*\Bigr)
 \in B_{m+1}^1$ for $m$ even and 
\newl $f_m= \frac{1}{d(\sigma )}\B{1}_{\sigma (m)} 
\times \Bigl(R_{\sigma }\circ R_{\sigma }^*\Bigr)
 \in B_{m+1}^1$ for $m$ odd.
\item[(iv)] The subfactor $A\subset B$ has finite depth. The standard 
invariant 
\[ \begin{array}{rcccccccc}
\comp \B{1}=A'\cap A &\subset &A' \cap B & \subset &A'\cap B_1&\subset & A'\cap B_2 &\subset & 
\ldots \\
 & &\cup & &\cup & &\cup & & \\
\comp \B{1}& =&B'\cap B & \subset &B' \cap B_1& \subset &B' \cap B_2 &\subset
&\ldots 
\end{array} \]  
of $A\subset B$ is equal to
$$ \hspace*{-9mm}\begin{array}{rcccccccccc}
(\iota ,\iota ) &\subset &(\sigma ,\sigma ) & \subset &(\sigma \osigma , 
\sigma \osigma )
&\subset & (\sigma \osigma \sigma , 
\sigma \osigma \sigma ) &\subset & (\sigma \osigma \sigma \osigma ,
\sigma \osigma \sigma \osigma ) & \subset 
\ldots \\
 & &\cup & &\cup & &\cup & &\cup  \\
& &(\iota , \iota ) & \subset & (\osigma ,\osigma )& \subset &(\osigma \sigma ,
\osigma \sigma ) &\subset &(\osigma \sigma \osigma , \osigma \sigma \osigma )&
\subset
\ldots .
\end{array}   $$
Especially the subfactor $A\subset B$ is irreducible
(i.e. $A'\cap B =\comp \B{1}$) if and only if
$\sigma $ is an irreducible object. \end{enumerate}
\label{TUnt} \end{thm} 

One easily shows the following

\subsection{Remarks:}
\begin{enumerate}
\item[(1)] If one chooses another conjugate $\tilde{\sigma }$ of $\sigma $ instead of 
$\osigma $ the tower (~\ref{E21E2}) is isomorphic to the tower formed with
$\tilde{\sigma }$ instead of $\sigma $, as an easy consideration shows.
Hence the subfactor does not depend on the choice of the conjugate $\osigma $.

\item[(2)] If $\C{C}$ and $\C{D}$ are finite $C^*$-tensor categories and 
$F$ is an equivalence of the
$C^*$-tensor categories $\C{C}$ and $\C{D}$ then the subfactor
constructed with the object $\sigma $ of $\C{C} $ is isomorphic to the
subfactor constructed with the object $F(\sigma )$ of $\C{D}$.
\end{enumerate}
 
\abs
\B{Proof of Theorem ~\ref{TUnt}:\ }
(i) We consider the inclusion
 $$A^{2n+1} =\Bigl(\ovs{2n}, \ovs{2n}\Bigr) \subset 
     B^{2n+1} =\Bigr(\ovs{2n} \, \sigma , \ovs{2n} \, \sigma \Bigl)$$
for a sufficiently 
large $n$. The object $\ovs{2n} \, \sigma \osigma =(\sigma \osigma
)^{n+1}$ contains the
same irreducible subobjects as $\ovs{2n}$, 
hence $B_1^{2n+1}$ is the basic construction
for \linebreak $A^{2n+1} \subset B^{2n+1}$ according to Lemma ~\ref{Lfund} and
$\tr _{\ovs{2n} \, \sigma }$ is a faithful Markov trace of modulus
$d(\sigma )^2$ for this inclusion. Theorem 1.5 in \cite{We}
(along with Proposition 2.7.2 in \cite{GHJ}) yields $[B:A]=d(\sigma )^2$.

\abs
\noindent (ii) By using \cite{PP}, Corollary 1.8 we will show 
that $A\subset B$ is the basic construction
downwards: We get
$[B_1:B] =d(\sigma )^2$ as in the proof for $A\subset B$. Lemma ~\ref{Lfund}
implies
\begin{equation} d(\sigma )^2 \, \tr _{B_1}(xf_0) =  \tr _B(x)  
\label{gmar}\end{equation}
for every $x\in B$. 
Therefore
$$ \frac{1}{[B_1:B]} \,\tr _B(x) = \tr _{B_1}(xf_0) =\tr _B(E_B(xf_0)) =
  \tr _B(x E_B(f_0))$$ for every $x\in B$ and $E_B(f_0)= [B_1:B]^{-1} \B{1}$.
Hence there is an isomporphism $\iota $ from the basic construction
$B(P,B)$ for the subfactor $P:=\{f_0\}'\cap B_1 \subset B$ onto 
the $\zws $-factor $B_1$ such that 
$\iota\rest B = id _B$ and $\iota ^{-1}(f_0)$ is the Jones projection.
Obviously $A$ is contained in $P$ 
and  $[B:P] =[B_1:B] =$ \linebreak $[B:A] <\infty $. Hence  $A=P$,
and we get the assertion. 

\abs
\noindent (iii) follows by applying the argument from (ii) repeatedly.

\abs
\noindent (iv) Obviously, for every $m\in \nat $
$B_m^1= (\sigma (m+1), \sigma (m+1))$
is contained in $A'\cap B_m$.
Let $n$ be a sufficiently large odd number,
and let $q$ be a minimal projection of 
$A^{n}=\nobreak \Bigl(\ovs{n-1}, \ovs{n-1}\Bigr)$ 
which corresponds to the identity object
$\iota $ of $\C{C}$. 
The subobject of $\ovs{n-1}\sigma (m+1)$ corresponding to the projection
$p=q\times \B{1}_{\sigma (m+1)}$ is $\iota \cdot \sigma (m+1) =\sigma (m+1)$,
hence
$$p ((A^n)'\cap B^n_m) = (p A^np)'\cap pB^n_mp =(\comp p)'\cap pB^n_mp
 =pB^n_mp$$
is isomorphic
to $(\iota \sigma (m+1), \iota \sigma (m+1))$ (compare Remark ~\ref{Resds} 
(1) and (5)). Theorem 1.6 in \cite{We} implies
$ \dim \, A'\cap B_m \leq \dim \, (\sigma (m+1), \sigma (m+1))
 =\dim B_m^1$. So $B_m^1 =A'\cap B_m$ has been shown. 
The same method yields
$$B'\cap B_m = (\iota \,\osigma \,\sigma (m-1), \iota \,\osigma \,\sigma (m-1) )  \subset B_m^1$$
for $m\geq 1$. 

Since the number of direct summands in the sequence $(A'\cap B_m)_m$
is bounded, the inclusion $A\subset B$ has
finite depth. \blacksquare

\subsection{Example \label{X21}}
We regard the subfactor $A\subset B$ if $\C{C}$ is the finite $C^*$-tensor
category $\C{B}_{N\subset M}$ for a subfactor $N\subsetneqq M$ of finite depth
and if $\sigma $ is the bimodule $_NL^2(M)_N$.

For an $(N,N)$-bimodule $\C{H}$,
$\C{L}_{-,N}(\C{H})$ denotes the von Neumann algebra of all operators of
$\B{L}(\C{H})$ commuting with $\rho (N)$ and
$\C{L}_{N,N}(\C{H})$ the von Neumann algebra of all operators commuting 
with $\rho (N)$ and $\lambda (N)$.

We have $\sigma =\osigma $ and
$(\sigma ^k,\sigma ^k) =\C{L}_{N,N}(L^2(M)^{\tn ^k})$. 
Corollary 2.3 in \cite{Sch3} tells us that for every $n>1$
there is an isomorphism 
$J_n:M_{2n-1}\pfeil \C{L}_{-,N}(L^2(M)^{\otimes _N^n})$
such that
\begin{eqnarray*}
J_n(m) &=& J_k(m)\tn id_{L^2(M)^{\otimes _N^{n-k}}} \qquad \hbox{for
$m\in M_{2k-1}$,}\\
J_n(N'\cap M_{2k-1})&=&\C{L}_{N,N}(L^2(M)^{\otimes _N^k})\tn  
\comp id_{L^2(M)^{\otimes _N^{n-k}}}, \qquad \hbox{and}\\
J_n(M_1'\cap M_{2k-1})&=&\comp id_{L^2(M)}\tn 
\C{L}_{N,N}(L^2(M)^{\otimes _N^{k-1}})\tn  
\comp id_{L^2(M)^{\otimes _N^{n-k}}}  
\end{eqnarray*} 
hold for $1\leq k \leq n$. Hence
the standard invariant of $A\subset B$ is isomorphic to
$$ \begin{array}{rcccccccc}
\comp \B{1}=N'\cap N &\subset &N' \cap M_1 & \subset &N'\cap M_3&\subset & 
N'\cap M_5 &\subset & 
\ldots \\
 & &\cup & &\cup & &\cup & &  \\
\comp \B{1}& =&M_1'\cap M_1 & \subset &M_1' \cap M_3& \subset &M_1' \cap M_5 &
\subset
&\ldots .
\end{array}   $$ 
But this tower of finite dimensional von Neumann algebras is also isomorphic
to the standard invariant of the subfactor $N\subset M_1$: 
From \cite{PP2} (compare also \cite{Sch3}) we conclude that
there is an isomorphism $J_2$ from $M_3$ onto the 
basic construction $B(N,M_1)$ such that the restriction of $J_2$ onto
$M_1$ is the identity. The same argument shows that we may identify
$M_5$ with the basic construction $B(M_1,M_3)$ and so on. Hence the Jones
tower for $N\subset M_1$ is isomorphic to
$$N\subset M_1\subset M_3 \subset M_5 \subset \ldots .$$

Now let us assume that $M$ is the hyperfinite $\zws $-factor. $[M_1:M]<\infty $
implies that $M_1$ is also isomorphic to the hyperfinite $\zws $-factor
(see Lemma 2.1.18 in \cite{Jo}). Since  $A\subset B$ and $N\subset M_1$ have
the same standard invariant and are of finite depth, they are isomorphic
(according to Popa's result 
\cite{P1}).

\abs
We modify the preceding construction of the subfactors by always using 
$\sigma $ instead of $\sigma $ and
$\osigma $ alternately. For example we obtain Wenzl's Hecke algebra
subfactors by this method, if we take a $C^*$-tensor category
$\C{C}$ which, roughly speaking, consists of finite-dimensional
representations of the quantum group $U_qsl_k$, where $q$ is a root
of unity. (This seems to be known, although the author does not know 
any detailed
reference. The author intends to publish a detailed approach elsewhere.)

Here the factors $A$ and $B$ are given by the tower  
$$\hspace{-3mm} \begin{array}{ccccccccc}
B^1 =(\sigma ,\sigma ) &\subset &B^2=( \sigma ^2  ,\sigma ^2)&
\subset
 & B^3=(\sigma ^3 ,\sigma ^3)   &\subset &
B^4=(\sigma ^4, \sigma ^4)
& \subset &\ldots \\
\cup & &\cup & &\cup & &\cup & & \\
A^1=(\iota , \iota )& \subset &A^2 =(\sigma ,\sigma )& \subset &
A^3= (\sigma ^2, \sigma ^2)  &\subset & A^4= (\sigma ^3, \sigma ^3)
&\subset &\ldots ,
\end{array}$$
where the embedding of $A^n$ into $B^n$ is given by
$T\Pfeil T\times \B{1}_{\sigma }$ and the embedding of $A^n$ (resp. $B^n$) into
$A^{n+1}$ (resp. $B^{n+1}$) by $T\Pfeil \B{1}_{\sigma } \times T$.

Additionally to the assumption $d(\sigma )>1$, we require 
the following condition for the object $\sigma $:

\begin{ass} There is a $\nu \in \nat \setminus \{1\}$ such that $\osigma $ is 
contained in $\sigma ^{\nu -1}$ (and consequently $\iota \leq \sigma ^{\nu }$
holds). \label{Ass} \end{ass}

\B{Remark: } One easily shows that $\sigma $ fulfils the condition if for
every irreducible subobject $\tau $ of $\sigma $ there exists  a 
$\nu _{\tau }\in \nat $ such that $\iota \leq \tau ^{\nu _{\tau }}$.

\abs

We get a similar result as before:
\begin{thm} Let $\sigma $ be an object of $\C{C}$
satisfying Assumption ~\ref{Ass}. The 
sequences $(A^n)_{n\in \nat }$ and $(B^n)_{n\in \nat }$ satisfy the
the assumptions of Wenzl's Theorem 1.5 in \cite{We}.
The subfactor $A\subset B$
defined by this tower has the following properties:
\begin{enumerate}
\item[(i)] $[B:A] =d(\sigma )^2$.
\item[(ii)] The 
tower
$$\displaylines{ B_m^1 =(\sigma (m+1),\, \sigma (m+1)) \subset \hfill\cr
\qquad \subset  B_m^2 =(\sigma\, \sigma (m+1) ,\, \sigma \, \sigma (m+1)) 
\subset
  \hfill\cr
\qquad \qquad \subset B_m^3 =(\sigma ^2\, \sigma (m+1),\, \sigma ^2 \, \sigma (m+1))  
\subset \hfill\cr
\qquad \qquad \qquad \subset
B_m^4 =(\sigma ^3\, \sigma (m+1),\, \sigma ^3 \, \sigma (m+1)) \subset 
\ldots \hfill } $$
defines a unique $\zws $-factor $B_m$.
The Jones tower for $A\subset B$ is isomorphic to
$$B_{-1}=A \subset B_0 =B \subset B_1 \subset B_2 \subset B_3 \subset
 \ldots ,$$ 
where the Jones projection $f_m\in B_{m+1}$ is given by \newl
$f_m= \frac{1}{d(\sigma )}\,  
 \B{1}_{\sigma (m)}\times \Bigl(\oR _{\sigma }\circ \oR _{\sigma }^*\Bigr)
 \in B_{m+1}^1$ for $m$ even and \newl
$f_m= \frac{1}{d(\sigma )}\, \B{1}_{\sigma (m)} \times \Bigl(R_{\sigma }\circ
R_{\sigma }^* \Bigr)
 \in B_{m+1}^1$ for $m$ odd.
\item[(iii)] The standard 
invariant for $A\subset B$ is
$$ \hspace*{-9mm}\begin{array}{rcccccccccc}
(\iota ,\iota ) &\subset &(\sigma ,\sigma ) & \subset &(\sigma \osigma , 
\sigma \osigma )
&\subset & (\sigma \osigma \sigma , 
\sigma \osigma \sigma ) &\subset & (\sigma \osigma \sigma \osigma ,
\sigma \osigma \sigma \osigma ) & \subset 
\ldots \\
 & &\cup & &\cup & &\cup & &\cup  \\
& &(\iota , \iota ) & \subset & (\osigma ,\osigma )& \subset &(\osigma \sigma ,
\osigma \sigma ) &\subset &(\osigma \sigma \osigma , \osigma \sigma \osigma )&
\subset
\ldots .
\end{array}   $$
\end{enumerate} \label{Tvar} \end{thm}
Hence the standard invariant of this subfactor is the same as the standard 
invariant of the subfactor constructed with the same
object $\sigma $ in Theorem ~\ref{TUnt}. 
It follows from the main result of \cite{P1} that the two subfactors are 
isomorphic.

\B{Proof:\ }The proof is similar to that of the preceding results.
\newl (1) 
If $\rho $ is an object of $\C{C}$, we have
$$ \rho \leq \sigma ^{\nu } \rho \leq \sigma ^{2\nu } \rho \leq
   \sigma ^{3\nu } \rho \leq \ldots  $$
with the consequence that
\begin{equation} [[\C{C}, \rho ]] \subset [[\C{C}, \sigma ^{\nu } \rho ]] 
\subset [[\C{C}, \sigma ^{2\nu } \rho ]]\subset \ldots . \label{sket}
\end{equation}
Since $\C{C}$ is finite, there is a $k_{\rho }\in \nat $ such that the tower
(~\ref{sket}) becomes stationary from $[[\C{C}, 
\sigma ^{k_{\rho }\nu }\rho ]]$ on.
By putting $\rho = \sigma ^n \sigma (m+1)$ and applying Observation 
~\ref{obs} (5),
we obtain that the towers
$$B_m^1 \subset B_m^2 \subset B_m^3 \subset B_m^4 \subset \ldots $$
($m\in \{-1,0\} \cup \nat $) are periodic with period 
$\nu $.

\abs
\noindent (2) 
We have $\osigma \sigma \leq \sigma ^{\nu }$ as well as $\sigma \osigma \leq
\sigma ^{\nu }$ and obtain
\begin{equation} \sigma ^n \leq \sigma ^n\sigma (2m) \leq \sigma ^n \sigma ^{\nu m} 
 \label{E221} \end{equation}
for $m,\, n\in \nat \cup \{0\}$.
From Part (1) we conclude  
$[[\C{C}, \sigma ^{n+\nu m}]] =[[\C{C}, \sigma ^{n}]]$,
if $n$ is sufficiently large. Relation (~\ref{E221}) implies   
$[[\C{C}, \sigma ^{n}]] = [[\C{C}, \sigma ^{n} \sigma (2m)]]$.
Therefore the numbers of the direct summands of
$B_{2m-1}^{n+1}$ and $B_{2m+1}^{n+1}$ agree for a sufficiently large $n$.
Lemma ~\ref{Lfund} implies that
$B_{2m+1}^{n+1}$ is the basic construction for $B_{2m-1}^{n+1} \subset 
B_{2m}^{n+1}$.
If $[[\C{C}, \sigma ^{n} \sigma (2m)]]=[[\C{C}, \sigma ^{n} \sigma (2m+2)]]$
then 
$[[\C{C}, \sigma ^{n} \sigma (2m)\osigma ]]=
[[\C{C}, \sigma ^{n} \sigma (2m+2)\osigma ]]$. By applying Lemma ~\ref{Lfund}
again, we obtain that $B_{2m+2}^{n+1}$ is the basic construction for
$B_{2m}^{n+1}\subset B_{2m+1}^{n+1}$.

We will verify in Step (3) and (4)
that the inclusion matrix for $B_m^n \subset B_m^{n+\nu }$ is primitive
if $n$ is sufficiently large.
After having finished, we are able to proceed as in
the proof of Theorem ~\ref{TUnt}, and the proof of the Theorem is completed. 

\abs
\noindent (3) Let $s$ be a multiple of $\nu $ such that
$[[\C{C},\sigma ^s]]\supset [[\C{C},\sigma ^{\nu k}]]$ for every $k\in\nat $
and let $\pi $ and  $\rho $ be irreducible objects of $\C{C}$ such that
$[\pi ],\,[\rho ]\in [[\C{C},\sigma ^s]]$. Applying the Frobenius reciprocity
law twice, we get
$$0 < \dim (\rho ,\sigma ^s) =\dim (\iota ,\sigma ^s \orho )= 
      \dim (\ov{\sigma ^s}, \orho ) =\dim (\osigma ^s, \orho ). $$
Since $\osigma ^s$ is a subobject of $\sigma ^{(\nu -1)s}$, $\orho $ is a 
subobject of $\sigma ^{(\nu -1)s}$. So 
$[\orho] \in [[\C{C},\sigma ^s]]$ and
$0 <\dim (\orho , \sigma ^s) = \dim (\iota , \sigma ^s\rho )$.
Hence $\iota \leq \sigma ^s\rho $ and $\pi \leq \sigma ^s \iota $ 
yields \newl
$\pi \leq \sigma ^{2s} \rho $.

\abs
\noindent (4) From $d(\sigma )>1$ we conclude $B^n_m\not= B^{n+s}_m$.
We intend to see that the matrix $D$ describing the
inclusion
$$\displaylines{B_m^n=(\sigma ^{n-1}\sigma (m+1) ,\sigma ^{n-1}\sigma (m+1))\,
\subset \hfill\cr \hfill \subset\,
B_m^{n+\nu }=(\sigma ^{n+\nu -1}\sigma (m+1) ,\sigma ^{n+\nu -1}
\sigma (m+1))}$$
is primitive for $n > s$ and $m\in \nat \cup \{-1,0\}$.
$D$ is quadratic for $n>s$, as 
$[[\C{C},\sigma ^{n-1}\sigma (m+1)]] =[[\C{C},\sigma ^{n+\nu -1}\sigma (m+1)]].$

Let $\rho $ and $\pi $ be irreducible subobjects of 
$\sigma ^{n-1}\sigma (m+1)$.
Since $\osigma \leq \sigma ^{\nu -1}$, there is an $r\in \nat $ such that
$\rho $ and $\pi $ are subobjects of $\sigma ^{s+r}$. $(\pi ,\sigma ^{s+r})
\neq 0$ implies $(\osigma ^r \pi ,\sigma ^s)\neq 0$. Hence there is an 
irreducible subobject $\tilde{\pi }$ of $\osigma ^r\pi $ such that
$(\tilde{\pi } , \sigma ^s)\neq 0$. We also have an irreducible subobject
$\tilde{\rho }$ of $\osigma ^r\rho $ such that $(\tilde{\rho },
\sigma ^s)\neq 0$.
Part (3) implies $(\tilde{\pi }, \sigma ^{2s} \tilde{\rho })\neq 0$
with the consequence that 
$(\osigma ^r\pi ,\sigma ^{2s}\osigma ^r \rho ) \neq 0$
and $(\pi ,\sigma ^r \sigma ^{2s}\osigma ^r \rho )\neq 0$.
Using $\osigma \leq \sigma ^{\nu -1}$ we get
$(\pi ,\sigma ^{2s+\nu r}\rho ) \neq 0$. Since $\sigma ^s$ has the same 
irreducible
subobjects as $\sigma ^{2s+\nu r}$, $\pi $ is a subobject of $\sigma ^s \rho $.

The matrix $D^{s/\nu }$ describes the inclusion $B_m^n\subset 
B_m^{n+s}$. The simple direct summands of $B^n_m$ and $B^{n+s}_m$ are in
one to one correspondence with the irreducible subobjects of
$\sigma ^{n-1}\sigma (m+1)$ (compare Observation ~\ref{obs} (5)) and 
the coefficient of $D^{s/\nu }$ associated with the 
irreducible subobjects $\rho $ of $\sigma ^{n-1}\sigma (m+1)$ and
$\pi $ of $\sigma ^{n+s-1}\sigma (m+1)$ is
the number of times that $\pi $ is contained in $\sigma ^s \rho $. According
to the considerations from above, this number is always greater than $0$.
So $D$ is primitive. \blacksquare 

\abs
We illustrate the construction of Section ~\ref{SSf}
in case $\C{C}$ is the $C^*$-tensor category $\C{U}_H$
of the finite dimensional unitary corepresentations of a finite dimensional
Hopf-$*$-algebra $H$. Let $\sigma :V\pfeil V \otimes H$ be an object
of $\C{U}_H$ and $A\subset B$ the associated subfactor. 
We introduce the abbreviations
$$V^1 :=\comp, \, V^2:=\ov{V},\,
V^3:=V \otimes \ov{V},\, V^4:= \ov{V}\otimes V\otimes \ov{V},\hbox{ and so on,}
$$
and consider the tower
\begin{equation} \B{L}(V^1) \subset \B{L}(V^2)\subset 
\B{L}(V^3) \subset \ldots 
  \label{ECH1} \end{equation}
of finite dimensional $*$-algebras. 
Let $N$ be the unique $\zws$-factor containing $N^{\infty }:=
\bigcup _{n=1}^{\infty } \B{L}(V^n)$ as an 
ultra-strongly dense $*$-subalgebra.

Since the $*$-algebra $A^n$ is contained in $\B{L}(V^n)$ for every 
$n\in \nat $,
$A$ may be regarded as a subfactor of $N$. We will describe $A$ as a
fixed point algebra $N^K$ under an action of 
the finite dimensional Hopf-$*$-algebra $K:=(H^o)^{cop}$, where 
$H^o:=\{f:H\pfeil \comp:\, f \hbox{ linear}\}$ denotes
the Hopf-$*$-algebra dual to $H$. (Details of the definition of $H^o$ 
can be found in \cite{Koo}.)

There is a one to one correspondence between the finite dimensional
unitary corepresentations of $H$ and the finite dimensional non-degenerate
$*$-representations of $H^o$. If $\rho :V_{\rho }\pfeil V_{\rho } \otimes H$
is a unitary corepresentation of $H$, the associated $*$-representation
$\rho ^o:H^o\pfeil \B{L}(V_{\rho })$ is defined by
$\rho ^o(f)w = (id_{V_{\rho }}\otimes f)\, \rho (w)$
(compare \cite{Koo}).
$\rho ^o$ is also a $*$-representation of $K$,
as 
$H^o$ and $K$ are the same $*$-algebras. We write $\rho ^c$ for $\rho ^o$, 
whenever we 
regard $\rho ^o$ as a $*$-representation of $K$. 

\begin{lem} \begin{enumerate}
\item[(i)] Let $\mu $ be a nondegenerate
$*$-representation of the Hopf-$*$-algebra
$K$ on the Hilbert space $V_{\mu }$. Then 
an action $\alpha _{\mu }$ of $K$ on $\B{L}(V_{\mu })$ is defined by
$$\alpha _{\mu } (a)x: =\sum _{(a)}\mu (a^{(1)})\, x \, \mu (S(a^{(2)}))$$
($\Delta (a) =\sum _{(a)} a^{(1)}\otimes a^{(2)}$,
Sweedler notation for $K$) for every $a\in K$ and $x\in \B{L}(V_{\mu })$.
The relation
$$\B{L}(V_{\mu })^K=\mu (K)' $$
holds. 
\newl In particular, the fixed point algebra $\B{L}(V_{\rho })^K$ under the action
$\alpha _{\rho ^c}$ of $K$ is equal to $(\rho ,\rho )$.
\item[(ii)] $\alpha _{\rho ^c}$ is unitarily equivalent to the 
representation $(\orho \otimes
\rho )^o$ of $H^o$, if we regard $\alpha _{\rho ^c}$ as a 
representation of $H^o$
and endow $\B{L}(W)$ with the inner product given by the normalized trace 
$\tr $ of $\B{L}(W)$. In particular, the representation
$\alpha _{\rho ^c}$ respects the $*$-operation. \end{enumerate}
\label{LCH1} \end{lem}

\B{Proof:\ } (i) $ \alpha _{\mu }$ is an action of $K$ by 
Example 2.5 (2) in \cite{PeSz}. The proof of $\mu (K)'\subset
\B{L}(V_{\mu })^K$ is easy. For $x\in \B{L}(V_{\mu })^K$ and $a\in K$ we get
$$\displaylines{
x\cdot \mu (a)\, =\sum _{(a)} \epsilon (a^{(1)})x\cdot
  \mu (a^{(2)})\, =  \hfill\cr
\sum _{(a),(a^{(1)})} \mu (a^{(1)(1)})\cdot x \cdot \mu (S(a^{(1)(2)}))
\cdot \mu (a^{(2)})\, =\hfill }$$
$$\displaylines{
\sum _{(a),(a^{(2)})} \mu (a^{(1)})\cdot x \cdot \mu (S(a^{(2)(1)}))
\cdot \mu (a^{(2)(2)})\, =\hfill\cr
\sum _{(a)} \mu (a^{(1)})\cdot x \cdot \mu \Bigl(\sum _{(a^{(2)})} 
S(a^{(2)(1)})
\cdot a^{(2)(2)}\Bigr)\, =\hfill\cr 
\sum _{(a)} \mu (a^{(1)})\cdot x \cdot \mu (\epsilon (a^{(2)})\B{1})\, 
=\hfill\cr
\sum _{(a)} \mu \Bigl(\epsilon (a^{(2)})\, a^{(1)} \Bigr)\cdot x =\, \, \mu (a)
\cdot x. \hfill } $$
Therefore $x$ belongs to $\mu (K)'$.

\abs

\noindent (ii) Let $(a_{ij})_{i,j}$ be the matrix coefficients of
$\rho $ with respect to an orthonormal base
$\C{B}:=(w_1,\ldots ,w_s)$ of $V_{\rho }$. We get
\begin{equation} 
\rho ^o(f)\, w_j \,=\sum _{i=1}^s f(a_{ij})\, w_i \qquad \hbox{and}
\qquad \orho ^o(f)\, \ov{w_j} \,=\sum _{i=1}^s \, f(S(a_{ji}))\, \ov{w_i}.
\label{ECH2} \end{equation}
Let $\epsilon _{ij}\in \B{L}(W)\, ( i,j=1,\ldots ,s)$ be defined by
$\epsilon _{ij}\, w_k= \delta _{j,k} w_i$.
The unitary operator
$U:\B{L}(V_{\rho }) \pfeil \ov{V_{\rho }}\otimes V_{\rho }$
determined by $\sqrt{s} \, \epsilon _{ij} \Pfeil \ov{w_j} \otimes w_i$ 
intertwines the representations $\alpha _{\rho ^c}$ and
$(\orho \otimes \rho )^o$ of $H^o$:
$$\displaylines{
\sqrt{s}\, U\, \alpha _{\rho ^c}(f) \, \epsilon _{ij} \, =\hfill\cr
\sqrt{s}\, U \sum _{(f)} \rho ^o(f^{(2)})\cdot \epsilon _{ij}
\cdot \rho ^o(S^o(f^{(1)}))
\, = \hfill \hbox{(by using (~\ref{ECH2}))} \cr
\sqrt{s}\, U\sum_{k,l=1}^s f^{(2)}(a_{ki})\, S^o(f^{(1)})(a_{jl})\,
\epsilon _{kl}\,
 =\hfill\cr
\sum_{k,l=1}^s S^o(f^{(1)})(a_{jl})\, \ov{w_l}\, \otimes \,
f^{(2)}(a_{ki})\, w_k\, =\hfill\cr
\sum_{l=1}^s f^{(1)}(S(a_{jl}))\, \ov{w_l} \otimes    
\sum_{k=1}^s f^{(2)}(a_{ki})\, w_k  \, =\hfill
\hbox{(by using (~\ref{ECH2}))}
\cr
(\orho \otimes \rho )^o(f)\, \ov{w_j} \otimes w_i \,= \qquad  \qquad 
\qquad\qquad 
(\orho \otimes \rho )^o(f) \, U\, \sqrt{s}\, \epsilon _{ij} \hfill
}$$
for $f\in H^o$ and $i,j=1,\ldots ,s$.
(For the first '$=$' observe that the Sweedler notation uses the
comultiplication of $H^o$ and not of $K$.)
\blacksquare

\begin{lem} There is an action $\alpha _{\infty } $ of $K$ on $N_{\infty }$
given by 
\begin{equation} \alpha _{\infty }(f)\, x\, =\alpha _{\ov{\sigma (n)}^c }(f)\, x \qquad
\hbox{for $x\in \B{L}(V^{n+1})$.} 
\label{ECH3}\end{equation}
\label{LCH2}
\end{lem}

\B{Proof:\ } It is necessary to show that Definition (~\ref{ECH3}) is 
compatible with the embedding in the tower (~\ref{ECH1}). 
(To this purpose we need $K=(H^o)^{cop}$ 
instead of $H^o$).
Let
$\rho:V_{\rho }\pfeil V_{\rho }\otimes H$ and 
$\tau :V_{\tau }\pfeil V_{\tau }\otimes H$ be unitary corepresentations of
$H$. We prove 
$$\alpha _{(\rho \otimes \tau )^c}(f) (\B{1}_{V_{\rho }}\otimes x)\, =\,
\B{1}_{V_{\rho }}\otimes \alpha _{\tau ^c}(f)x $$
for $f\in K$ and $x\in \B{L}(V_{\tau })$.
An easy computation shows
$$(\rho \otimes \tau )^c(f) \,=\sum _{(f)} \rho ^c(f^{(2)}) \otimes
\tau ^c(f^{(1)}) \qquad \hbox{for $f\in K$}   $$
(where the Sweedler notation for $K$ is used). From now on $S$ denotes the
antipode of $K$. We obtain
$$\displaylines{
\alpha _{(\rho \otimes \tau )^c}(f) (\B{1} \otimes x)\, =\hfill\cr
\sum _{(f)} \Bigl(\sum _{(f^{(1)})}\rho ^c(f^{(1)(2)}) \otimes
\tau ^c(f^{(1)(1)}) \Bigr) \,\cdot \Bigl(\B{1}\otimes x\Bigr) \,\cdot 
\hfill\cr  
\hfill \Bigl( \sum _{(S(f^{(2)}))}\rho ^c\bigl( (S(f^{(2)}))^{(2)} \bigr) 
\otimes
\tau ^c\bigl( (S(f^{(2)}))^{(1)} \bigr) \Bigr) \, =\cr
\sum _{(f)} \sum _{(f^{(1)})} \sum _{(f^{(2)})} 
\rho ^c(f^{(1)(2)}) \cdot \B{1}
\cdot \rho ^c(S(f^{(2)(1)})) \otimes 
\tau ^c(f^{(1)(1)}) \cdot x
\cdot \tau ^c(S(f^{(2)(2)}))\, =\hfill\cr
\sum _{(f)} \sum _{(f^{(1)})} \sum _{(f^{(1)(2)})} \rho ^c(f^{(1)(2)(1)})\cdot 
\rho ^c(S(f^{(1)(2)(2)})) \otimes 
\tau ^c(f^{(1)(1)}) \cdot x
\cdot \tau ^c(S(f^{(2)}))\, =\hfill\cr
\sum _{(f)} \sum _{(f^{(1)})} \rho ^c(\epsilon (f^{(1)(2)}) \,
\B{1} ) \otimes 
\tau ^c(f^{(1)(1)}) \cdot x
\cdot \tau ^c(S(f^{(2)}))\, =\hfill\cr
\B{1}\otimes \sum _{(f)}\, \tau ^c\Bigl(\sum _{(f^{(1)})} \epsilon (f^{(1)(2)})
f^{(1)(1)}\Bigr) \cdot x
\cdot \tau ^c(S(f^{(2)}))\, =\qquad \qquad
\B{1}\otimes \alpha _{\tau ^c}(f)x.
 }$$
(For the second $'='$ observe
$\Delta (S(a))=\sum _{(a)} S(a^{(2)})\otimes S(a^{(1)})$.) \blacksquare 

\abs 
\noindent We extend $\alpha _{\infty }$ to an action $\alpha $ of $K$ on $N$:
by Lemma ~\ref{LCH1} (ii), we know that $\alpha _{\ov{\sigma (n)}^c}$ is a
$*$-representation of the finite dimensional $C^*$-algebra $K$ on the 
Hilbert space $\B{L}(V^n)$ for every $n\in \nat $. Hence
$||\alpha _{\ov{\sigma (n)}^c}(f)|| \leq ||f||$ for $f\in K$, 
where $||f||$ denotes the $C^*$-norm
of $f$. We conclude that the linear operator $\alpha _{\infty }(f):N^{\infty }
\pfeil N^{\infty }$ is continuous and satisfies the relation 
$||\alpha _{\infty }(f)||\leq ||f||$, if we regard $N^{\infty }$ as a 
subspace of the Hilbert space $L^2(N)$. Hence $\alpha _{\infty }(f)$ has a 
unique extension
to a continuous linear operator $u(f):L^2(N) \pfeil L^2(N)$. It is
straightforward to verify that $f\in K \pfeil u(f)$ is a 
nondegenerate $*$-representation
of $K$ on $L^2(N)$.

Now
\begin{equation} \alpha (f)\, x =\sum _{(f)} u(f^{(1)})\, x \, u(S(f^{(2)})) \qquad
\hbox{($f\in K$)}
\label{ECH4} \end{equation}
defines an action $\alpha $ of $K$ on $\B{L}(L^2(N))$ (see Lemma ~\ref{LCH1}
(i)).
For \linebreak 
$n\in N_{\infty }\subset N\subset \B{L}(L^2(N))$ we have
$\alpha (f)\, n=\alpha _{\infty }(f)\, n$, as the computation
$$\displaylines{
\lambda (\alpha (f)n) \, \ov{m}\,=
\sum _{(f)} \ov{\alpha _{\infty }(f^{(1)})\Bigl(n\cdot 
\alpha_{\infty }(S(f^{(2)}))m\Bigr)}
\,= \hfill\cr
\sum _{(f)}\ov{\alpha _{\infty }\Bigl(f^{(1)(1)}\Bigr)n\cdot
 \alpha _{\infty }\Bigl(f^{(1)(2)}\, S(f^{(2)})\Bigr)m }\, =\hfill\cr
\sum _{(f)}\ov{\alpha _{\infty }(f^{(1)})n\cdot
 \alpha _{\infty }\Bigl(f^{(2)(1)}\, S(f^{(2)(2)})\Bigr)m }\, =\hfill\cr
\sum _{(f)}\ov{\alpha _{\infty }(f^{(1)})n\cdot \,
 \alpha _{\infty }(\epsilon (f^{(2)})\B{1})m} \, =\hfill\cr
\ov{\alpha _{\infty }\Bigl(\sum _{(f)} \epsilon (f^{(2)})\, f^{(1)}\Bigr)n\cdot
 m }\,= \qquad \qquad \qquad 
\lambda (\alpha _{\infty }(f)n) \,\, \ov{m}\hfill      }$$
for $m\in N^{\infty }$ shows.
Since $x\in \B{L}(L^2(N))\Pfeil \alpha (f)x \in \B{L}(L^2(N))$ is 
strongly continuous and
$\alpha (f)x\in N$ for $x\in N^{\infty }$, we find that $\alpha (f)$
maps $N$ into $N$ for every $f\in K$. So the restriction of $\alpha $ to an 
action of 
$K$ on $N$ extends $\alpha _{\infty }$. We denote this restriction by 
$\alpha $ again.

\begin{lem} The fixed point algebra $N^K$ under the action $\alpha $ is equal
to $A$. \end{lem}

\B{Proof:\ } There is a unique minimal central projection $e$ of $K$ such that 
$fe=\epsilon (f)e$ for every $f$ of $K$ (see \cite{PeSz}, Theorem 2.2 (4)).
(If $K$ is the group algebra $\comp [G]$ for a finite
group $G$, then $e: =1/|G|\, \sum _{g\in G} 
g$.) By Proposition 2.12 in \cite{PeSz},
$$N^K=\{\alpha (e)x:\, x\in N\}. $$
From Lemma ~\ref{LCH1} (i) we conclude $A^{\infty } =N^{\infty } \cap
N^K$. For $x\in N$ there is a sequence $(x_i)_{i\in \nat }\subset N^{\infty }$ 
strongly converging to $x$ (in $L^2(N)$). Then the sequence 
$(\alpha (e)x_i)_i$ converges strongly to $\alpha (e)x$, hence
$\alpha (e)x\in A$ and  $N^K\subset A$. 
The other inclusion is obvious. \blacksquare

\begin{lem} $(N^K)'\cap N\, =\comp \B{1}$. \end{lem}

\B{Proof: }It is not difficult to see that the towers  
$(A^n)_n$ and $(\B{L}(V^n))_n$ describing the subfactor $N^K\subset N$
fulfil the periodicity assumptions of \cite{We}, Theorem 1.5. 
Note that
$$E_n:\B{L}(V^n)\pfeil A^n,\, x\Pfeil \alpha (e)x, $$
is the conditional expectation from $\B{L}(V^n)$ onto the fixed point
algebra \linebreak
$\B{L}(V^n)^K =A^n$ corresponding to the unique normalized
trace on $\B{L}(V^n)$ (use Proposition 2.12 in \cite{PeSz} and observe
Lemma ~\ref{LCH1} (ii)).

We will apply Theorem 1.6 in \cite{We} in order to estimate the dimension
of $(N^K)'\cap N$.
For every $n\in \nat $ there is a minimal projection $p$ in $A^{2n+1}$ 
such that the
restriction $\ov{\sigma (2n)}\rest pV^{2n+1} $ of the corepresentation 
$\ov{\sigma (2n)}$ is equivalent to the identity corepresentation.
So $\dim \, pV^{2n+1}=1$, and 
$$\dim \,  (N^K)'\cap N \leq \hbox{dim }p\,
\Bigl((A^{2n+1})'\cap \B{L}(V^{2n+1})\Bigr)\,=1$$
follows. \blacksquare

\abs

Similarly, it is possible to write $B$ as a fixed point algebra under an 
action of $K$:
\newl $N\otimes \B{L}(V)$ is the $\zws$-factor containing
$$\bigcup _{n=1}^{\infty } \B{L}(V^n\otimes V) \,= 
\bigcup _{n=1}^{\infty } \B{L}(V^n)\otimes \B{L}(V) $$ as an ultra-strongly 
dense $*$-subalgebra. As before we introduce an action $\beta $ of $K$ on
$N\otimes \B{L}(V)$ such that
$$\beta (f)x\, =\alpha _{\ov{\sigma (n)}\otimes \sigma }(f)\, x \qquad
\hbox{for $x\in
\B{L}(V^n\otimes V)$.}   $$
We get $(N\otimes \B{L}(V))^K=B$ as well as 
$((N\otimes \B{L}(V))^K)' \cap N\otimes \B{L}(V)=\comp \B{1}$.
Hence the subfactor $A\subset B$ may be written as an inclusion
of fixed point algebras
\begin{equation} N^K\subset (N\otimes \B{L}(V))^K, \label{ECH6} \end{equation}
where the embedding of $N$ into $N\otimes \B{L}(V)$ is given by
$n\Pfeil n\otimes \B{1}$.

The restriction of the action $\beta $
onto $N$ is in general not equal to $\alpha $, as the proof of Lemma
~\ref{LCH2} does not work for $x\otimes id_{V_{\tau}}$ ($x\in
\B{L}(V_{\rho })$) instead of $id_{V_{\rho }}\otimes x$ ($x\in
\B{L}(V_{\tau })$) except in the group case.
If $H$ is equal to the commutative Hopf-$*$-algebra
$\Fun (G):=\{f:G\pfeil \comp \}$ for a finite group $G$, 
$\sigma $ may be considered
as a unitary finite dimensional representation of $G$ and the 
action $\beta $ of $G$ is equal to 
$$\displaylines{
\alpha \otimes \Ad \sigma :G\pfeil \Aut N\otimes \B{L}(V), \hfill \cr
(\alpha \otimes \Ad \sigma )(g) (n\otimes x)=
 \alpha (g)n\otimes \sigma (g)x\sigma (g)^{-1} \hbox{\ for\ }
 g\in G,\, n\in N,\, x\in \B{L}(V).  \hfill        }$$

\abs
By applying Lemma ~\ref{LCH1} (ii), one  obtains
\begin{lem} The action $\alpha $ (resp. $\beta $) is an outer action of $K$
on $N$ (resp. $N\otimes \B{L}(V)$)
if there is an $n\in \nat $ such that $(\sigma \osigma )^n$ 
(resp. $(\osigma \sigma )^n$) contains every irreducible unitary
corepresentation of $K$. \label{L21alph} \end{lem}

\abs 
\noindent In particular, the exposition from above describes a method for 
constructing an 
outer action of a finite dimensional Hopf-$*$-algebra $K$ on the hyperfinite
$\zws $-factor, if we use a unitary corepresentation $\sigma $ of
$H=(K^{cop})^o$ as in Lemma ~\ref{L21alph}. 
(For example, choose the comultiplication $\Delta $ of $H$ as 
the corepresentation $\sigma $.) 
 
\section{The $C^*$-tensor category associated with the constructed
subfactors} 
In this Section we determine the $C^*$-tensor category of the 
$(A,A)$-bimodules for
the subfactors $A\subset B$ from Section 3. We assume the
notation used there.

Let $\rho $ be an object of the finite $C^*$-tensor category
$\C{C}$. $\C{C}_{\rho }$ denotes the following
full $C^*$-tensor subcategory of $\C{C}$:
\newl the objects of $\C{C}_{\rho }$ form the smallest subset $\C{O}_{\rho }$
of $\Ob \C{C}$
with the following properties:
\begin{enumerate}
\item[(a)] $\iota ,\,\rho \in \C{O}_{\rho }$.
\item[(b)] If $\tau \in \C{O}_{\rho }$ then every object equivalent to
$\tau $ and every subobject of $\tau $ belongs to $\C{O}_{\rho }$.
\item[(c)] If $\tau ,\,\phi \in \C{O}_{\rho }$ then $\ov{\tau }\in 
\C{O}_{\rho }$ and $\tau \phi \in \C{O}_{\rho }$.
\item[(d)] Finite direct sums of objects of $\C{O}_{\rho }$ are objects of
$\C{O}_{\rho }$.
\end{enumerate}
Obviously, $\C{C}_{\rho }$ is a finite $C^*$-tensor category.

\begin{thm} Let $A\subset B$ be the subfactor from Section ~\ref{SSf} or
Theorem ~\ref{Tvar}. There is an equivalence $G:\C{C}_{\sigma \osigma }\pfeil
\C{B}_{A\subset B}$ of the $C^*$-tensor categories $\C{C}_{\sigma \osigma }$
and $\C{B}_{A\subset B}$ such that $_AL^2(B)_A$ is equivalent to 
$G(\sigma \osigma )$. \label{TCB} \end{thm}

\begin{cor} For every finite $C^*$-tensor category $\C{C}$ there is a subfactor
$A\subset B$ of the hyperfinite factor $B$ with finite index such that 
the $C^*$-tensor category $\C{B}_{A\subset B}$ is equivalent to $\C{C}$. 
\end{cor}

\B{Proof:\ } In Theorem ~\ref{TCB} take an object $\sigma $, which contains
every irreducible object of $\C{C}$ as a subobject. \blacksquare

\begin{cor} Let $N\subset M$ and $P\subset Q$ be two inclusions of
$\zws $-factors with finite index and finite depth. If the standard 
invariants of $N\subset M$ and $P\subset Q$ are isomophic, the
$C^*$-tensor categories $\C{B}_{N\subset M}$ and $\C{B}_{P\subset Q}$ are
equivalent. \end{cor}

\B{Proof:\ } Assume that the standard invariants of $N\subset M$ and
$P\subset Q$ are isomorphic.
 Let $A\subset B$ (resp. $C\subset D$)
be the subfactor from Section ~\ref{SSf} 
with $\C{C}=\C{B}_{N\subset M}$ (resp. $\C{C}=\C{B}_{P\subset Q}$)
and $\sigma = _NL^2(M)_N$ (resp. $_PL^2(Q)_P$).
According to Example ~\ref{X21}, the standard invariants of $N\subset M_1$ and
$A\subset B$ are isomorphic as well as the the standard invariants of
$P\subset Q_1$ and $C\subset D$. Moreover, the standard invariant
of $N\subset M_1$ is isomorphic to the standard invariant of $P\subset Q_1$. 
$B$ and $D$ are hyperfinite $\zws $-factors, so Popa's classification
\cite{P1} of finite depth subfactors implies that the subfactors
$A\subset B$ and $C\subset D$ are isomorphic. Theorem ~\ref{TCB} shows that 
the $C^*$-tensor category $\C{B}_{A\subset B} = \C{B}_{C\subset D}$ is
equivalent to $\C{B}_{N\subset M}$ as well as to $\C{B}_{P\subset Q}$.
\blacksquare

\abs Intending to work with a modified
version of the $C^*$-tensor category $\C{C}_{\sigma \osigma }$ in the proof
of Theorem ~\ref{TCB}, we carry
out this modification in the following observations.

\subsection{Observations \label{O23}}
(1) Let $\C{D}$ be a strict $C^*$-tensor category and let $\C{D}^P$ denote
the following $C^*$-tensor category:
\newl The objects of $\C{D}^P$ are pairs $(\rho ,P )$, where $\rho $ is an 
object of $\C{D}$ and $P$ a projection of $(\rho ,\rho )\setminus \{0\}$. 
The space of
morphisms is
$$((\rho ,P), (\psi ,Q)):=\{T\in (\rho ,\psi ):\,QTP =T\}$$
for any objects $(\rho ,P)$ and $(\psi ,Q)$ of $\C{D}^P$. Endowed with the
Banach space structure, the composition law, and the $*$-operation
from $\C{D}$, $\C{D}^P$ is a $C^*$-category
with subobjects, as one easily sees.

Let the product of objects be given by $(\rho ,P)\cdot (\psi ,Q)=
(\rho \psi, P\times Q)$. The product $S_1\times S_2$ of morphisms
$S_i\in ((\rho _i,P_i), (\psi _i,Q_i))$, $i=1,2$, is defined as in $\C{D}$.
Obviously $S_1\times S_2$ belongs to 
$\bigl( (\rho _1 \rho _2,P_1\times P_2),
(\psi _1 \psi _2,Q_1 \times Q_2) \bigr)$. The unit object in $\C{D}^P $ is
$(\iota ,\B{1}_{\iota })$. Endowed with this product structure, $\C{D}^P$
is a strict $C^*$-tensor category.

Let $F:\C{D}\pfeil \C{D}^P$ be the functor given by $F(\rho )=(\rho ,
\B{1}_{\rho })$ for $\rho \in \Ob \C{D}$ and
$F(T)=T\in ((\rho ,\B{1}_{\rho }),(\psi ,\B{1}_{\psi }))$ for
$T\in (\rho ,\psi )$.  
Obviously $F$ is a full and faithful strict $C^*$-tensor functor. 
It is a $C^*$-tensor 
equivalence if $\C{D}$ has subobjects. (Assume that $(\rho ,P)$ is an object of
$\C{D}^P$ and $\phi $ is a subobject of $\rho $ corresponding to $P$. Then
$F(\phi )$ is equivalent to $(\rho ,P)$.)

\abs \noindent (2) Let $\C{D}$ be a strict $C^*$-tensor category and let
$\C{S}$
be a subset of $\Ob \C{D}$ satisfying the following properties:
\begin{enumerate}
\item[(a)] Each equivalence class of $\Ob \C{D}$ contains an element of 
$\C{S}$,
and the unit object $\iota $ belongs to $\C{S}$.
\item[(b)] $\rho ,\phi \in \C{S}$ implies $\rho \phi \in \C{S}$.
\end{enumerate}

We get a modified $C^*$-tensor category $\C{D}_{\C{S}}$ if we reduce the
object set of $\C{D}$ to the set $\Ob \C{D}_{\C{S}}:=\C{S}$ and take the
remaining structure of $\C{D}$. Obviously the functor
$G:\C{D}_{\C{S}} \pfeil \C{D}$ defined by
$$G(\phi )=\phi \hbox{\ for $\phi \in \C{S}$\ \ \ and\ \ \ }
G(T)=T \hbox{ for every morphism $T \in \C{D}_{\C{S}}$} $$  
is a strict $C^*$-tensor equivalence.

\abs \noindent (3) We will apply Observation (2) to the $C^*$-tensor category
$\C{D} = (\C{C}_{\sigma \osigma })^P$ if
$d(\sigma ) \neq 1$. Let
$$\C{S} \,=\,\{((\sigma \osigma )^m,P): \,m\in \nat \cup \{0\}, \,
   P \hbox{ a projection $\neq 0$ of\ } 
((\sigma \osigma )^m, (\sigma \osigma )^m) \}.$$

It is not difficult to see that $\C{S}$ satisfies the 
Properties (a) and (b) from Observation (2). Therefore we are able  
to introduce the $C^*$-tensor category
$\C{D}_{\sigma \osigma }:=  ((\C{C}_{\sigma \osigma })^P)_{\C{S}}$,
which we use in the proof of the Theorem. 
By using the Observations (1) and (2) and Proposition ~\ref{PCTK} (i) and (ii),
we see that there is an equivalence $H:\C{C}_{\sigma \osigma }\pfeil  
\C{D}_{\sigma \osigma }$ of the $C^*$-tensor categories 
$\C{C}_{\sigma \osigma }$ and $\C{D}_{\sigma \osigma }$
such that
$H(\sigma \osigma )$ is equivalent to $(\sigma \osigma , 
\B{1}_{\sigma \osigma })$.

\abs
\B{Proof of the Theorem:\ }(1) We only deal with the case of Theorem 
~\ref{TUnt}, as the proofs
do not differ. 
We fix a standard pair $(R_{\phi }, \oR _{\phi })$ of conjugation operators
for every object $\phi $ of $\C{C}$.
For all objects $\rho,\,
\phi ,\,\psi $ of $\C{C}$ we introduce a map
$f(\rho ,\phi ,\psi )$ from $(\phi \ophi ,\psi \opsi )$ into the set
$\B{L}((\rho \phi ,\rho \phi ),(\rho \psi ,\rho \psi ))$
of all linear maps
from $(\rho \phi ,\rho \phi )$ into $(\rho \psi ,\rho \psi )$:
\newl For $K\in (\phi \ophi ,\psi \opsi )$ let
$$\displaylines{
f(\rho ,\phi ,\psi )(K):(\rho \phi , \rho \phi )\pfeil 
 (\rho \psi , \rho \psi ), \hfill\cr
 S\Pfeil \sqrt{\frac{d(\psi )}{d(\phi )}}\,
\B{1}_{\rho \psi }\times R_{\psi }^*\, \circ \, \B{1}_{\rho }\times K
\times \B{1}_{\psi } \, \circ \, S\times \B{1}_{\ophi \psi } \,\circ \,
\B{1}_{\rho } \times \oR _{\phi }\times \B{1}_{\psi }.   \hfill }$$ 
We will often abbreviate $f(\rho ,\phi ,\psi )(K)$ to $f(\rho )(K)$.
We obtain the following properties:  
\begin{itemize}
\item[(a)] $f(\rho ,\phi ,\psi )(K)$ is $(D,D)$-linear for every $K\in
(\phi \ophi ,\psi \opsi )$ where \newl $D:=(\rho ,\rho )$.
\item[(b)] $K\in (\phi \ophi ,\psi \opsi ) \Pfeil f(\rho )(K)$ is linear and
injective.
\item[(c)] $f(\rho )(L\circ K)=f(\rho )(L) \circ f(\rho )(K)$ for every 
$K\in (\phi \ophi ,\psi \opsi )$ and \newl $L\in (\psi \opsi ,\tau \otau )$
($\phi,\, \psi ,\, \tau \in \Ob \C{C}$).
\item[(d)] $f(\rho )(K)^* =f(\rho )(K^*)$ for every 
$K\in (\phi \ophi ,\psi \opsi )$,
where the linear space $(\rho \phi ,\rho \phi )$
(resp. $(\rho \psi ,\rho \psi )$) is endowed with the inner product \newl
$\ska{S}{T} =\tr _{\rho \phi }(ST^*)$ (resp. $\tr _{\rho \psi }(ST^*)$). 
\item[(e)] $f(\rho ,\phi ,\phi )(\B{1}_{\phi }) = 
 id_{(\rho \phi ,\rho \phi )}$. 
\item[(f)] 
$f(\rho )(K)(\B{1}_{\rho _1} \times S) = \B{1}_{\rho _1} \times
 f(\rho _2)(K)(S)$ for
$\rho =\rho _1 \rho _2,\,\, K\in (\phi \ophi ,\psi \opsi )$, and
$S\in (\rho _2 \phi ,\rho _2 \phi )$.
\end{itemize}

\abs The Properties (a) and (f) are obvious. Property (e) is an easy 
consequence of the defining relations for the conjugation operator. We show
the remaining relations:

\abs \noindent ad (b): The linearity is obvious. We prove $f(\rho )(K) \neq 0$
for $K\neq 0$.
First we consider the case that $\phi =\psi $  and $K\neq 0$ is positive. The 
computation 
$$\displaylines{
\tr_{\rho \phi \ophi \phi }\Bigl(\B{1}_{\rho \phi }\times R_{\phi }\,\circ \,
f(\rho )(K)(\B{1}_{\rho \phi }) \,\circ \,
\B{1}_{\rho }\times \oR _{\phi }^* \times \B{1}_{\phi } \Bigr)\,= \hfill\cr
\tr_{\rho \phi \ophi \phi }\Bigl(\B{1}_{\rho \phi }\times ( R_{\phi }\circ
R_{\phi }^*)\,\circ \,\B{1}_{\rho }\times K \times \B{1}_{\phi } 
\,\circ \,\B{1}_{\rho }\times (\oR _{\phi }\circ \oR _{\phi }^*) 
\times \B{1}_{\phi }
 \Bigr)\,= \hfill \hbox{(by Equation (~\ref{Mark}))}
 \cr
d(\phi )^{-1} \, \tr _{\rho \phi \ophi }\Bigl(\B{1}_{\rho }\times K \,\circ \,
\B{1}_{\rho }\times (\oR _{\phi }\circ \oR _{\phi }^*)\Bigr)\,= 
\hfill \hbox{(by Equation (~\ref{zwtr}))} \cr
d(\phi )^{-3} \, \tr _{\rho }(\B{1}_{\rho }\times \oR _{\phi }^* \,\circ \,
\B{1}_{\rho }\times K \,\circ \,\B{1}_{\rho }\times \oR _{\phi })\,=
\hfill\cr
d(\phi )^{-2} \tr _{\rho }(\Psi _{\rho }^{\phi }(K)) \,=\,
d(\phi )^{-2} \tr _{\rho \phi }(K) \,>0  \hfill    }$$
yields $f(\rho )(K)\neq 0$.
Now let $\phi ,\, \psi $ and $K\neq 0$ be arbitrary. We get
$f(\rho )(K^*)f(\rho )(K)=
f(\rho )(K^*\circ K) \neq 0$ by applying Property (c).

\abs \noindent ad (c): For $S\in (\rho \phi ,\rho \phi )$ we have 
$$\displaylines{
f(\rho )(L)\Bigl(f(\rho )(K)(S)\Bigr) =\hfill\cr
\sqrt{\frac{d(\tau )}{d(\phi )}}\,\, \B{1}_{\rho \tau }\times 
R_{\tau }^*\,\circ \,
\B{1}_{\rho }\times L\times \B{1}_{\tau } \,\circ\,
\B{1}_{\rho \psi }\times R_{\psi }^*\times \B{1}_{\opsi \tau }\,\circ \,
\B{1}_{\rho }\times K \times \B{1}_{\psi \opsi \tau } \,\circ \hfill \cr
  \hfill \circ S\times \B{1}_{\ophi \psi \opsi \tau } \,\circ \,
\B{1}_{\rho }\times \oR _{\phi }\times \B{1}_{\psi \opsi \tau } 
 \,\circ \,
\B{1}_{\rho }\times \oR _{\psi }\times \B{1}_{\tau } .}$$
The interchange law implies
$$  \B{1}_{\rho }\times \oR _{\phi }\times \B{1}_{\psi \opsi \tau } \,\circ \,
\B{1}_{\rho }\times \oR _{\psi }\times \B{1}_{\tau } \,=
\B{1}_{\rho \phi \ophi }\times \oR _{\psi }\times \B{1}_{\tau } \,\circ \,
\B{1}_{\rho }\times \oR _{\phi }\times \B{1}_{\tau }.$$
By applying this equation and shifting $\oR _{\psi }$ to the left, we obtain
$$\displaylines{
f(\rho )(L)\Bigl(f(\rho )(K)(S)\Bigr) =\hfill\cr
\sqrt{\frac{d(\tau )}{d(\phi )}} \, 
\B{1}_{\rho \tau }\times R_{\tau }^*\,\circ \,
\B{1}_{\rho }\times L\times \B{1}_{\tau } \,\circ\,
\B{1}_{\rho \psi }\times R_{\psi }^*\times \B{1}_{\opsi \tau }\,\circ \,
\B{1}_{\rho \psi \opsi }\times \oR _{\psi }\times \B{1}_{\tau } \,\circ \hfill
\cr
\hfill \circ \,\B{1}_{\rho }\times K \times \B{1}_{\tau} \, \circ \,
S\times \B{1}_{\ophi \tau } \, \circ \,
\B{1}_{\rho }\times \oR _{\phi }\times \B{1}_{\tau }\,= \cr
 f(\rho )(L\circ K)(S) . \hfill }$$ 
By doing so we used the conjugation relation (~\ref{RR2}) for   
$R _{\psi }$ and $\oR _{\psi }$.

\abs \noindent ad (d): We have to verify
\begin{equation} \tr _{\rho \psi }\Bigl(f(\rho )(K)(S)\,\, T^*\Bigr)\,=\,
     \tr _{\rho \phi }\Bigl(S\,\, f(\rho )(K^*)(T)^*\Bigr) \label{contr} 
   \end{equation}
for $S\in (\rho \phi ,\rho \phi )$ and $T\in (\rho \psi ,\rho \psi )$.

Applying the conjugation relation (~\ref{RR2}) for $R_{\psi }$ and 
$\oR _{\psi }$,
we get \newl
$T=\B{1}_{\rho \psi }\times R_{\psi }^* \circ T\times \B{1}_{\opsi
\psi } \circ \B{1}_{\rho }\times \oR _{\psi } \times \B{1}_{\psi }$ and
conclude
$$\displaylines{
\tr_{\rho \psi }\Bigl( f(\rho )(K)(S)\,\, T^*\Bigr) = \hfill \cr
\sqrt{\frac{d(\psi )}{d(\phi )}}\,\, \tr _{\rho \psi }\Bigl(\B{1}_{\rho \psi }
\times R_{\psi }^* \,\circ \, \B{1}_{\rho }\times K \times \B{1}_{\psi }
\,\circ \, S\times \B{1}_{\ophi \psi } \,\circ \, 
\B{1}_{\rho }\times \oR _{\phi } \times \B{1}_{\psi }\,\circ \Bigr.\hfill \cr
\hfill \Bigl. \circ \,
\B{1}_{\rho }\times \oR _{\psi}^*\times \B{1}_{\psi } \,\circ \,
T^*\times \B{1}_{\opsi \psi } \,\circ \, \B{1}_{\rho \psi }\times R_{\psi }
  \Bigr) = \hfill \hbox{(by Equation (~\ref{zwtr}))}  \cr  
d(\phi )^{-1/2} \, d(\psi )^{5/2} \,\, \tr _{\rho \psi \opsi \psi }
\Bigl(\B{1}_{\rho \psi }\times (R_{\psi }\circ R_{\psi }^*) \,\circ \,
\B{1}_{\rho }\times K \times \B{1}_{\psi }
\,\circ \,S\times \B{1}_{\ophi \psi } \,\circ \Bigr.\hfill \cr
  \hfill \Bigl.\circ \, 
 \B{1}_{\rho }\times \oR _{\phi } \times \B{1}_{\psi }\,\circ \,
\B{1}_{\rho }\times \oR _{\psi}^*\times \B{1}_{\psi } \,\circ \,
T^*\times \B{1}_{\opsi \psi } \,
  \Bigr) = \hfill \hbox{(by Equation (~\ref{Mark}))}  \cr
d(\phi )^{-1/2} \, d(\psi )^{3/2} \,\,\tr _{\rho \psi \opsi }
(\B{1}_{\rho }\times K 
\,\circ \, S\times \B{1}_{\ophi } \,\circ \,
\B{1}_{\rho }\times \oR _{\phi } \,\circ \,
\B{1}_{\rho }\times \oR _{\psi }^* \,\circ \,
T^*\times \B{1}_{\opsi } ).   \hfill }$$      
The same calculation for
$$\tr _{\rho \phi }\, \Bigl(S\, f(\rho )(K^*)(T)^*\Bigr) \,=\,
  \ov{\tr _{\rho \phi }\, \Bigl(f(\rho )(K^*)(T) \,S^* \Bigr)  }$$
yields
$$\displaylines{
 \tr _{\rho \phi }\, \Bigl(S\, f(\rho )(K^*)(T)^* \Bigr) \,=\hfill \cr
d(\psi )^{-1/2} \, d(\phi )^{3/2} \,\,
\ov{\tr _{\rho \phi \ophi }
(\B{1}_{\rho }\times K^* 
\,\circ \, T\times \B{1}_{\opsi } \,\circ \,
\B{1}_{\rho }\times \oR _{\psi } \,\circ \,
\B{1}_{\rho }\times \oR _{\phi }^* \,\circ \,
S^*\times \B{1}_{\ophi } )} \hfill \cr
=\, d(\psi )^{-1/2} \, d(\phi )^{3/2} \,\, \tr _{\rho \phi \ophi } 
(S\times \B{1}_{\ophi } \,\circ \,
\B{1}_{\rho }\times \oR _{\phi } \,\circ \,
\B{1}_{\rho }\times \oR _{\psi }^* \,\circ \,
T^*\times \B{1}_{\opsi } \,\circ \, \B{1}_{\rho }\times K) . }$$
So by applying Equation (~\ref{zwtr}) 
we find that the results of both computations
coincide, and that Equation (~\ref{contr}) has been established.

\abs \noindent (2) We will use the maps $f(\rho )$ from Part (1) of the proof
in order to define an $(A,A)$-linear map $F_{m,l}(K) =F(K)$ from
$L^2(B_{m-1})$ into $L^2(B_{l-1})$ for 
$K\in \Bigl(\sigma (m) \ovs{m}, \sigma (l) \ovs{l}\Bigr) =
\Bigl((\sigma \osigma )^m, (\sigma \osigma )^l\Bigr)$, 
$m,l\in \nat \cup \{0\}$.

First we define a linear map $F_{m,l}(K)=F(K)$ from $B_{m-1}^{\infty }\subset
L^2(B_{m-1})$ into $B_{l-1}^{\infty }\subset
L^2(B_{l-1})$ by
$$F(K) \, \ov{S} = \ov{f\Bigl(\ovs{n}, \sigma (m), \sigma (l)\Bigr)(K)(S)}
\, \in B_{l-1}^{n+1}\subset B_{l-1}^{\infty }
$$
for $S\in B_{m-1}^{n+1}$ ($n\in \nat \cup \{0,-1\}$).
In order to do that,
we choose a standard pair $(R_{\sigma (m)}, \oR _{\sigma (m)})$
of conjugation operators for $\sigma (m)$ and $\ovs{m}$,
$m\in \nat \cup \{0\}$, in the following way: we fix a standard pair
$(R_{\sigma }, \oR _{\sigma })$ for $\sigma $ and $\osigma $ and put
$R_{\sigma (0)}:=\oR _{\sigma (0)}: =\B{1}_{\iota }$, 
$R_{\sigma (1)}:= R_{\sigma }$ and $\oR _{\sigma (1)}:=\oR _{\sigma }$.
For $m>1$ we define $(R_{\sigma (m)}, \oR _{\sigma (m)})$ by several 
times applying 
Equation (~\ref{conpr}), in which
we use $(\oR _{\sigma },R_{\sigma })$
as a standard pair of conjugation operators for $\osigma $ and $\sigma $.

By the Property (f) from Part (1), $F(K)$ is well defined. 
According to the rules
in (1),
$$L\in ((\sigma \osigma )^m\, ,(\sigma \osigma )^m)\Pfeil
 f\Bigl(\ovs{n},\sigma (m),\sigma (m)\Bigr)(L)\in 
 \Bigl(\ovs{n} \sigma (m),\ovs{n} \sigma (m)\Bigr)  $$
is a representation of the finite dimensional $C^*$-algebra
$((\sigma \osigma )^m\, ,(\sigma \osigma )^m)$.
It follows
$\Bigl\|f\Bigl(\ovs{n} \Bigr)(L) \Bigr\| \leq ||L||$ for 
every $L\in ((\sigma \osigma )^m\, ,(\sigma \osigma )^m)$. So
\begin{eqnarray*} \Bigl\|f\Bigl(\ovs{n}\Bigl)(K)\Bigr\|^2 &=&
\Bigl\|f\Bigl(\ovs{n}\Bigr)(K)^*\, f\Bigl(\ovs{n}\Bigr)(K) \Bigr\| \,\,= \\
\Bigl\|f\Bigl(\ovs{n}\Bigr)(K^*K)\Bigr\| &\leq & 
||K^*K|| \,\,=\,\, ||K||^2 \end{eqnarray*}
and $F(K)$ is continuous on $B_{m-1}^{\infty }$. Property (a) for $f(\rho )$ 
in (1) shows that $F(K)$ is $(A^{\infty }, A^{\infty })$-linear.
As $B_{m-1}^{\infty }$ is dense in $L^2(B_{m-1})$, $F(K)$ has a unique
extension to a continuous $(A,A)$-linear map from $L^2(B_{m-1})$ onto 
$L^2(B_{l-1})$. We denote this extension by $F(K)$ again. 

\abs In order to note the properties of the assignment 
 $K\Pfeil F(K)$ we have to
introduce unitary $(A,A)$-linear operators
$$W_{m,l}:L^2(B_{m-1}) \ta L^2(B_{l-1}) \pfeil L^2(B_{m+l-1}), \, l,m\in \nat
  \cup \{0\}.$$ 
According to \cite{Sch3},
there is a unique $(A,A)$-linear unitary operator \linebreak
$U_k:L^2(B)^{\ta ^k} \pfeil L^2(B_{k-1})$
for $k\in \nat $ such that
\begin{eqnarray*} \lefteqn{U_k(\ov{x_1}\ta \cdots \ta \ov{x_k}) =} \\
        &=&d(\sigma ) ^{k(k-1)/2}\, \ov{ x_1f_0x_2f_{1,0}x_3f_{2,0}x_4 \ldots
        x_{k-1}f_{k-2,0}x_k} \\
  &=&d(\sigma ) ^{k(k-1)/2}\, \ov{x_1f_{0,k-2}x_2f_{0,k-3}x_3 \ldots 
x_{k-2}f_{0,1}x_{k-1}f_0x_k}  \end{eqnarray*}
for $x_1,\ldots ,x_k\in B$. The abbreviation $f_{r,s},\, 
r,s\in \nat \cup \{0\}$, is defined by
\[ f_{r,s}: =\left\{ \begin{array}{ll}
                    f_r\cdot f_{r+1}\cdot \ldots \cdot f_{s-1}\cdot f_s & 
           \mbox{for $r<s$,} \\
                    f_r                                                &  
           \mbox{for $r=s$,} \\
                    f_r\cdot f_{r-1}\cdot \ldots \cdot f_{s+1}\cdot f_s &  
            \mbox{for $r>s$.}
                    \end{array}  \right. \]
Now we put
\[ W_{m,l}:=\left\{ \begin{array}{ll}
             U_{m+l}\circ (U_m^{-1} \ta U_l^{-1}) & 
               \mbox{for $m,l\in \nat $,} \\
                   r_{L^2(B_{m-1})}            & \mbox{for $l=0$,} \\
                   l_{L^2(B_{l-1})}  & \mbox{for $m=0$,}
                    \end{array}  \right. \]
where $l_{L^2(B_{l-1})}$ (resp.
$r_{L^2(B_{l-1})}$) is the canonical unitary operator from \linebreak
$L^2(B_{l-1})\ta L^2(A)$  (resp. $ L^2(A)\ta L^2(B_{l-1})$) onto 
$L^2(B_{l-1})$.
Especially we have
$$\displaylines{
W_{m,l} \,\,  \ov{x_1 f_0 x_2 \ldots f_{m-2,0}x_m} \, \ta
\ov{y_1 f_0 y_2 \ldots f_{l-2,0}y_l} \, =\hfill \cr
d(\sigma )^{ml} \,\, \ov{x_1f_0 \ldots f_{m-2,0}x_m f_{m-1,0} y_1 f_{m,0} 
\ldots f_{m+l-2,0} y_l}\, = \hfill \cr
d(\sigma )^{ml} \,\, \ov{x_1f_{0,m+l-2} \ldots x_m f_{0,l-1} y_1 f_{0,l-2} 
\ldots
f_0 y_l}  \hfill }$$
for $x_1,\ldots x_m,\, y_1, \ldots ,y_l \in B$ and $m,l\in \nat $.
The maps 
\begin{eqnarray*} F=F_{m,l}:\Bigl((\sigma \osigma )^m, (\sigma \osigma )^l)
\Bigr) 
&\pfeil &
   \C{L}
 _{A,A}\Bigl(L^2(B_{m-1}), \,L^2(B_{l-1})\Bigr),\\
      K&  \Pfeil &F(K), \end{eqnarray*}
satisfy the following relations:
\begin{itemize}
\item[(a)] The map $F_{m,l}$ is linear and bijective for $m,l\in \nat \cup
\{0\}$. 
\item[(b)] $F(L\circ K) = F(L) \circ F(K)$ for every $K\in 
  ((\sigma \osigma )^m, (\sigma \osigma )^l))$ and \newl $L\in
  ((\sigma \osigma )^l, (\sigma \osigma )^k))$ ($m,\, l,\, k \in 
\nat \cup \{0\}$).
\item[(c)] $F(K)^* =F(K^*)$ for every $K\in ((\sigma \osigma )^m, 
(\sigma \osigma )^l))$ ($m,l\in\nat \cup \{0\}$).
\item[(d)] $F(\B{1}_{(\sigma \osigma )^m}) = id _{L^2(B_{m-1})}$ for 
 $m\in \nat \cup \{0\}$.
\item[(e)] $F(K\times L) \,= \, W_{l,r} \circ (F(K) \ta F(L)) \circ W_{m,s}^*$ 
for every 
$K\in ((\sigma \osigma )^m,
(\sigma \osigma )^l)$ and   
$L\in ((\sigma \osigma )^s,
(\sigma \osigma )^r)$ ($m,\, l,\, r,\, s \in \nat \cup \{0\}$).
\end{itemize}
Property (e) and the surjectivity of the map $F_{m,l}$ will be shown in Step
(3) of the proof. The other properties are obvious or immediate consequences
of the Properties (a) - (f) from Step (1).

Let $\C{D}_{\sigma \osigma }$ be the finite $C^*$-tensor category 
introduced in Observation ~\ref{O23}. Using Property (b) we find that there
is a functor $\C{G}:\C{D}_{\sigma \osigma }\pfeil \C{B}_{A\subset B}$ such
that 
\begin{eqnarray*}
\C{G}((\sigma \osigma )^m,P)) &=&F(P)\, L^2(B_{m-1})
\qquad \hbox{and}\\
\C{G}(K) &=& F(K) \qquad \hbox{for\ } K\in 
\Bigl( ((\sigma \osigma )^m,P),\,((\sigma \osigma )^l,Q) \Bigr),
\end{eqnarray*}
where $\Bigl( ((\sigma \osigma )^m,P),\,((\sigma \osigma )^l,Q) \Bigr)$
is regarded as a subspace of the morphism space 
$((\sigma \osigma )^m,(\sigma \osigma )^l)$ 
in the category $\C{C}$.
For objects $\phi =((\sigma \osigma )^m,P)$ and
$\psi =((\sigma \osigma )^l,Q)$ in $\C{D}_{\sigma \osigma }$, let
$U_{\phi \psi }$ be the restriction of $W_{m,l}$ to a unitary operator from
$F(P)L^2(B_{m-1})\ta F(Q)L^2(B_{l-1})$ onto
$F(P\times Q)L^2(B_{m+l-1})$. 
(Observe that $F(P\times Q) = W_{m,l} \circ (F(P)\ta F(Q))\circ W_{m,l}^*$
holds according to Property (e).)

Now we conclude that $(\C{G}, (U_{\phi \psi })_{\phi ,\psi }, id _{L^2(A)})$ 
is a
$C^*$-tensor equivalence. $\C{G}$ is a $C^*$-functor according to the 
Properties (a) and (c). Definition ~\ref{DCequ} (ii) is 
satisfied, Part (a) follows from Property (e), Part
(c) is obvious, and Part (b) is an easy consequence of the fact that
\[ W_{m+l,k} \circ (W_{m,l} \ta id_{L^2(B_{k-1})}) =
W_{m,k+l} \circ (id_{L^2(B_{m-1})} \ta W_{l,k}) \]
is satisfied for  
$k,l,m \in \nat \cup \{0\}$.

From Property (a) we conclude that Definition ~\ref{DCequ} (iii) is
satisfied. Using the results from Observation ~\ref{O23} and Proposition
~\ref{PCTK}, we are able to replace $\C{D}_{\sigma \osigma }$ by the
equivalent $C^*$-tensor category $\C{C}_{\sigma \osigma }$ and get the 
assertion of the Theorem.

\abs \noindent (3) We will prove Property (e) and the surjectivity in
Property (a) from Step (2). In order to facilitate some computations we
develop a special notation:
\newl We fix an object $\rho $ of $\C{C}$ and define the operators
\[ R_{[\rho ]l}^{[m]}:=\left\{ \begin{array}{ll}
     \B{1}_{\rho \sigma (l)} \times \oR _{\sigma }\times \B{1}_{\sigma (m-l)}
      & \mbox{if $l$ is even,} \\
     \B{1}_{\rho \sigma (l)} \times R _{\sigma }\times 
                   \B{1}_{\osigma \sigma (m-l-1) }
      & \mbox{if $l$ is odd and $m\geq l+1$,} \\
      \B{1}_{\rho \sigma (l)} \times R _{\sigma } &
      \mbox{if $l$ is odd and $m= l$}
                    \end{array}  \right.   \]
(belonging to $(\rho \sigma (m),\rho \sigma (m+2))$ ) for $m\geq l\geq 0$,
\[ R_{[\rho ]l}^{*[m]}:=\left\{ \begin{array}{ll}
     \B{1}_{\rho \sigma (l)} \times \oR _{\sigma }^*\times \B{1}_{\sigma 
      (m-l-2)}
      & \mbox{if $l$ is even,} \\
     \B{1}_{\rho \sigma (l)} \times R _{\sigma }^*\times 
             \B{1}_{\osigma \sigma (m-l-3)}
      & \mbox{if $l$ is odd and $m\geq l+3$,} \\
      \B{1}_{\rho \sigma (l)} \times R _{\sigma }^* &
      \mbox{if $l$ is odd and $m= l+2$}
                    \end{array}  \right.   \]
(belonging to $(\rho \sigma (m),\rho \sigma (m-2))$) for $m\geq l+2$
and $l\in \nat \cup \{0\}$,
\[ S^{[m]}:=\left\{ \begin{array}{ll}
     S\times \B{1}_{\sigma (m-l)}
      & \mbox{if $l$ is even,} \\
     S\times \B{1}_{\osigma \sigma (m-l-1)}
      & \mbox{if $l$ is odd and $m\geq l+1$,} \\
      S &
      \mbox{if $l$ is odd and $m= l$}
                    \end{array}  \right.\,\,   \in 
         (\rho \sigma (m),\rho \sigma (m)) \]
for any operator $S\in (\rho \sigma (l), \rho \sigma (l))$ and $m\geq l\geq 0$ 
and
\[K_{[\rho ]l}^{[m]}:= \B{1}_{\rho \sigma (l)} \times K
   \times \B{1}_{\sigma (m-l-2r)} \quad \in \quad 
(\rho \sigma (m), \rho \sigma (m-r+s))\]
for $K\in ((\sigma \osigma )^r, (\sigma \osigma )^s)$, $l$ even and
$m\geq l+2r$. The indices put in brackets \newline $[\hbox{\ \ }]$ are usually 
omitted if
the meaning is clear from the context. 
For instance we have
$$\B{1}_{\rho }\times \oR _{\sigma (m)} \,= R_{m-1}^{2m-2}\circ \ldots \circ
  R_1^2\circ R_0^0 \, = 
  R_{m-1} \ldots    R_1  R_0^0 $$  
(we often omit the sign $\circ $) as well as
$$f(\rho )(K)(S) =d(\sigma )^{\frac{l-m}{2}}\, R_l^{*l+2} R_{l+1}^{*l+4} 
  \ldots  R_{2l-1}^{*3l} K_0^{l+2m}  S^{l+2m} 
   R_{m-1}^{l+2m-2} \ldots  R_1^{l+2}  R_0^l $$
for $K\in ((\sigma \osigma )^m, (\sigma \osigma )^l)$ and
$S\in (\rho \sigma (m),\,\rho \sigma (m))$ ($m,l\in\nat $).

\abs \noindent (a) We prove Property (e) from Step (2) for 
$m=l=r=s=1$.
For $n\in \nat \cup \{0\}$, $\rho = \ovs{n}$ and $S,T\in 
(\rho \sigma , \rho \sigma ) =B^{n+1}$
we compute
$$\displaylines{
f(\rho )(K\times L)(Sf_0T) \,= \hfill\cr
\frac{1}{d(\sigma )}R_2^{*4}\, R_3^{*6}\, L_2^6\, K_0^6\, S^6\, R_0^4\, 
R_0^{*6}\, T^6\, R_1^4\,
R_0^2 \, =\hfill \cr
\frac{1}{d(\sigma )}\, R_2^{*4}\, R_3^{*}\, L_2\, R_1^*\, R_2\, K_0\,
S\, R_0\, R_0^*\, R_1\, T \, R_0^2 \,
=\hfill \cr
\frac{1}{d(\sigma )}\, R_2^{*4}\, R_3^{*}\, L_2\, R_1^*\, K_0\, S\, R_2\,
R_0\, T\, R_0^2 \, =\hfill \cr
\frac{1}{d(\sigma )}\, R_2^{*4}\, R_3^{*}\, L_2\, R_1^*\, K_0\, S\, R_0\,
R_0\, T\, R_0^2 \, = \qquad \hbox{(with $C:=f(\rho )(K)(S)$)} \hfill \cr
\frac{1}{d(\sigma )}\, R_2^{*4}\, R_3^{*}\, L_2\, C\, R_0\, T\, R_0^2 \, =
\hfill\cr 
\frac{1}{d(\sigma )}\, R_2^{*4}\, C\, R_3^*\, L_2\, R_0 \, T\, R_0^2 \, =
\hfill\cr
\frac{1}{d(\sigma )}\, R_2^{*4}\, C\, R_0\, R_1^*\, L_0\, T\, R_0^2 \,=
\qquad \qquad \qquad \hbox{(with $D:=f(\rho )(L)(T)$)} \hfill \cr
\frac{1}{d(\sigma )}\, C^2\, R_0\, R_0^*\, D^2 \,=\qquad \qquad \qquad
f(\rho )(K)(S)\, f_0\, f(\rho )(L)(T). \hfill }$$
During the computation we used the rule
\begin{equation} R_{k\pm 1}^{*m+2}\circ R_k^m \,=\, \B{1} _{\rho \sigma (m)} 
\label{ERc}\end{equation}
several times. (The case '-' follows from the defining
relations (~\ref{RR1}) and (~\ref{RR2}) of the conjugation operators, for the 
case
'+' we have to apply $*$ to the Equations (~\ref{RR1}) and (~\ref{RR2})). 
Moreover, the  
interchange law from Definition ~\ref{DC} (ii) (a) was used very often.
We present a detailed proof of the equation 
$R_2^4 \, R_0^2\, =\, R_0^4\, R_0^2$ (which has been used in line 5 of 
the preceding computation) as an example:
\begin{eqnarray*}
R_2^4 \, R_0^2\,& =& \B{1}_{\rho }\times \B{1}_{\sigma \osigma }\times
\oR _{\sigma } \, \circ \, \B{1}_{\rho }\times \oR _{\sigma }\times 
\B{1}_{\iota }
\, \, = \\
\B{1}_{\rho }\times \oR _{\sigma }\times \oR _{\sigma } &=&
\B{1}_{\rho }\times \oR _{\sigma }\times \B{1}_{\sigma \osigma } \,\circ \,
\B{1}_{\rho }\times \B{1}_{\iota } \times \oR _{\sigma } \,\,=\,\,
R_0^4 \, R_0^2.    \end{eqnarray*}
We get 
\begin{eqnarray*}
F(K\times L)\, W_{1,1}\, \ov{S}\ta \ov{T} &=&
d(\sigma )\, \ov{f\Bigl(\ovs{n}\Bigr) (K\times L)(Sf_0T)} \,\,= \\
d(\sigma ) \,\ov{f\Bigl(\ovs{n}\Bigr)(K)(S)\, f_0\, 
f\Bigl(\ovs{n}\Bigr)(L)(T) } &=&
W_{1,1}\, \ov{F(K)(S)} \ta \ov{F(L)(T)} \end{eqnarray*}
for $S,T \in B^{n+1}\subset L^2(B) \, (n\in \nat )$, and  (e) is established
for $m=l=s=r=1$.

\abs \noindent (b) If we replace $\sigma $ by $\sigma (m)$ for $m>0$ in (a),
we obtain
\begin{equation} f(\rho )(K\times L)(Sg_m T)\,=\,
 f(\rho )(K)(S)\, g_m \, f(\rho )(L)(T) \label{Eb1} \end{equation}
for $K,L \in ((\sigma \osigma )^m,\,(\sigma \osigma )^m ),\, 
S,T\in (\rho \sigma (m),
\rho \sigma (m))$ ($\rho =\ovs{n}$ )
and 
$$g_m:=\frac{1}{d(\sigma )^m} \B{1}_{\rho }\times 
     \left(\oR _{\sigma (m)}\circ \oR _{\sigma (m)}^*\right) 
\,\, \in \,\, \Bigl(\rho (\sigma \osigma )^m, \rho 
(\sigma \osigma )^m\Bigr)
     = B_{2m-1}^{n+1} \subset B_{2m-1}. $$
We get Property (e) from Step (2) for $m=l=s=r\in \nat $ as before
if we are able 
to prove
\begin{equation} W_{m,m} \, \ov{S} \ta \ov{T} \,=\, d(\sigma )^m \, \ov{Sg_mT}
 \qquad \hbox{for $S,T\in B_{m-1}$.} \label{Eb2} \end{equation}
First we will show
\begin{equation} g_m \,=\, d(\sigma )^{(m-1)m} f_{m-1,0}\cdot f_{m,1} \cdot
  \ldots \cdot f_{2m-2,m-1}. \label{Eb3} \end{equation}
We have
\begin{eqnarray}
f_{m+k-1,k} &=& \frac{1}{d(\sigma )^m} \,
  R_{m+k-1}^{2m-2}\, R_{m+k-1}^*\,R_{m+k-2}\, R_{m+k-2}^*\circ \ldots
  \circ R_{k+1}^* \,R_{k}\, R_{k}^{*2m}  \nonumber \\
&=&\frac{1}{d(\sigma )^m} \,
  R_{m+k-1}^{2m-2}\,R_k^{*2m}. \label{Eb4}  \end{eqnarray}
By induction on $k$ we will prove
\refstepcounter{equation} \label{Eb5}
$$\displaylines{
R_{m-1}^{2m-2} \circ \ldots \circ R_{m-k}\, R_{m-k-1}\, R_0^*\, R_1^* 
 \circ
\ldots \circ R_k^{*2m}\, =\,\hfill \cr
 \hfill R_{m-1}^{2m-2}\, R_0^*\,R_m\,R_1^*\,R_{m+1}\, R_2^*\, 
\circ \ldots \circ
R_{m+k-1}\, R_k^{*2m} \qquad (~\ref{Eb5}) }$$  
for $k= 0,\ldots ,m-1$. The case $k=m-1$ shows Equation (~\ref{Eb3}), as an
application of Equation (~\ref{Eb4}) yields. 

The case $k=0$ in Equation (~\ref{Eb5}) is obvious, the following 
computation shows the implication $k\rightarrow k+1,\, (k< m-1) $:
$$\displaylines{
R_{m-1}^{2m-2}\, R_0^* \circ \ldots \circ
R_{m+k-1}\, R_k^* R_{m+k}\, R_{k+1}^{*2m}\, = \hfill \hbox{(by induction
hypothesis)} \cr 
R_{m-1}^{2m-2} \circ \ldots \circ R_{m-k}\, R_{m-k-1}\, R_0^*\, R_1^* 
 \circ
\ldots \circ R_k^{*}\, R_{m+k}\, R_{k+1}^{*2m}\, = \hfill \cr
R_{m-1}^{2m-2} \circ \ldots \circ R_{m-k-1}\, R_{m-k-2}\, R_0^* 
 \circ
\ldots \circ R_k^{*}\, R_{k+1}^{*2m}.\hfill   }$$

It suffices to verify Equation (~\ref{Eb2}) for
$$T= y_1f_0y_2\ldots y_{m-1}f_{m-2,0}y_m, $$
where $y_1,\,y_2,\ldots ,y_m\in B$. The following computation deals with
this case:
$$\displaylines{
W_{m,m}\, \ov{S} \ta \ov{T}\,= \hfill \cr
d(\sigma )^{m^2} \, \ov{S f_{m-1,0}y_1 f_{m,0}y_2 \ldots y_{m-1}
  f_{2m-2,0} y_m} \, = \hfill \cr
d(\sigma )^{m^2} \, \ov{S f_{m-1,0}f_{m,1}y_1 f_0y_2 f_{m+1,0} y_3
 \ldots y_{m-1}f_{2m-2,0}y_m} \, = \hfill \cr
d(\sigma )^{m^2} \, \ov{S f_{m-1,0}f_{m,1}f_{m+1,2}y_1 f_0y_2 f_{1,0} y_3
 \ldots y_{m-1}f_{2m-2,0}y_m} \, = \hfill \cr
\ldots \, = \, \ldots \,=\, \hfill \cr
d(\sigma )^{m^2} \, \ov{S f_{m-1,0}f_{m,1}\, \ldots \,f_{2m-2,m-1}y_1 f_0y_2 
 f_{1,0} y_3
 \ldots y_{m-1}f_{m-2,0}y_m} \, = \hfill \cr
d(\sigma )^m \, \ov{Sg_mT}.\hfill  }$$

\abs \noindent (c) For $S\in B_{m-1}$ we prove
\begin{eqnarray} 
F_{0,1}(\oR _{\sigma })\,\ov{S}\!\! &= 
d(\sigma )^{\frac{1}{2}}\, \ov{S} \in L^2(B)  
 \hspace{18mm} &\mbox{for $m=0$,\ \ }\hspace{2mm} \label{Ec1}\\
F_{m,m+1}(\B{1}_{(\sigma \osigma )^m}\times \oR _{\sigma })\,\ov{S}\!\! &=
 d(\sigma )^{m +\frac{1}{2}} \ov{Sf_{m-1,0}} \in L^2(B_m)
  &\mbox{for $m\geq 1$\ \ }\hspace{2mm} \label{Ec2}\\
\lefteqn{\qquad \hbox{and}} \nonumber \\
F_{m,m+1}(\oR _{\sigma }\times \B{1}_{(\sigma \osigma )^m})\,\ov{S} \!\!&=
 d(\sigma )^{m +\frac{1}{2}} \ov{f_{0,m-1}S} \in L^2(B_m)
  &\mbox{for $m\geq 1$.\ \ } \hspace{2mm}\label{Ec3}   \end{eqnarray}
The first equation is a consequence of
$$f\Bigl(\ovs{n}\Bigr)(\oR _{\sigma })(S)\,=
\, d(\sigma )^{1/2}R_1^{*3}\, R_0^1\, S^1 \,=
  d(\sigma )^{1/2} S^1$$
for $S\in \Bigl(\ovs{n}, \ovs{n}\Bigr) =A^{n+1} \,( n\in \nat \cup \{0\})$.
Intending to derive Equation (~\ref{Ec2}) we state 
$$\displaylines{
f(\rho )(\B{1}_{(\sigma \osigma )^m}\times \oR _{\sigma })(S)\, =\hfill \cr
\hfill =\, d(\sigma )^{1/2} \,R_{m+1}^{*m+3}\circ \ldots \circ R_{2m}^{*3m+1}\,
R_{2m+1}^{*3m+3}\, R_{2m}^{3m+1}\, S^{3m+1}\, R_{m-1}^{3m-1}\circ
\ldots \circ R_0^{m+1}}$$
for $\rho = \ovs{n}$ and $S\in B_{m-1}^{n+1}$. Using 
$R_{2m+1}^{*3m+3}R_{2m}^{3m+1} \,=\B{1}_{\rho \sigma (3m+1)}$ and shifting 
$R_{2m}^{*3m+1}$ to the right,
we get
$$\displaylines{
f(\rho )(\B{1}_{(\sigma \osigma )^m}\times \oR _{\sigma })(S)\, = \hfill\cr
\hfill d(\sigma )^{1/2} \,R_{m+1}^{*m+3}\circ \ldots \circ R_{2m-1}^{*3m-1}\,
S^{3m+1}\, R_{m-1}^{3m-3}\circ
\ldots \circ R_1^{m+1}\, R_0^{m-1}\, R_0^{*m+1}.  }$$
Next we shift $R_{2m-1}^{*3m-1}$ to the right and so on such that
$$\displaylines{
f(\rho )(\B{1}_{(\sigma \osigma )^m}\times \oR _{\sigma })(S)\,=\hfill\cr
d(\sigma )^{1/2} \,R_{m+1}^{*m+3}\circ \ldots \circ R_{2m-2}^{*}\,
S\, R_{m-1}\circ
\ldots \circ R_2\, R_1\, R_1^*\, R_0\, R_0^{*m+1} \,=\hfill \cr
\ldots \,=\, \ldots \, = \hfill \cr
d(\sigma )^{1/2} \,S^{m+1}\, R_{m-1}\, R_{m-1}^*\circ \ldots \circ
R_1\, R_1^*\, R_0\, R_0^{*m+1} \,= \,d(\sigma )^{m +\frac{1}{2}} \,
Sf_{m-1,0}\hfill
          }$$
follows.
So Equation (~\ref{Ec2}) has been shown. 

We use the following computation in order to
prove Equation (~\ref{Ec3}):
$$\displaylines{
f(\rho )(\oR _{\sigma }\times \B{1}_{(\sigma \osigma )^m})(S)= 
d(\sigma )^{1/2} \,R_{m+1}^{*m+3}\circ \ldots \circ 
R_{2m+1}^*\, R_0\, S\, R_{m-1}\circ
\ldots \circ R_0^{m+1} \,=\hfill \cr
d(\sigma )^{1/2} \,R_0^{m-1}\,R_{m-1}^{*}\circ \ldots \circ 
R_{2m-2}^*\, R_{2m-1}^*\, S\, R_{m-1}\circ
\ldots \circ R_1\, R_0^{m+1} \,=\hfill \cr
d(\sigma )^{1/2} \,R_0^{m-1}\,R_{m-1}^{*}\circ \ldots \circ 
R_{2m-2}^*\, S\, R_{m-1}\circ
\ldots \circ R_1\, R_1^*\, R_0^{m+1} \,=\hfill \cr
d(\sigma )^{1/2} \,R_0^{m-1}\,R_{m-1}^{*}\circ \ldots \circ 
R_{2m-3}^*\, S\, R_{m-1}\circ
\ldots \circ R_2\, R_2^*\, R_1\, R_1^*\, R_0^{m+1} \,=\hfill \cr
\ldots \,= \, \ldots \,= \hfill\cr
d(\sigma )^{1/2} \,R_0^{m-1}\,R_{m-1}^{*}
S\, R_m^*\, R_{m-1}\circ
\ldots \circ R_2^*\, R_1\, R_1^*\, R_0^{m+1} \,=\hfill \cr
d(\sigma )^{1/2} \,R_0^{m-1}\,R_0^*\, R_1\, R_1^* \circ \ldots \circ
R_{m-1}\, R_{m-1}^{*}\, S^{3m+1}\, =
\qquad d(\sigma )^{m+\frac{1}{2}} \, f_{0,m-1} \,S.\hfill     }$$
In line 2 we shifted $R_0^{3m+1}$ to the left, in line 3 - 5 the computation
is similar to that of Equation (~\ref{Ec2}), in line 7 we used Equation 
(~\ref{ERc})
several times. 

\abs \noindent (d) Now we are able to verify Property (e) in Step (2) for any
arbitrary \linebreak $m,\, l,\, s,\, r\in \nat \cup \{0\}$. Let $M$ be a
natural number satisfying 
$M\geq \linebreak \max \{2,m,l,s,r\}$. For an integer $\nu >0$ let
$\oR _{\sigma }^{\times \nu }$ denote the $\nu $-fold product
$\oR _{\sigma }\times \ldots \times \oR _{\sigma }$ and let
$\oR _{\sigma }^{\times 0}:=\B{1}_{\iota }$. 
\begin{eqnarray*}
\tilde{K} &:=& (\oR _{\sigma }^{\times (M-l)} \times 
\B{1}_{(\sigma \osigma )^l}) \circ K \circ 
(\oR _{\sigma }^{\times (M-m)} \times 
\B{1}_{(\sigma \osigma )^m})^* \qquad \hbox{and} \\
\tilde{L} &:=& (\B{1}_{(\sigma \osigma )^r}\times
\oR _{\sigma }^{\times (M-r)}) 
 \circ L \circ 
(\B{1}_{(\sigma \osigma )^s}\times \oR _{\sigma }^{\times (M-s)} 
)^*   \end{eqnarray*}    
are operators of $((\sigma \osigma )^M, (\sigma \osigma )^M)$ such that the 
relation
\begin{equation}
F(\tilde{K} \times \tilde{L}) \, =\, 
W_{M,M} \circ (F(\tilde{K}) \ta F(\tilde{L})) \circ W_{M,M}^* 
    \label{Ed1} \end{equation}
is satisfied. Using $\oR _{\sigma }^* \circ \oR _{\sigma } =
d(\sigma )\B{1}_{\iota }$ we get
\begin{equation}
K = \frac{1}{d(\sigma )^{2M-m-l}}
\Bigl(\oR _{\sigma }^{\times (M-l)} \times 
\B{1}_{(\sigma \osigma )^l}\Bigr)^* \circ \tilde{K} \circ 
\Bigl(\oR _{\sigma }^{\times (M-m)} \times 
\B{1}_{(\sigma \osigma )^m}\Bigr)  \label{Ed2} \end{equation}
\qquad \hbox{and}
\begin{equation} L =\frac{1}{d(\sigma )^{2M-s-r}}
\Bigl(\B{1}_{(\sigma \osigma )^r}\times
\oR _{\sigma }^{\times (M-r)}\Bigr)^* 
 \circ \tilde{L} \circ 
\Bigl(\B{1}_{(\sigma \osigma )^s}\times \oR _{\sigma }^{\times (M-s)} 
\Bigr).   \label{Ed3}\end{equation}
We intend to show
\refstepcounter{equation} \label{Ed4}
$$\displaylines{
 F(\oR _{\sigma }^{\times (M-m)} \times 
\B{1}_{(\sigma \osigma )^{m+s}}) \,=\hfill\cr \hfill
W_{M,s}\circ \Bigl( F(\oR _{\sigma }^{\times (M-m)} \times 
\B{1}_{(\sigma \osigma )^m})\ta id_{L^2(B_{s-1})}\Bigr) \circ W_{m,s}^*. 
\qquad (~\ref{Ed4})   }$$
The cases $s=0$ and $M=m$ are obvious, the cases $s\geq 1$ and $0<m<M$ follow
from the following computation for $T\in B_{s-1}$ and \newl
$S = x_1 f_{0,m-2}x_2 \ldots x_{m-1} f_0 x_m 
           \in B_{m-1}\qquad (x_1,\ldots ,x_m \in B)$:
$$\displaylines{ 
F(\oR _{\sigma }^{\times (M-m)} \times 
\B{1}_{(\sigma \osigma )^{m+s}}) \circ W_{m,s}\,\, \ov{S}\ta \ov{T} \,= \hfill
\hbox{(see (c))}   \cr
d(\sigma )^{Ms +\frac{M^2 -m^2}{2}}
\ov{f_{0,M+s-2}\cdot \ldots \cdot f_{0,m+s-1}x_1 f_{0,m+s-2} x_2\cdot\ldots }
\hfill\cr
\hfill \ov{ \ldots \cdot x_{m-1}f_{0,s} x_m f_{0,s-1} T} \, =  \cr
d(\sigma )^{\frac{M^2-m^2}{2}} \, W_{M,s} \,
\ov{f_{0,M-2}\cdot \ldots \cdot f_{0,m-1}x_1 f_{0,m-2} x_2\cdot \ldots \cdot
x_{m-1}f_0 x_m} \ta \ov{T}= \hfill\cr \vspace{1mm} 
\hfill \hbox{(see (c))}  \cr
W_{M,s}\,\,  \Bigl(F(\oR _{\sigma }^{\times (M-m)} \times 
\B{1}_{(\sigma \osigma )^m})\, \, \ov{S}\Bigr) \, \ta \, \ov{T}.\hfill      }$$
The case $s\geq 1$ and $m=0 <M$ requires a separate computation
which we leave to the reader.
If we make use of Formula (~\ref{Ec2}) instead of Formula (~\ref{Ec3})
a computation similar as above yields
\refstepcounter{equation} \label{Ed5}
$$\displaylines{ 
F(\B{1}_{(\sigma \osigma )^{M+s}}\times \oR _{\sigma }^{\times (M-s)} 
) \, = \hfill\cr
\hfill W_{M,M}\,\circ \Bigl(id_{L^2(B_{M-1})} \ta 
F(\B{1}_{(\sigma \osigma )^s}\times 
\oR _{\sigma }^{\times (M-s)} 
)\Bigr) \circ W_{M,s}^*.\qquad (~\ref{Ed5}) }$$
By applying the Relations (~\ref{Ed2}), (~\ref{Ed3}), (~\ref{Ed4}),
(~\ref{Ed1})
and (~\ref{Ed5}) we get
$$\displaylines{ F(K\times L) \,= \hfill\cr
\frac{1}{d(\sigma )^{4M- (m+l+s+r)}}\, 
F(\oR _{\sigma }^{\times (M-l)} \times 
\B{1}_{(\sigma \osigma )^{l+r}})^* \circ
F(\B{1}_{(\sigma \osigma )^{M+r}}\times \oR _{\sigma }^{\times (M-r)} 
)^*
\circ \hfill \cr
F(\tilde{K}\times \tilde{L})\, \circ \, 
F(\B{1}_{(\sigma \osigma )^{M+s}}\times \oR _{\sigma }^{\times (M-s)})
\circ F(\oR _{\sigma }^{\times (M-m)} \times 
\B{1}_{(\sigma \osigma )^{m+s}})  \qquad = \hfill\cr
\frac{1}{d(\sigma )^{4M- (m+l+s+r)}}\, W_{l,r} \circ
\left(F(\oR _{\sigma }^{\times (M-l)} \times 
\B{1}_{(\sigma \osigma )^l})^* \ta id_{L^2(B_{r-1})} \right) \,\circ
\hfill\cr
\left(id_{L^2(B_{M-1})} \ta 
F(\B{1}_{(\sigma \osigma )^r}\times \oR _{\sigma }^{\times (M-r)})^* 
\right)\,\circ \hfill\cr
\left(
F(\tilde{K}) \ta F(\tilde{L}) \right) \circ \,
\left(id_{L^2(B_{M-1})} \ta  
F(\B{1}_{(\sigma \osigma )^s}\times \oR _{\sigma }^{\times (M-s)}) \right)\, 
\circ \hfill\cr
\left( F(\oR _{\sigma }^{\times (M-m)} \times 
\B{1}_{(\sigma \osigma )^m}) \ta id_{L^2(B_{s-1})} \right) 
\circ W_{m,s}^* \qquad = \hfill \cr
W_{l,r}\circ \left( F(K) \ta F(L) \right)\circ  W_{m,s}^*. \hfill      }$$

\abs \noindent (e) It remains to show that the linear maps
$$F_{m,l}:((\sigma \osigma )^m, (\sigma \osigma )^l) \pfeil
\C{L} _{A,A} (L^2(B_{m-1}),L^2(B_{l-1}))$$
 are surjective for
$m,\,l\in \nat \cup \{0\}$. It is well known that there is a canonical 
linear isomorphism from $A'\cap B_{2m-1}$ onto 
$\C{L} _{A,A} (L^2(B_{m-1}),L^2(B_{m-1}))$ for $m\in \nat $ (see
\cite{Sch3} for example). Furthermore, Theorem ~\ref{TUnt} (iv)
shows \linebreak
$A'\cap B_{2m-1}\cong ((\sigma \osigma )^m, (\sigma \osigma )^m)$
for $m\in \nat \cup \{0\}$, and $F_{m,m}$ is surjective.

We regard the case $l<m$: 
consider an operator 
\newline $X\in \C{L} _{A,A} (L^2(B_{m-1}),L^2(B_{l-1}))$.  
Since $F_{m,m}$ is surjective, there is
an operator $K\in ((\sigma \osigma )^m, \, (\sigma \osigma )^m)$  such that
$$F_{m,m}(K) = F_{l,m}(\B{1}_{(\sigma \osigma )^l}\times 
\oR _{\sigma }^{\times (m-l)}) \circ X. $$
Using $\oR _{\sigma }^* \circ \oR _{\sigma } =d(\sigma )\B{1}_{\iota }$, 
we get
$$\displaylines{
F_{m,l} \Bigl(\frac{1}{d(\sigma )^{m-l}} 
\Bigl(\B{1}_{(\sigma \osigma )^l}\times 
\oR _{\sigma }^{*\times (m-l)}\Bigr)
\circ K\Bigr)=  \hfill \cr
F_{l,l} \Bigl(\frac{1}{d(\sigma )^{m-l}} \Bigl(\B{1}_{(\sigma \osigma )^l}
\times 
(\oR _{\sigma }^{*\times (m-l)} \circ \oR _{\sigma }^{\times (m-l)}) \Bigr)
\Bigr)\,\circ X \,= \, \,X. \hfill   }$$

At last we reduce the case $l>m$ to the case $l<m$ by applying the
$*$-operation. \blacksquare

\subsection{Finite dimensional Hopf-$*$-algebras as examples}
We consider the case that $\C{C}$ is the $C^*$-tensor category $\C{U}_H$,
where $H$ is a finite dimensional Hopf-$*$-algebra. For the sake of
simplicity we
assume $\C{C}_{\sigma \osigma }\, =\C{C}$. That means that there
is an $n\in \nat $ such that every irreducible corepresentation of $H$ is
contained in $(\sigma \osigma )^n$. 

Theorem ~\ref{TCB} tells us that $\C{B}_{A\subset B}$ is equivalent to
$\C{U}_H$.

On the other hand, the considerations in Section 3 imply that 
$A\subset B$ is equal to an 
inclusion $N^K\subset (N\otimes \B{L}(V))^K$ of fixed point algebras,
where the action $\alpha $ of
$K:=(H^o)^{cop}$ on $N$ is outer. There is an outer action of the dual 
Hopf-$*$-algebra
$K^o$ on $M:=N^K$ such that the subfactor
$N^K\subset N$ is isomorphic to $M\subset M\kreuz K^o$.
One easily sees
that $\C{B}_{A\subset B}= \C{B}_{N^K\subset (N\otimes \B{L}(V))^K}$
is a full $C^*$-tensor subcategory of 
$\C{B}_{N^K\subset N\otimes \B{L}(V)}\, =\C{B}_
{N^K\subset N}$. Hence Proposition ~\ref{PH1} yields that
$\C{B}_{N^K\subset N}$ is equivalent to $\C{U}_{(K^o)^{cop}}$.
Clearly, 
$$(K^o)^{cop}= ((H^o)^{cop})^{o\, cop}\, =
((H^{o\,o })^{op})^{cop}\,=H^{op\,
cop}$$ 
holds. The Hopf-$*$-algebras
$H$ and $H^{op\, cop }$ are isomorphic, as an application of the 
antipode $S$ of
$H$ shows. So we conclude that $\C{B}_{A\subset B}$ is equivalent to
a full $C^*$-tensor subcategory of $\C{U}_H$, and
Theorem ~\ref{TCB} has been confirmed for this special case.
 
\appendix
\section{Appendix: Proof of Proposition 1.5 (ii)}
We will proceed similarly as for the equivalences 
of arbitrary
categories (for example see \cite{Ka}, Proposition XI.1.5), but more details
have to be verified.

Let an equivalence $(F,(U_{\rho \sigma })_{\rho ,\sigma },J)$ of the
$C^*$-tensor categories $\C{C}$ and $\C{D}$ be given. For each object $\phi
\in \Ob \C{D}$ we choose an object $G(\phi )$ of $\C{C}$ such that 
$F(G(\phi ))$ is equivalent to $\phi $ and a unitary operator 
$W_{\phi }\in (\phi ,F(G(\phi )))$. Especially we put
\begin{equation} G(\iota _{\C{D}}):=\iota _{\C{C}} \qquad \hbox{and} \qquad
     W_{\iota _{\C{D}}}:=J^*\in (\iota _{\C{D}}, \, F(\iota _{\C{C}})).
     \label{ECS1} \end{equation}     
The linear map $F_{\rho ,\sigma }:T\in (\rho ,\sigma )\Pfeil F(T)\in 
 (F(\rho ), F(\sigma ))$ is invertible 
for all objects $\rho $ and $\sigma $. For $S\in (\phi ,\psi )$ ($\phi ,\,
\psi \in \Ob \C{D}$) we define 
\begin{equation} G(S):= F_{G(\phi ), G(\psi )}^{-1}(W_{\psi }\circ 
S\circ W_{\phi }^*)\in
(G(\phi ), G(\psi )). \label{ECS2} \end{equation}
One easily checks that a full and faithful $C^*$-functor 
$G:\C{D} \pfeil \C{C}$ is defined
in this way.

We note that Definition (~\ref{ECS2}) implies
\begin{equation} F(G(S))\circ W_{\phi } = W_{\psi } \circ S \qquad 
\hbox{for $S\in (\phi ,
     \psi )$.}  \label{ECS3} \end{equation}     
In order to define a $C^*$-tensor equivalence we introduce unitary operators
$V_{\phi \psi }:G(\phi )\, G(\psi ) \pfeil G(\phi \psi )$ by putting
$$V_{\phi \psi }:= F_{G(\phi )G(\psi ), G(\phi \psi )}^{-1} (W_{\phi \psi}
\circ W_{\phi }^* \times  W_{\psi }^* \circ U_{G(\phi ), G(\psi )}^*)$$
for all objects $\phi $ and $\psi $ of $\C{D}$. We will prove that
$(G, (V_{\phi \psi })_{\phi ,\psi }, \B{1}_{\iota _{\C{C}}})$ is a 
$C^*$-tensor equivalence with $[G] \, =[F]^{-1}$. For this purpose
we have to verify (a), (b) 
and (c) in Definition ~\ref{DCequ} (ii), the remaining assumptions are 
obviously satisfied.

\abs
\noindent (a) We consider the following diagram for objects $\phi ,\, \phi ',\,
\psi $ and $\psi '$ of $\C{D}$ and morphisms $R\in (\phi ,\phi ')$ and
$S\in (\psi ,\psi ')$:
\begin{center}
\xext=1800 \yext=1800
\adjust[`F(G(R)\times G(S));F(G(\phi )F(G(\psi ))`;F(G(\phi ')F(G(\psi '))`;
`F(G(R\times S))]
\begin{picture}(\xext,\yext)(\xoff,\yoff)
\setsqparms[1`1`1`1;1800`600]
\putsquare(0,1200)[F(G(\phi )G(\psi ))`F(G(\phi ')G(\psi '))`
F(G(\phi )) F(G(\psi ))`F(G(\phi ')) F(G(\psi '));
F(G(R)\times G(S))`U^*_{G(\phi ),G(\psi )}`U^*_{G(\phi '),G(\psi ')}`
F(G(R)) \times F(G(S))]
\setsqparms[1`1`1`1;1800`600]
\putsquare(0,0)[\phi \psi `\phi '\psi '`F(G(\phi \psi ))`F(G(\phi '\psi '));
R\times S`W_{\phi \psi }`W_{\phi '\psi '}`F(G(R\times S))]
\putmorphism(0,1200)(0,-1)[\phantom{F(G(\phi ))F(G(\psi ))}`
   \phantom{\phi \psi }`W_{\phi }^*\times W_{\psi }^*]{600}{1}{l}
\putmorphism(1800,1200)(0,-1)[\phantom{F(G(\phi '))F(G(\psi '))}`
   \phantom{\phi '\psi '}`W_{\phi '}^*\times W_{\psi '}^*]{600}{1}{r}
\put(900,1500){\makebox(0,0){(1)}}
\put(900,900){\makebox(0,0){(2)}}
\put(900,300){\makebox(0,0){(3)}}
\end{picture} \end{center}
Definition ~\ref{DCequ} (ii) (a) implies that diagram (1) commutes.
Furthermore the diagrams (2) and (3)
commute according to Relation (~\ref{ECS3}). Hence the exterior diagram is 
commuting, and the application of $F^{-1}$ shows
$$G(R\times S)\circ V_{\phi \psi } \, = 
  V_{\phi '\psi '}\circ G(R) \times G(S).$$

\abs
\noindent
(b) We have to check 
\refstepcounter{equation} \label{ECS4}
$$\displaylines{ G(a(\phi ,\psi ,\eta )) \circ V_{\phi ,\psi \eta } \circ 
     \B{1}_{G(\phi )} \times V_{\psi \eta } \, = \hfill\cr
  \hfill   V_{\phi \psi ,\eta } \circ
     V_{\phi \psi }\times  \B{1}_{G(\eta )}\circ 
     a(G(\phi ),G(\psi ), G(\eta )) \qquad (~\ref{ECS4}) }$$
for all objects $\phi ,\,\psi ,\eta \in \Ob \C{D}$.   
We know
\begin{equation} F(\B{1}_{G(\phi )}\times V_{\psi \eta })\, =
   U_{G(\phi ),G(\psi \eta )}\circ \B{1}_{F(G(\phi ))}\times F(V_{\psi \eta })
   \circ U_{G(\phi ),G(\psi )G(\eta )}^*. \label{ECS5} \end{equation}
By applying $F$ on both sides of Equation (~\ref{ECS4}) and using Equation
(~\ref{ECS5}) as well as the corresponding relation for 
$F(V_{\phi \psi }\times \B{1}_{G(\eta )})$, we find that (~\ref{ECS4}) is 
equivalent to the equation 
\refstepcounter{equation} \label{ECS6} 
$$\displaylines{
F(G(a(\phi ,\psi ,\eta )))\circ W_{\phi (\psi \eta )}\circ W_{\phi }^*\times
W_{\psi \eta }^* \circ U_{G(\phi ),G(\psi \eta )}^* \circ 
U_{G(\phi ),G(\psi \eta )}\circ \B{1}_{F(G(\phi ))}\times W_{\psi \eta }
\circ \hfill \cr
 \circ \B{1}_{F(G(\phi ))}\times (W_{\psi }^* \times W_{\eta }^*) \circ
\B{1}_{F(G(\phi ))}\times U_{G(\psi ),G(\eta )}^* \circ 
U_{G(\phi ),G(\psi )G(\eta )}^* \, = \hfill \cr
 W_{(\phi \psi )\eta }\circ W_{\phi \psi }^*\times
W_{\eta }^* \circ U_{G(\phi \psi ),G(\eta )}^* \circ 
U_{G(\phi \psi ),G(\eta )} \circ W_{\phi \psi }\times \B{1}_{F(G(\eta ))}
\circ \hfill \cr
\circ (W_{\phi }^* \times W_{\psi }^*)\times \B{1}_{F(G(\eta ))} \circ
U_{G(\phi ),G(\psi )}^*\times  \B{1}_{F(G(\eta ))} \circ
U_{G(\phi )G(\psi ),G(\eta )}^*\circ \hfill \cr
\circ F(a(G(\phi ),G(\psi ),G(\eta ))).
\hfill (~\ref{ECS6})           }$$
(~\ref{ECS3}) implies
$$ F(G(a(\phi ,\psi ,\eta )))\, = W_{(\phi \psi )\eta }\circ 
  a(\phi ,\psi ,\eta )\circ W_{\phi (\psi \eta )}^*, $$
moreover
$$\displaylines{
 F(a(G(\phi ),G(\psi ),G(\eta ))) \circ U_{G(\phi ),G(\psi )G(\eta )}\circ
\B{1}_{F(G(\phi ))}\times U_{G(\psi ),G(\eta )}\, =\hfill \cr
\hfill U_{G(\phi )G(\psi ),G(\eta )}\circ
U_{G(\phi ),G(\psi )} \times  \B{1}_{F(G(\eta ))} \circ
a(F(G(\phi )),F(G(\psi )), F(G(\eta ))).    }$$
Inserting these relations into (~\ref{ECS6}) we see that Equation 
(~\ref{ECS4}) is equivalent to
$$\displaylines{
a(\phi ,\psi ,\eta )\, \circ \, W_{\phi }^*\times W_{\psi \eta }^* \, \circ \,
\B{1}_{F(G(\phi ))}\times W_{\psi \eta } \, \circ \,
\B{1}_{F(G(\phi ))}\times (W_{\psi }^*\times W_{\eta }^*)\, = \hfill \cr
W_{\phi \psi }^* \times W_{\eta }^* \, \circ \,
W_{\phi \psi }\times \B{1}_{F(G(\eta ))} \, \circ \,
(W_{\phi }^* \times W_{\psi }^*)\times \B{1}_{F(G(\eta ))} \, \circ \hfill \cr
\hfill \circ \, a(F(G(\phi )),F(G(\psi )), F(G(\eta )))    }$$  
and to
$$\displaylines{
\qquad 
a(\phi ,\psi ,\eta )\, \circ \, W_{\phi }^*\times 
(W_{\psi }^*\times W_{\eta }^*)\, = \hfill \cr
\hfill (W_{\phi }^* \times W_{\psi }^*)\times W_{\eta }^*
\, \circ \, a(F(G(\phi )),F(G(\psi )), F(G(\eta ))). }$$ 
The last relation is satisfied, because $a(\phi ,\psi ,\eta )$ is natural in
$\phi ,\, \psi $ and $\eta $.

\abs 

\noindent (c) By Definition (~\ref{ECS1}) we have to establish the equation
\begin{equation}
l_{G(\phi )}\, = G(l_{\phi })\circ V_{\iota _{\C{D}},\phi } \label{AE1}
\end{equation}
for every object $\phi \in \Ob \C{D}$. Equation (~\ref{AE1}) is equivalent to
\begin{equation} F(l_{G(\phi )})\, = F(G(l_{\phi }))\circ 
W_{\iota _{\C{D}}\phi }
\circ J\times W_{\phi }^* \circ U_{\iota _{\C{C}}, G(\phi )}^*. 
\label{ECS7} \end{equation}
From Definition ~\ref{DCequ} (ii) (c) we conclude
$$F(l_{G(\phi )})\circ U_{\iota _{\C{C}}, G(\phi )}\, = l_{F(G(\phi ))}\circ
  J\times \B{1}_{F(G(\phi ))}.  $$
By inserting this relation into (~\ref{ECS7}) we obtain that Equation 
(~\ref{ECS7}) is equivalent to 
\begin{equation} l_{F(G(\phi ))}\circ
  J\times \B{1}_{F(G(\phi ))}=F(G(l_{\phi }))\circ W_{\iota _{\C{D}}\phi }
\circ J\times W_{\phi }^*. \label{ECS8} \end{equation}
Since $l_{\phi }$ is natural in $\phi $,
$$ l_{F(G(\phi ))}\circ \B{1}_{\iota _{\C{D}}}\times W_{\phi } \, = 
   W_{\phi }\circ l_{\phi } $$
holds, and Equation (~\ref{ECS3}) yields
$$ W_{\phi }\circ l_{\phi } \, = 
F(G(l_{\phi }))\circ W_{\iota _{\C{D}}\phi }.$$
By applying those two equations we see that Equation (~\ref{ECS8})
is satisfied. \blacksquare

\end{document}